\documentclass{amsart}[draft]

\usepackage{amsmath, amsfonts, amsthm, amssymb}
\usepackage{array}
\usepackage{lmodern}
\usepackage{soul}
\usepackage{enumerate}
\usepackage{mathtools}
\usepackage{csquotes}
\usepackage{cite}
\usepackage{tikz}
\usepackage{tikz-cd}
\usetikzlibrary{matrix,arrows,decorations.pathmorphing}
\usepackage{hyperref}
\usepackage[final]{pdfpages}
\usepackage{wrapfig}
\usepackage{bm}

\usepackage{cancel}

\renewcommand{\Re}{\operatorname{Re}}
\renewcommand{\Im}{\operatorname{Im}}

\newcommand{\R}{\ensuremath\mathbb{R}}
\newcommand{\E}{\ensuremath\mathcal{E}}
\let\L\relax
\let\C\relax
\newcommand{\A}{\ensuremath\mathcal{A}}
\newcommand{\C}{\ensuremath\mathbb{C}}
\newcommand{\L}{\ensuremath\mathcal{L}}
\newcommand{\Z}{\ensuremath\mathbb{Z}}
\newcommand{\Q}{\ensuremath\mathbb{Q}}
\newcommand{\N}{\ensuremath\mathbb{N}}
\renewcommand{\H}{\ensuremath\mathbb{H}}
\newcommand{\Hb}{\ensuremath\mathbb{H}}

\newcommand{\cD}{\ensuremath\mathcal{D}}
\newcommand{\cV}{\ensuremath\mathcal{V}}
\newcommand{\cN}{\ensuremath\mathcal{N}}
\newcommand{\cG}{\ensuremath\mathcal{G}}
\newcommand{\1}{\mathbf 1}

\newcommand{\eps}{\ensuremath\varepsilon}

\newcommand{\bt}{\bm t}

\renewcommand{\mod }{\text{ {\rm mod} }}
\newcommand{\modulo}{\text{ \rm mod }}
\renewcommand{\setminus}{\smallsetminus}
\newcommand{\abs}[1]{\left|#1\right|}
\newcommand{\floor}[1]{\lfloor #1\rfloor}
\newcommand{\PP}{\ensuremath\mathbb{P}}
\newcommand{\EE}{\ensuremath\mathbb{E}}
\newcommand{\VV}{\ensuremath\mathbb{V}}

\newcommand{\df}{\mathop{}\!\mathrm{d}}

\DeclareMathOperator{\sgn}{sgn}
\DeclareMathOperator*{\Res}{Res}
\DeclareMathOperator{\e}{e}
\DeclareMathOperator{\id}{id}
\DeclareMathOperator{\SL}{SL}
\DeclareMathOperator{\GL}{GL}
\DeclareMathOperator{\PSL}{PSL}
\DeclareMathOperator{\sym}{sym}

\DeclareMathOperator{\den}{den}

\newtheorem{thm}{Theorem}[section]

\newtheorem{cor}[thm]{Corollary}
\newtheorem{conj}[thm]{Conjecture}
\newtheorem{lemma}[thm]{Lemma}
\newtheorem{prop}[thm]{Proposition}
\theoremstyle{remark}
\newtheorem{rem}[thm]{Remark}

\newcommand\numberthis{\stepcounter{equation}\tag{\theequation}}

\numberwithin{equation}{section}

\title{Central values of additive twists of Maa{\ss} form $L$-series}

\author{Sary Drappeau}
\address{I2M, Université d'Aix-Marseille, CNRS \\ Campus de Luminy, 13288 Marseille Cedex 9, France}
\curraddr{LMBP, Institut Universitaire de France, Université Clermont-Auvergne \\ 3 place Vasarely, 63178 Aubière Cedex, France}
\email{sary.drappeau@uca.fr}

\author{Asbjørn Christian Nordentoft}
\address{University of Copenhagen, Universitetsparken 5, 2100 Copenhagen Ø, Denmark}
\email{nordentoft@math.ku.dk}

\date{\today}



\thanks{This work was supported by CNRS through IRN GANDA (Geometry and Arithmetic) and IRN MaDeF (Mathematics in Denmark and France). The authors thank CNRS for financial support, and the Universities of Witswatersrand and Aix-Marseille for their hospitality. The authors are indebted to Fabien Pazuki and Adeline Leclercq Samson for making this support possible, and to Morten Risager for remarks on a earlier version on the manuscript. The second named author's research was supported by the Independent Research Fund Denmark DFF-1025-00020B. The first named author's research was supported by the joint FWF-ANR project Arithrand: FWF: I 4945-N and ANR-20-CE91-0006.}

\begin{document}

\begin{abstract}
  In the present paper we study the central values of additive twists of Maa{\ss} form $L$-series. In the case of the modular group, we show that the additive twists (when averaged over denominators) are asymptotically normally distributed. This supplements the recent work of Petridis--Risager which settled an averaged version of a conjecture of Mazur--Rubin concerning modular symbols. The methods of the present paper combine dynamical input due to Bettin and the first named author with the new fact that the additive twists define quantum modular forms in the sense of Zagier. This latter property is shown for a general discrete, co-finite group with cusps. Our results also has a number of arithmetic applications; in the case of Hecke congruence groups the quantum modularity implies certain reciprocity relations for twisted moments of twisted $\GL_2$-automorphic $L$-functions, extending results of Conrey and the second named author. In the case of cuspidal Maa{\ss} forms for the modular group, we also obtain a calculation of certain wide moments of twists of the $L$-function of the Maa{\ss} form.
\end{abstract}

\maketitle

\section{Introduction}

\subsection{Central values}

Let~$\Gamma \subset \SL(2, \R)$ be a discrete, co-finite group with cusps, and~$\phi:\H \to \C$ be a Maa{\ss} form for~$\Gamma$.
To~$\phi$ we associate the sequence~$(a(n))_{n\neq 0}$ of its Fourier coefficients, by means of which we form the~$L$-series
$$ L(\phi, s) := \sum_{n=1}^\infty \frac{a(n)}{n^{s-1/2}}, $$
initially convergent on same half-plane, and analytically continued.
The normalization here is such that when~$\Gamma = \Gamma_0(N)$ is a Hecke congruence subgroup, the Ramanujan-Peterson conjecture~\cite[p.~95]{IwaniecKowalski2004} predicts~$\abs{a(n)} \leq C_\eps n^{-1/2+\eps}$ for~$\phi$ cuspidal, all~$\eps>0$ and some~$C_\eps>0$.
When~$\Gamma$ is a Hecke congruence subgroup, and~$\phi$ is a Hecke newform, then~$L(s, \phi)$ satisfies additionally a functional equation relating~$s$ to~$1-s$~\cite{DukeEtAl2002}, and has an Euler product factorization~\cite[chapter~5.11]{IwaniecKowalski2004}. It is an instance of a rank~$2$ $L$-function in the Selberg class, and it is conjectured that all~$L$ functions in the Selberg class of rank~$2$ over $\Q$ can be obtained in this way~\cite{Perelli}. This statement for Artin $L$-functions of degree~$2$ representations is a particular case of the Langlands conjectures.

To simplify the exposition we assume here that~$\Gamma = \Gamma_0(N)$ is a Hecke congruence subgroup. The present work concerns the analytic properties of the twisted~$L$-value
$$ L(\phi, s, x) := \sum_{n=1}^\infty \frac{a(n)}{n^{s-1/2}} \e^{2\pi i n x} $$
as a function of~$x\in\R$. The sum converges uniformly on compacts inside~$\{\Re(s)>1\}$.
At the edge $\Re(s)=1$ of this half-plane, the regularity with respect to~$x$ of~$L(\phi, s, x)$ is an old theme, see for instance~\cite{Wilton1933}, and the references in~\cite{BalazardMartin2019}.

When~$x\in\Q$, the map~$s\mapsto L(\phi, s, x)$ has an analytic continuation to~$\C$ minus a possible simple pole at~$s=1$. This meromorphic continuation has a functional equation relating~$(s, a/q)$ to~$(1-s, -\bar{a}/q)$ for~$N\mid q$ and~$a\bar{a} \equiv 1\pmod{q}$~\cite[eqs.~(A.10)-(A.13)]{KoMiVa02} (a more complicated relation holds in cases when~$N\nmid q$).
We note that one important special case is~$s = \tfrac12\pm s_\phi$, where~$s_\phi$ is the Laplace eigenvalue associated with~$\phi$. It was shown by Lewis and Zagier that the map~$x\mapsto L(\phi, \tfrac12+s_\phi, x)$ extends to~$\R$ and enjoys distinctive analytic properties, which generalize in a way those of the Eichler-Shimura map, see~\cite{LewisZagier2001, Bruggeman2007, BruggemanEtAl2015}.
In general, for~$\Re(s)>1/2$, it seems possible to extend~$L(\phi, s, x)$ to a continuous function of~$x\in\R$ outside a set of zero Lebesgue measure by arguments similar to the case~$\Re(s)=1$ (see~\cite{BettinDAppendixA}), but this breaks down at~$\Re(s)=1/2$.

We are interested in the central values
$$ L_\phi(x) := L(\phi, 1/2, x) \qquad (x\in \Q). $$
In the situation when~$\phi$ is associated to a holomorphic cusp form~$f$ of weight~$k$, meaning that~$\phi(z) = y^{k/2} f(z)$, then there is a constant~$c_f$ such that
\begin{equation}\label{eq:Lphi-Eichler-intro}
  L_\phi(x) = c_f \int_x^\infty f(z) (z-x)^{k/2-1} \df z.
\end{equation}
In particular, for~$k=2$, the value~$L_\phi(x)$ is essentially the modular symbol~$\langle x\rangle_f$ associated to~$f$. On the other hand, by orthogonality relations, the central $L$-value of the multiplicative twist
\begin{equation}\label{eq:holomorphicadditive} L(f\otimes \chi, s) = \sum_{n=1}^\infty \frac{a(n)}{n^{s-1/2}} \chi(n) \end{equation}
by a primitive Dirichlet character~$\chi\pmod{q}$ can be expressed in terms of a weighted average of~$\langle a/q\rangle_f$ with~$a$ varying over classes mod~$q$ (see Proposition~\ref{BS}).
In turn, when~$f$ is of weight~$2$ and is associated to an elliptic curve, then it is expected that the value~$L(f\otimes \chi, 1/2)$ encodes geometric information on the underlying elliptic curve, in particular its rank over abelian number fields, see~\cite[Proposition~2.2]{MazurRubin}. Motivated by conjectures about elliptic curves, Mazur and Rubin~\cite{MazurRubin} were led to conjecture, among other phenomena, that the multisets
\begin{equation}\label{eq:MRconj-qfixed}
  \cD_\phi(q) := \Big\{ \frac{L_\phi(a/q)}{(\log q)^{1/2}} : a\in (\Z/q\Z)^\times \Big\}
\end{equation}
become distributed as~$q\to\infty$ according to a centered normal law, when~$\phi$ is a weight~$2$ form associated to an elliptic curve~\cite[Conjecture~4.3]{MazurRubin}.
We believe that this holds in general for any Maa\ss{} form.
\begin{conj}[Additive twists conjecture]\label{conj:addtwist}
  Let $\phi$ be a Hecke-Maa{\ss} cusp form for~$\Gamma_0(N)$. Then the multi-sets~$\cD_\phi(q)$ become asymptotically normally distributed as $q\rightarrow \infty$.
\end{conj}

Note that this would be in contrast to the expected log-normal behaviour of central values of families of $L$-functions having an Euler product as in Selberg's Central Limit Theorem, see for instance~\cite{Selberg1946, KeatingSnaith2000, RadziwillSound2015}.

Regarding this conjecture, the furthest achievement so far is a power-saving estimate for the second moment~\cite{BlFoKoMiMiSa18}, which lies at the edge of current known techniques in analytic number theory. On average over~$q$, however, the corresponding statement is now known for forms~$\phi$ which are associated with a holomorphic form, by Petridis and Risager~\cite{PetridisRisager2018} and the second named author~\cite{Nordentoft2021}.

\begin{thm}[Additive twists on average for holomorphic forms, \cite{Nordentoft2021}]\label{thm:addtwist-cusp}
  Let $f$ be a holomorphic Hecke cusp form of integer weight~$k$ for~$\Gamma$ discrete co-finite with cusps, and~$\phi(z) = (\Im z)^{k/2} f(z)$. Then the multi-sets
  $$ \left\{ \frac{L_\phi(a/q)}{(\log q)^{1/2}} : a\in (\Z/q\Z)^\times, q\leq Q \right\}, $$
  become asymptotically normally distributed as $q\rightarrow \infty$.
\end{thm}
The proofs use spectral theory of automorphic forms relying on certain twisted Eisenstein series introduced by Goldfeld \cite{Go99}, \cite{Go99.2} whose analytic properties (poles, growth on vertical lines, \emph{etc.}) are studied using analytic properties of automorphic resolvent operators (a technique pioneered by Petridis-Risager \cite{PeRi2}). The holomorphicity of $f$ plays a crucial role in the key \emph{automorphic completion}-step, which enables one to express the non-automorphic Eisenstein series of Goldfeld in terms of (a finite sum of) automorphic Poincar\'{e} series (we refer to \cite[Section 2]{Nordentoft2021} for more details on the method of proofs). More precisely, one does a contour shift in the integral representation~\eqref{eq:Lphi-Eichler-intro} using crucially that the integrand is holomorphic. In the case of non-holomorphic  Maa{\ss} forms the corresponding approach breaksdown which makes Conjecture \ref{conj:addtwist} a particularly interesting challenge.  We note that one advantage of the automorphic approach in~\cite{PetridisRisager2018} is that it applies to general Fuchsian groups of the first kind, basically without extra effort. 

In the two special cases~$k=2$ and~$N=1$, another proof of Theorem~\ref{thm:addtwist-cusp} was obtained respectively by Lee-Sun~\cite{LeeSun2019} and by Bettin and the first named author~\cite{BettinDrappeau2019} around the same time as ~\cite{PetridisRisager2018} and ~\cite{Nordentoft2021}.
Both proofs ultimately rely on the ``quantum modularity'' property of~$L_\phi$ in these two cases, in the sense of Zagier~\cite{Zagier2010}, by which we mean that for any~$\gamma \in \Gamma_0(N)$, the map
\begin{equation}\label{eq:def-hgamma-intro}
  h_\gamma : x\mapsto L_\phi(\gamma x) - \chi(\gamma) L_\phi(x)
\end{equation}
extends to a map on~$\R\setminus\{\gamma^{-1}\infty\}$ which is constant when~$k=2$~\cite[Section~1.3]{LeeSun2019}, and Hölder-continuous with a uniform exponent for~$N=1$~\cite[Lemma~9.3]{BettinDrappeau2019}.
It was also shown in~\cite{Bettin16} that this property holds when~$\phi$ is a certain non-holomorphic Eisenstein series of level~$1$, for which~$L_\phi(x)$ is the Estermann function. Combined with \cite[Theorem 3.1]{BettinDrappeau2019}, this leads to a normal distribution result for the Estermann function (which definitely lies beyond the scope of the automorphic approach since one cannot apply the method of moments). More generally, the dynamical framework in \cite{BettinDrappeau2019} applies to general quantum modular forms satisfying certain growth conditions. This raises the question of the regularity and growth of~$h_\gamma$ in general.

\subsection{Statement of results: regularity and growth of \texorpdfstring{$h_\gamma$}{h-gamma} in general}

In the present paper we study the analytic properties of the maps~$h_\gamma$ in the general context of a Maa\ss{} form (not necessarily cuspidal) for a cofinite Fuchsian group.

Let~$\Gamma\subset \SL(2, \R)$ be a cofinite Fuchsian group with cusps, and let~$\phi$ be a Maa\ss{} form of weight~$k\in\Z_{\geq 0}$ and multiplier~$\chi$ for~$\Gamma$ (not necessarily cuspidal). The precise definition will be presented in Section~\ref{sec:setting} below. Let~$C(\Gamma)\subset\R\cup\{\infty\}$ be the set of cusps of~$\Gamma$.
For any~$\gamma \in \Gamma$, we consider the map~$h_\gamma : C(\Gamma)\setminus\{\infty, \gamma^{-1}\infty\} \to \C$ given by
$$ h_\gamma(x) := L_\phi(\gamma x) - \chi(\gamma) \sgn(x-\gamma^{-1}\infty)^k L_\phi(x). $$

\begin{thm}[Regularity of~$h_\gamma$]\label{thm:intro-reg}
  The map~$h_\gamma$ extends to a~$(1/2-\eps)$-Hölder continuous function of~$x$ in~$\R\setminus\{\gamma^{-1}\infty\}$.
\end{thm}

\begin{thm}[Growth of~$h_\gamma$]\label{thm:intro-growth}
  If~$\lambda_\phi$ denotes the eigenvalue of~$\phi$ and~$s_\phi(1-s_\phi) = \lambda_\phi$, then there exist constants~$A_+, A_-, B_+, B_-$ and~$C \in \C$ such that, as~$\abs{x}\to\infty$ with~$\sgn(x) = \pm 1$, there holds
  $$ h_\gamma(x) = \chi(\gamma)(A_\pm \abs{x-\gamma^{-1}\infty}^{s_\phi} + B_\pm \abs{x-\gamma^{-1}\infty}^{1-s_\phi} + C) + O_{\phi, \eps, \gamma}(\abs{x}^{-1+\eps}), \qquad (\lambda_\phi \neq 1/4). $$
  If~$\lambda_\phi = 1/4$, this estimate holds with the term involving~$B_\pm$ replaced by~$B_\pm \abs{x-\gamma^{-1}\infty}^{1/2}\log\abs{x-\gamma^{-1}\infty}$.
  We have~$A_\pm = B_\pm = 0$ if~$\phi$ is cuspidal.
\end{thm}

Before the present work, Theorems~\ref{thm:intro-reg} and~\ref{thm:intro-growth} were known in essentially two cases:
\begin{itemize}
\item For forms of weight~$2$, the map~$h_\gamma$ is constant. This is a well-known property of modular symbol, see~\cite[Lemma~1.2.(iii)]{MazurRubin}.
\item When~$\phi$ is associated to a holomorphic form of even weight, this was proved in~\cite[Lemma~9.3]{BettinDrappeau2019} in the special case~$\Gamma=\SL(2, \Z)$ and in~\cite[Theorem~4.4]{Nordentoft20} in full generality. In these cases, the maps~$h_\gamma$ are bounded.
\item When~$\phi$ is the central Eisenstein series of level~$1$, this was proved in~\cite[Lemma~10]{Bettin16} using a functional equation for the associated Dirichlet series.
\end{itemize}

In all other cases, Theorems~\ref{thm:intro-reg} and~\ref{thm:intro-growth} are new. Compared with these earlier works, the main point of Theorem~\ref{thm:intro-reg} is that it does not require holomorphicity, nor does it use the functional equation for the associated $L$-series. In particular it does not require the presence of Fricke involutions.

\subsection{Normal distribution for \texorpdfstring{$\SL(2, \Z)$}{SL(2,Z)}}

For reasons that will be clarified below, we restrict to~$\Gamma = \SL(2, \Z)$ in this section. In this precise case, it was shown in~\cite{BettinDrappeau2019} using dynamical methods that a statement of the kind given by Theorem~\ref{thm:intro-reg} yields the limiting distribution for the multisets~$\cD_\phi(q)$ on average over~$q$ (consult Section \ref{sec:context} for an overview of the dynamical input). We will deduce the following, which proves the averaged version of Conjecture~\ref{conj:addtwist} for~$\Gamma=\SL(2,\Z)$.

\begin{thm}[Additive twists on average for~$\SL(2,\Z)$]\label{thm:intro-normal-level1}
  Let~$\phi$ be a Maa\ss{} cusp form for~$\Gamma = \SL(2, \Z)$. Then the multisets
  $$ \Big\{ \frac{L_\phi(a/q)}{(\log q)^{1/2}}: a\in (\Z/q\Z)^\times, q\leq Q\Big\} $$
  become distributed, as~$Q\to\infty$, to a centered normal law.
\end{thm}

On the other hand, we believe the restriction to~$\Gamma = \SL(2, \Z)$ to be artificial. More precisely, we conjecture the following.
\begin{conj}
  The statement of Theorem~\ref{thm:intro-normal-level1} holds when~$\Gamma$ is replaced by an arbitrary cofinite Fuchsian group with cusps, and~$\phi$ is a Maa\ss{} cusp form for~$\Gamma$ of integer weight.
\end{conj}

As we have mentioned already, this conjecture is known when~$\phi$ is associated with a holomorphic form. At the present time we lack a proper analogue of the methods of~\cite{BettinDrappeau2019, LeeSun2019} which would allow to handle the general case. This is the subject of work in progress of Bettin, Lee and the first named author.
\begin{rem}
The cuspidality hypothesis is not essential to the method, however the non-cuspidal Maa\ss{} forms for~$\SL(2, \Z)$ are essentially spanned by an Eisenstein series, and for this series we can reduce to the case of the Estermann function treated in~\cite[section~9.2]{BettinDrappeau2019}.
\end{rem}
\begin{rem}
Given a finite orthogonal family~$(\phi_j)$ of Hecke-Maa\ss{} cusp forms for~$\SL(2, \Z)$, we obtain more generally the joint convergence to independent normal distributions of the values~$(L_{\phi_j}(a/q))_j$, see Corollary~\ref{cor:joint-CLT} below.
\end{rem}

\subsection{Arithmetic applications}

In the special case where $\phi$ is a Hecke--Maa{\ss} form and~$\Gamma = \Gamma_0(N)$ is a congruence group, the additive twists $L_\phi(x)$ are connected to the central values of the twisted $L$-functions $L(\phi,\chi,1/2)$ using orthogonality of characters (known as the Birch--Stevens formula) alluded to above. Here
$$L(\phi,\chi,s)=\sum_{n\geq 1} \frac{\lambda_\phi(n)\chi(n)}{n^s}, \qquad (\Re s>1),$$
and elsewhere by meromorphic continuation, where $\chi$ is a Dirichlet character and $\lambda_\phi(n) = a(n)\sqrt{n}$ denotes the $n$-th Hecke eigenvalue of $\phi$ (we reserve the notation~$\lambda_\phi(n)$ to cases when~$\phi$ is a Hecke eigenform). Using this connection we get a number of applications to twisted $L$-functions of Theorems~\ref{thm:intro-reg} and \ref{thm:intro-growth}

\subsubsection{Reciprocity formul\ae} 

In an unpublished preprint \cite{Conrey07}, Conrey proved a certain {\lq\lq}reciprocity relation{\rq\rq} for twisted second moments of Dirichlet $L$-functions relating the following two quantities;
\begin{equation}\label{eq:conrey-recip}
  \sum_{\chi \mod p} |L(\chi,1/2)|^2\chi(\ell) \rightsquigarrow  \sum_{\chi \mod \ell} |L(\chi,1/2)|^2\chi(-p),
\end{equation}
where $L(\chi,s)=\sum_{n\geq 1}\chi(n) n^{-s}$ for $\Re s>1$ and elsewhere by analytic continuation. The results were later extended by Young \cite{Young11} and Bettin \cite{Bettin16}. This can be seen as the $\GL_2\times \GL_1$-case (with the $\GL_2$-form being an Eisenstein series) of the phenomena of \emph{spectral reciprocity} investigated in \cite{AndersenKiral18},\cite{BlomerKhan19},\cite{BlomerLiMiller19}. The quantum modularity results we obtain resolve completely the $\GL_2\times \GL_1$ (over $\Q$) case and we obtain a relation of the type
$$ \sum_{\chi \mod p} \tau(\overline{\chi})L(\phi,\chi,1/2)\chi(\ell) \rightsquigarrow   \sum_{\chi \mod N\ell} \tau(\overline{\chi})L(\overline{\phi},\chi,1/2)\chi(-p),$$
where $\phi$ is a ($\GL_2$) Hecke--Maa{\ss} newform of level $N$. We refer to Theorem \ref{recipr} (cuspidal case) and Theorem \ref{thm:reciprnoncuspidal} (Eisenstein case) for the exact statements. In the special case of $\phi$ being cuspidal of level $1$, we get the following result. 
\begin{cor}  \label{cor:recilevel1.1}
  Let $\phi$ be a Hecke--Maa{\ss} cuspform for~$\Gamma=\SL(2, \Z)$ whose Fourier coefficients satisfy  $a_\phi(-n)=\epsilon_\phi a_\phi(n)$ with $\epsilon_\phi\in \{\pm 1\}$. Then for any pair of primes $0<p<\ell$ and any choice of sign~$\eta \in \{\pm 1\}$, we have
  \begin{align}
    \label{eq:recilevel1}&\frac{2}{p-1}\,\,\, \sideset{}{^{\ast,\eta}}\sum_{\chi \modulo p}\tau(\overline{\chi}) L(\phi, \chi, 1/2)\chi(\ell)\\
    \nonumber &  \qquad \qquad - \eta \frac{2}{\ell-1}\,\,\, \sideset{}{^{\ast,\eta }}\sum_{\chi \modulo \ell}\tau(\overline{\chi}) L(\phi, \chi,1/2)\chi(p) = M_{\phi,\eta}+O_{\phi,\eps}((p/\ell)^{1-\eps}+p^{\theta-1+\eps}),
  \end{align}
  where
  $$M_{\phi,\eta}=
  \begin{cases}  
    - \epsilon_\phi L(\phi,1/2), & \eta=+1,\\
    0, & \eta=-1,\\
  \end{cases}$$
  and  $\theta= \frac{7}{64}$ is the best bound towards the Ramanujan--Petersson conjecture for Maa{\ss} forms due to Kim and Sarnak \cite{KiSa03}.
  Here the decorations on the sums means that we restrict to primitive characters with $\chi(-1)=\eta$, and $L(\phi,s)$ denotes the (standard) $L$-function of $\phi$. 
\end{cor}

\begin{rem}
  When~$\phi$ is not cuspidal, a similar statement holds with an altered right-hand side.
  Choosing~$\phi$ to be the Eisenstein series~$E_{1, 1}^*(z, \frac12)$ defined in Section~\ref{sec:examples-eisenstein} below, we find that~$a(n) = d(n)$, the divisor function, for~$n>0$, and therefore, for all~$\chi\pmod{p}$ primitive,
  $$ \tau(\overline \chi) L(\phi, \chi, 1/2) = \tau(\overline \chi) L(\chi, 1/2)^2 = i^a \sqrt{p} \abs{L(\chi, 1/2)}^2 $$
  where~$a\in\{0, 1\}$ depends on~$\chi(-1)$. In this way we can deduce a form of Conrey's reciprocity formula~\cite{Conrey07} alluded to above~\eqref{eq:conrey-recip}.
\end{rem}

\subsubsection{Wide moments of Dirichlet twists}%%%%%%%%%%%%%%%%%%%%%%

In the case of $\phi$ a Hecke--Maa{\ss} cuspform of level $1$, we not only obtain the normal distribution result for the additive twists in Theorem~\ref{thm:intro-normal-level1}, but furthermore a convergence of moments (see Proposition \ref{prop:moments-level1}). Using the Birch--Stevens relations this implies certain new moment calculations for the twisted $L$-functions $L(\phi,\chi,1/2)$ which have not been obtained by the standard {\lq\lq}approximate functional equation{\rq\rq}-approach. This fits into the framework of \emph{wide moments} of families of $\GL_1$-twists of automorphic $L$-functions as in \cite{Nordentoft21}, \cite{Nordentoft21.2}, \cite{Be17}, \cite[Corollary 1.9]{Nordentoft2021} (see also Section \ref{sec:arith-applications} below for more background). We state here the moment calculation in the simplest version and refer to Corollary \ref{cor:widemoments} for the most general statement which is a rare example of a moment calculation with a huge amount of cancellation.

\begin{cor} \label{cor:widemomentssimplified2}
  Let $\phi$ be a Hecke--Maa{\ss} cusp form for~$\Gamma = \SL(2, \Z)$ and  $n\in  2 \N_{\geq 0}$.  Then we have as $Q\rightarrow \infty$
  \begin{align}
    \label{eq:widemoment2}\sum_{\substack{0<c\leq Q}} \frac{1}{\varphi(c)^{n-1}} \sum_{\substack{\chi_i \modulo c,\\ 1\leq i\leq n:\\ \chi_1\dotsb \chi_n=\mathbf{1}}}\, \prod_{i=1}^n  \nu_\phi(\chi_i) L(\phi, \chi_{i},1/2)
    = P(\log Q) Q^2+O_{\phi,n}(Q^{2-\delta}),
  \end{align}
  for some $\delta>0$, where $P$ is a degree $n/2$ polynomial with leading coefficient 
  $$ 2^{n/2} (n/2)! L(\sym^2 \phi, 1)^{n/2}. $$ 
  Here $\mathbf{1}$ denotes the principal character (of the relevant modulus suppressed in the notation), and the factors $\nu_\phi(\chi)$ are certain local weights essentially of size $c^{1/2}$ for $\chi \modulo c$.
\end{cor} 

We refer to Corollary \ref{cor:widemoments} and equation~\eqref{eq:weightBS} for precise expressions of~$\nu_\phi(\chi)$. When~$\chi\pmod{c}$ is primitive, then we simply have~$\nu_\phi(\chi) = \tau(\overline \chi)$.

\begin{rem}
  Assuming the Lindel\"{o}f bound $L(\phi,\chi,1/2)\ll_{\phi, \eps} c^\eps$ for $\chi \modulo c$, together with the Ramanujan--Petersson conjecture $\lambda_\phi(n)\ll_{\phi, \eps} n^\eps$ one gets the {\lq\lq}trivial{\rq\rq} bound 
  $O_{\phi,\eps}(Q^{n+1+\eps})$ for the left-hand side of~\eqref{eq:widemoment2} (using also the bound $\nu_\phi(\chi)\ll_\eps c^{1/2+\eps}$ coming from \eqref{eq:boundnu} which is essentially sharp). Thus we see that for $n>1$, there is massive cancellation in the sum. In particular, it appears to be very hard to obtain such a result using an {\lq\lq}approximate functional equation{\rq\rq}-approach (as in e.g. \cite{BlFoKoMiMiSa18}).
\end{rem}

\subsection{Context and overview}\label{sec:context}

Instead of~$L_\phi(x)$, we start our analysis with the Eichler integrals
$$ \E_\phi(x) \simeq \int_x^\infty \phi(z) \frac{\df z}{\Im(z)}, $$
see Section~\ref{sec:def-eichler}.
We establish first the analogues of Theorems~\ref{thm:intro-reg} and~\ref{thm:intro-growth} with~$L_\phi(x)$ replaced by~$\E_\phi(x)$, so that the function of interest becomes (we ignore in this sketch issues concerning nebentypus or regularizations)
\begin{equation}\label{eq:hE-diff}
    h^\E_\gamma(x) \simeq \int_{\gamma x}^\infty \phi(z) \frac{\df z}{\Im(z)} - \int_x^\infty \phi(z) \frac{\df z}{\Im(z)}.
\end{equation}
Our proof departs from the earlier proofs in~\cite{Bettin16, BettinDrappeau2019, Nordentoft20}:
\begin{itemize}
    \item the proofs in~\cite{BettinDrappeau2019, Nordentoft20} use the holomorphicity of the underlying modular form to perform a contour integral deformation (see~\cite[p.~1414]{BettinDrappeau2019} and~\cite[p.~158]{Nordentoft20}). This yields an expression for~$h_\gamma(x)$ as a linear combination of translates~$L_\phi(\phi, 1/2+j, x)$ for~$j\geq 1$, from which Theorems~\ref{thm:intro-reg} and~\ref{thm:intro-growth} follow easily in the holomorphic case.
    \item the proof in~\cite{Bettin16} proceeds with Mellin transformation to reduce the question to estimating an integral of~$L(\phi, s)$ over~$s$, and uses matching of Gamma factors on both sides of the functional equation, see \cite[pp.~6900--6903]{Bettin16}.
\end{itemize}

The proof we give of Theorem~\ref{thm:intro-reg} can be viewed as in-between~\cite{BettinDrappeau2019, Nordentoft20} and~\cite{Bettin16}. We use modularity in the integrals~\eqref{eq:hE-diff}, however corresponding paths in both resulting integrals in~\eqref{eq:hE-diff} do not systematically vanish (due to possible lack of holomorphy). We analyze the behaviour in~$x$ by a careful splitting-and-matching of these integrals. This is the subject of Section~\ref{sec:qmod-Eichler}.

The proof we provide of Theorem~\ref{thm:intro-growth}, in Section~\ref{sec:hgamma-at-infinity}, does not easily relate to earlier known cases; in the holomorphic cuspidal case~\cite{BettinDrappeau2019, Nordentoft20} the corresponding statement is straightforward, while in the level-$1$ Eisenstein case, the proof in~\cite{Bettin16} (which starts with a Mellin integral) is set up in a different way as the proof we follow. We refer to Section~\ref{sec:hgamma-at-infinity} for more details.

For both Theorems~\ref{thm:intro-reg} and~\ref{thm:intro-growth} (more precisely, their analogues for Eichler integrals), we 
also provide a geometric and conceptually simpler proof in the special case of forms of weight~$2$ using the Stokes formula, which unfortunately we couldn't generalize to arbitrary weight, but which helps motivate our later arguments.
At this point we obtain the analogues of Theorems~\ref{thm:intro-reg} and~\ref{thm:intro-growth} for the Eichler-integral analogue~\eqref{eq:hE-diff}.

To carry the arguments over to~$L_\phi(x)$ itself, we express~$L_\phi(x)$ as linear combinations of Eichler integrals. In the earlier known cases of Theorems~\ref{thm:intro-reg} and~\ref{thm:intro-growth}, this issue was nonexistent, due notably to the vanishing of negative-index Fourier coefficients of holomorphic forms; however a similar issue was met for instance in Section~8 of~\cite[p.~530]{DukeEtAl2002} where the authors express~$L_j(s)$ in terms of $\Psi_j(s)$ using a certain conjugation operator. In this paper, we proceed differently: we use level-raising and level-lowering operators to express~$L_\phi(x)$ as a combination of two Eichler integrals. This is the subject of Section~\ref{sec:qmod-L}, and the core of the proof is the analysis of certain integrals of Whittaker functions, which requires some work with hypergeometric functions (deferred in Appendix~\ref{app:hypergeom}). The approach we use has the advantage over~\cite{DukeEtAl2002} that it does not assume~$\Gamma$ is self-conjugate under~$(\begin{smallmatrix} -1& 0 \\ 0 & 1\end{smallmatrix})$. As a by-product of our computations with integrals of Whittaker functions, we also prove the functional equation for~$L(\phi, x, s)$ in Section~\ref{sec:fe}, which is new in the odd $k$ case.
We present a few examples (Hecke-Maa\ss{} cusp forms and Eisenstein series) in Section~\ref{sec:examples}.

In Section~\ref{sec:normal-distrib}, we prove Theorem~\ref{thm:intro-normal-level1}. The existence of the limit law follows easily from~\cite[Theorem~3.1]{BettinDrappeau2019}, using the quantitative versions of Theorems~\ref{thm:intro-reg} and~\ref{thm:intro-growth} proven in Sections~\ref{sec:qmod-Eichler}--\ref{sec:qmod-L}. To get an idea as to why quantum modularity (Theorem~\ref{thm:intro-reg}) is useful for obtaining normal distribution for $\Gamma=\SL_2(\Z)$, consider a rational number $x\in \Q$. Using the definition~\eqref{eq:def-hgamma-intro} for the matrix $\gamma =\begin{psmallmatrix}0&-1\\ 1&0\end{psmallmatrix}$, along with $1$-periodicity, we obtain
\begin{align}
        f(x)=-h_\gamma(x)+f(\gamma x)= -h_\gamma(x)+f(-T(x)),
\end{align}
where $T(x):= \{ \tfrac{1}{x}\} $ denotes the Gauss map (here $\{x\}\in [0,1)$ denotes the fractional part of $x$). Iterating this we arrive at 
\begin{align}
    f(x)= -h_S(x)+h_S(-T(x))-\ldots \pm h_S((-1)^{r(x)-1} T^{(r(x)-1)}(x)) + f(0)
\end{align}
where $r(x)$ denotes the length of the continued fraction expansion of $x$ (so that $T^{(r(x))}(x)=0$). Setting aside the minus signs in the arguments of~$h_\gamma$, we see that the value distribution of the quantum modular form $f$ is determined by the sum of a Hölder-continuous map along iterates of the Gauss map. This was studied in the celebrated work of Baladi and Vall\'{e}e \cite{BaladiVallee2005}. There, the analogue of $h_\gamma$ is called the cost function, and was assumed to be constant on intervals of the shape $(\tfrac{1}{n+1},\tfrac{1}{n})$. This is however not true for the cost functions $h_\gamma$ as above. Exactly motivated by this example, the extension to general regular cost functions (of moderate growth) was carried out by Bettin and the first named author in \cite{BettinDrappeau2019}, which is what we use here.
The computation of the parameters of the limit law (mean and variance) takes some additional work, using standard methods of analytic number theory. The functional equation proven in Section~\ref{sec:fe} is used there. 

In Section~\ref{sec:arith-applications}, we restrict to the case of congruence subgroups $\Gamma=\Gamma_0(N)$ and Hecke eigenforms to deduce consequences for the arithmetically interesting twisted $L$-functions $L(\phi,\chi,1/2)$, namely 1) spectral reciprocity for twisted $L$-functions as in Corollary \ref{cor:recilevel1.1} which follows from the quantum modularity for the Fricke involution combined with the Birch--Stevens formula from Section \ref{sec:BirchStevens} and 2) asymptotic evaluation of wide moments of twisted $L$-functions  as in Corollary \ref{cor:widemomentssimplified2} which follows from the asymptotic formula for the moments of the additive twist $L$-series combined with the Birch--Stevens formula.   

% \subsection{Structure of the paper}

% In Section~\ref{sec:setting}, we set the background and main definitions of the automorphic forms we will deal with.
% In Section~\ref{sec:qmod-Eichler}, we establish the quantum modularity for the discrepancy function~$h_\gamma$ associated with Eichler integrals.
% In Section~\ref{sec:hgamma-at-infinity}, we investigate the behaviour at infinity of~$h_\gamma$ for the Eichler integrals.
% In Section~\ref{sec:qmod-L}, we transfer the quantum modularity and behaviour at infinity, from the Eichler integrals, to the actual central~$L$-values, which will prove Theorems~\ref{thm:intro-reg} and \ref{thm:intro-growth}.
% In Section~\ref{sec:examples}, we apply our main results to a few examples in congruence groups: Hecke-Maa\ss{} cusp forms, and real-analytic Eisenstein series.
% In Section~\ref{sec:normal-distrib}, we deduce the convergence in law in Theorem~\ref{thm:intro-normal-level1}.
% In Section~\ref{sec:arith-applications}, we deduce the arithmetic applications to reciprocity formul\ae{} (Corollary~\ref{cor:recilevel1.1}) and wide moments (Corollary~\ref{cor:widemomentssimplified2}).
% Appendix~\ref{app:hypergeom} contains two lemmas about hypergeometric functions which are used in Section~\ref{sec:qmod-L}.

\subsection{Notation}

We use indistinctively the symbols~$X = O(Y)$ and~$X \ll Y$ to indicate the existence of a constant~$C>0$ such that~$\abs{X} \leq C Y$. The value of~$C$ may depend at most on variables which are either indicated in subscript, as in \emph{e.g.}~$X \ll_\eps Y$, or mentioned in the immediate context. The symbol~$X\asymp Y$ means~$X\ll Y$ and~$Y\ll X$. The letter~$\eps$ denotes an arbitrarily small quantity, which may differ between occurences.

\section{Background}\label{sec:setting}

\subsection{Maa\ss{} forms}
For a detailed account of the following material we refer to~\cite[Chapter~2]{Iw}, \cite[Section 4]{DukeEtAl2002} as well as the classical sources \cite{Roelcke67}, \cite{Selberg56}, \cite{Katok92}.
Fix $\Gamma \subset \SL(2,\R)$ a discrete, co-finite subgroup with a cusp at $\infty$, a character $\chi:\Gamma \rightarrow \C$ trivial on all parabolic elements of $\Gamma$ and an integer $k\in \Z_{\geq 0}$.
Denote by $C(\Gamma)\subset \mathbb{P}^1(\R)$ the cusps of $\Gamma$, and for~$x\in C(\Gamma)$ let~$\Gamma_x$ be the stabiliser of~$x$.
For~$x\in C(\Gamma)$, a scaling matrix~$\sigma_x$ for~$x$ is any matrix which satisfies~$\sigma_x^{-1} \Gamma_x \sigma_x = \{\begin{psmallmatrix} 1 & \Z \\ & 1\end{psmallmatrix} \}$.
We assume~$\sigma_\infty = \id$ (this can always be ensured by conjugating~$\Gamma$ by a diagonal element).

Let $k\in \Z_{\geq 0}$. We denote by $\A(\Gamma, \chi , k)$ the vector space of all \emph{weight $k$ automorphic forms of $\Gamma$ with nebentypus $\chi$}, i.e. smooth maps $\phi: \H \rightarrow \C$ satisfying:

\begin{itemize}
\item[$(H_1)$] For all~$\gamma \in \Gamma$ and~$z\in\H$, $\phi(\gamma z) = u_\gamma(z) \phi(z)$, where
  \begin{equation}
    u_\gamma(z)=j_\gamma(z)^k \chi(\gamma), \qquad j_\gamma(z) = \frac{j(\gamma, z)}{|j(\gamma,z)|}= \frac{cz+d}{|cz+d|}.\label{eq:def-ugamma}
  \end{equation}
\end{itemize}
Notice that if $-\id\in \Gamma$ then we must have the compatibility condition $\chi(-1)=(-1)^k$ for $\A(\Gamma, \chi , k)$ to be non-trivial. 
We borrow the analytic notations from~\cite{DukeEtAl2002}. 
Let
\begin{equation}\label{eq:levelraising} R_k = \frac k2 + (z-\bar z) \frac{\partial}{\partial z},\quad \Lambda_k = \frac k2 + (z-\bar z) \frac{\partial}{\partial \overline{z}}  \end{equation}
be respectively the weight~$k$ level raising and level lowering operator, as defined in~\cite[eqs.~(4.3)-(4.4)]{DukeEtAl2002}. These define maps 
$$R_k:\A(\Gamma, \chi , k)\rightarrow \A(\Gamma, \chi , k+2),\quad \Lambda_k:\A(\Gamma, \chi , k)\rightarrow \A(\Gamma, \chi , k-2).$$ The weight~$k$ Laplacian is defined by
$$ \Delta_k =-R_{k-2} \Lambda_{k}-\frac{k}{2}\left(1-\frac{k}{2}\right)=-\Lambda_{k+2} R_{k}+\frac{k}{2}\left(1+\frac{k}{2}\right)= y^2\Big(\frac{\partial^2}{\partial x^2} + \frac{\partial^2}{\partial y^2}\Big) - iky \frac{\partial}{\partial x}. $$
For $k\in \Z_{\geq 0}$ and $s\in \C$, we define  an operator 
\begin{equation}\label{eq:involution}
  Q_{s,k}:\A(\Gamma, \chi , k)\rightarrow \A(\Gamma, \chi, k),
\end{equation} 
as in~\cite[eq.~(4.65)]{DukeEtAl2002} by 
$$ (Q_{s,k}\phi)(z)=\frac{\Gamma(s-k/2)}{\Gamma(s+k/2)} (\Lambda_{-k+2} \cdots\Lambda_{k-2} \Lambda_{k} \phi)(-\overline{z}),   $$
where we put $Q_{s,k}=0$ if~$s\in k/2+\Z_{\leq 0}$. Notice that for $k=0$ we have $(Q_{s,0}\phi)(z)=\phi(-\overline{z})$ which is the usual reflection operator. The operator $Q_{s,k}$ preserves the eigenspace of $\Delta_k$ with eigenvalue $s(1-s)$, and is an involution for $s\not\in k/2+\Z_{\leq 0}$.
Similarly, we define $Q_{s,k}$ for negative $k$ using the raising operators.  

We say that $\phi \in \A(\Gamma, \chi , k)$ is a \emph{Maa{\ss} form} if it satisfies
\begin{itemize}
\item[$(H_2)$] For all~$x \in C(\Gamma)$, $\phi(\sigma_x z) = o(e^{2\pi y})$ as $y = \Im z\rightarrow \infty$.
\item[$(H_3)$] $\phi$ is an eigenfunction of $\Delta_k$ with eigenvalue $\lambda_\phi = s_\phi(1-s_\phi) = 1/4+t_\phi^2$ with $\Re s_\phi \geq 1/2$ and $s_\phi = 1/2+it_\phi$.
\item[$(H_4)$] $\phi$ is an eigenfunction of the operator~$Q_{s_\phi,k}$ with eigenvalue~$\epsilon_\phi \in \{\pm 1,0\}$ (this letter~$\epsilon_\phi$ will always be written with a subscript indicating the corresponding form, to avoid confusion with the notation~$\eps$ for an arbitrary small number). If $k=0$ then $\epsilon_\phi$ is the \emph{sign} of the Maa{\ss} form~\cite[p.~106]{Bump1997}, and if $\phi$ comes from a holomorphic form then $\epsilon_\phi=0$. Note that when~$k$ is odd, $\epsilon_\phi$ depends on the choice of sign of~$t_\phi$ in the definition of~$Q_{s_\phi, k}$.
\end{itemize}

The above conditions imply~\cite{Iw, DukeEtAl2002} that at any cusp $x\in C(\Gamma)$ we have the Fourier expansion
\begin{equation}\label{fexp}
  j_{\sigma_x}(z)^{-k} \phi(\sigma_x z)= \phi_x(\Im z)+\sum_{n\neq 0} a_x(n) e(n\Re z) W_{\frac k2 \sgn(n),it_\phi}(4\pi |n|\Im z),
\end{equation}
where $W_{\alpha,\beta}:\R_{>0}\rightarrow \C$ is the weight $\alpha$-Whittaker function with parameter $\beta$, i.e. the unique solution $W$ to
$$  \frac{d^2 W}{dy^2}+ \left(-\frac{1}{4}+\frac{\alpha}{y}+\frac{1/4-\beta^2}{y^2}\right)W=0, $$
satisfying $W(y)\sim y^{\alpha}e^{-y/2}$ as $y\rightarrow \infty$ (with $\alpha,\beta$ fixed). In particular we have $W_{0,it}(4\pi y)=2 y^{1/2}K_{it}(2\pi y)$ where $K_s(y)$ denotes the K-Bessel function. The constant term~$\phi_x(y)$ (with~$y=\Im z$ above) is given by
\begin{equation}
  \phi_x(y)= \begin{cases}A_x y^{s_\phi}+B_x y^{1-s_\phi}, \quad s_\phi\neq 1/2\\ A_x y^{1/2}+B_x y^{1/2}\log y , \quad s_\phi= 1/2. \end{cases}\label{eq:def-phi-cuspidalpart}
\end{equation}
We will say that $\phi$ is \emph{cuspidal at the cusp $x$} meaning that $\phi_x(y)=0$. Finally we say that $\phi$ is a \emph{Maa{\ss} cuspform} if it is cuspidal at all cusps $x\in C(\Gamma)$. 

When~$\phi$ is square-integrable, then from~\cite[Corollary~4.4]{DukeEtAl2002} we know that~$\phi$ may arise in two ways~:
\begin{enumerate}[(i)]
\item If~$1/2\leq \Re(s_\phi)<1$ and~$s_\phi\neq \frac12$, then~$\phi$ is obtained from repeated applications of level-raising or lowering operators from a weight~$0$ or~$1$ Maa{\ss} form, depending on the parity of~$k$,
\item Otherwise~$s_\phi = \ell/2$ with~$\ell\equiv k\pmod{2}$, and~$\phi$ is then associated, through level-raising or lowering operators, to a form~$\psi$ of weight~$\ell$ for which~$z\mapsto y^{-k/2} \psi(z)$ is holomorphic. In this case, we necessarily have~$B_x = 0$ in~\eqref{eq:def-phi-cuspidalpart}.
\end{enumerate}
We won't need to assume, in the present work, that~$\phi$ is square-integrable. However, in case it is not, we will assume that~$\Re(s_\phi) = 1/2$, as this will simplify some of our statements. This is summarized in the following hypothesis:
\begin{itemize}
  \item[$(H_5)$] Either~$\phi$ is square-integrable (and therefore (i) or (ii) above hold), or~$\Re(s_\phi) = 1/2$.
\end{itemize}

We will abbreviate throughout
$$ a(n) = a_\infty(n). $$
Moreover, we have by~\cite[eq.~(4.70)]{DukeEtAl2002}
\begin{equation}\label{eq:fouriercoeff-sign}
  a(-n) = \epsilon_\phi \frac{\Gamma(s_\phi+\frac k2)}{\Gamma(s_\phi-\frac k2)} a(n).
\end{equation}
For~$s_\phi = \ell/2$ with~$\ell\equiv k\pmod{2}$, this says that~$a(-n) = 0$ for~$n>0$.

Finally, we will assume the following bound.
\begin{itemize}
\item[$(H_6)$] The Fourier coefficients~$a(n)$ in~\eqref{fexp} at~$x=\infty$ satisfy
  \begin{equation}
    \sum_{1\leq n\leq X} |a(n)|\ll_{\phi, \eps} X^{1/2+\eps}.\label{eq:bound-an-rough}
  \end{equation}
\end{itemize} 
This condition is satisfied when~$\phi$ is square-integrable. This is obtained by Cauchy--Schwarz combined with a straightforward modification of~\cite[Theorem 3.2]{Iw}.
We also know that it holds true in the case of congruence groupe~$\Gamma = \Gamma_0(q)$ for Eisenstein series with~$\Re(s_\phi) = 1/2$, as then the coefficients~$a(n)$ are essentially of the shape~$d(n) / \sqrt{n}$, with~$d(n)$ being the divisor function or a variant thereof.

\begin{rem}\label{rmk:cuspidal-bound-1}
  When~$\phi$ is cuspidal, it follows from a classical argument~\cite[Theorem~8.1]{Iw} (Wilton's bound) that the following estimate
  \begin{equation}
    \sum_{1\leq n\leq X} a(n) \e(nt) \ll_{\phi, \eps} X^\eps \qquad (t\in \R, \phi \text{ cuspidal})\label{eq:bound-an-wilton}
  \end{equation}
  holds uniformly in~$t$.
\end{rem}

We will use the following rough but uniform bounds on~$\phi$.

\begin{lemma}\label{lem:bounds-phi}
  For all~$x, x', \eta \in \R$ and~$0<y\ll 1$, we have
  \begin{align}
    &|\phi(x+iy)| \ll_\eps y^{-1/2-\eps}, \label{eq:bound-phi} \\
    &|\phi(x+iy)-\phi(x'+iy)| \ll_\eps y^{-1/2-\eps} \min(1, \tfrac{|x-x'|}y), \label{eq:bound-phi-diff} \\
    &|\phi(x+iy)-\phi(x'+iy)-\phi(x+\eta+iy)+\phi(x'+\eta+iy)| \ll_\eps y^{-1/2-\eps} \min(1, \tfrac{|x-x'|}y) \min(1, \tfrac{\abs{\eta}}y). \label{eq:bound-phi-doublediff}
  \end{align}
  If, additionnally, $\phi$ is cuspidal, then all three of the bounds hold with the factor~$y^{-1/2-\eps}$ on the right-hand side replaced by~$y^{-\eps}$.
\end{lemma}

\begin{proof}
  The bound~\eqref{eq:bound-phi} follows at once by using hypothesis~\eqref{eq:bound-an-rough} in the Fourier expansion~\eqref{fexp}, partial summation and the exponential decay of the Whittaker function. 
  Similarly, the bound~\eqref{eq:bound-phi-diff} follows from~\eqref{eq:bound-phi} if~$|x-x'|>y$, and otherwise we write~$\e(nx)-\e(nx') = 2\pi i n \int_{x'}^x \e(nt)\df t$ in the Fourier expansion and conclude again by the Fourier decay of the Whittaker function.
  Finally, the bound~\eqref{eq:bound-phi-doublediff} follows from~\eqref{eq:bound-phi} and \eqref{eq:bound-phi-diff} if~$|x-x'|>y$ or~$\eta>y$, and otherwise we write
  $$ \e(nx)-\e(nx')-\e(n(x+\eta))+\e(n(x'+\eta) = (\e(nx) - \e(nx'))(1-\e(n\eta)) $$
  and use the partial summation again.

  If~$\phi$ is cuspidal, then instead of the bound~\eqref{eq:bound-an-rough}, we use Wilton's bound~\eqref{eq:bound-an-wilton} in the argument above. Note that in this case~$\phi$ is actually bounded.
\end{proof}

\subsection{The Eichler integral}\label{sec:def-eichler}

For $x\in C(\Gamma)\setminus\{\infty\}$, we define
\begin{equation}
  \E(\phi, x, s) := \int_0^\infty (\phi(x+iy)-\phi_\infty(y)) y^{s-1/2}\frac{dy}{y},\label{eq:def-eichlerint}
\end{equation}
which converges absolutely for $\Re s>1$ by the exponential decay of the Whittaker function in~\eqref{fexp}. This is a generalization of the original Eichler integrals~\cite{Eichler} which were associated to cuspidal holomorphic forms and had $s=(k-1)/2$.
This kind of integrals originate from Riemann's memoir~\cite{Riemann}. The special case~$x=0$ was considered by Hecke~\cite{Hecke36} to establish the functional equation of Hecke $L$-function of holomorphic cusp forms. The idea of studying the analytic properties in the variable~$x$ is due to Eichler~\cite{Eichler} and led to the development of the Eichler-Shimura isomorphism; see~\cite{DiamantisRolen2018} for references.

Here we will need to study the value at~$s=1/2$, whose existence we first deduce from analytic continuation.

\begin{prop}\label{prop:def-twisted-l}
  The Eichler integral $\E(\phi, x, s)$ admits meromorphic continuation to the entire complex plane with possible poles contained in $\{\pm i t_\phi, 1 \pm i t_\phi\}$. If~$\Re(s_\phi)\geq 1$, then the only possible poles are at~$\{1+it_\phi, -it_\phi\}$.
\end{prop}

\begin{proof}
  By the Fourier expansion of $\phi$ at $x$, we know that 
  $$ \phi(x+iy) = j_{\sigma_x}(\sigma_x^{-1}(x+iy))^k \phi_x(\Im(\sigma_x^{-1}(x+iy)))+ g_x(iy) = (i\sgn(c_x))^k \phi_x((c_x^2 y)^{-1})+ g_x(y),$$
  where $c_x$ is the bottom-left coefficient of a scaling matrix for~$x$ and, $g_x(y)\ll_A y^{A}$ for all $A>0$ as $y\rightarrow 0$. Let $\psi: (0,\infty)\rightarrow \R$ be smooth and decreasing with $\psi(y)=1$ for $y<1$ and   $\psi(y)=0$ for $y>2$. Then 
  \begin{equation}\label{eq:Ex-regular-part}
     s\mapsto \int_0^\infty (\phi(x+iy)-(1-\psi(y))\phi_\infty(y)-\psi(y) (i \sgn(c_x))^k \phi_x((c_x^2y)^{-1}))y^{s-1/2}\frac{dy}{y}
  \end{equation}
  extends to an entire function due to the rapid decay of the integrand as $y\to 0$ and~$y\to \infty$. Moreover we see that for $\Re s>1/2+\Re s_\phi$, we have by partial integration:
  $$  \int_0^\infty\psi(y)\phi_x((c_x^2 y)^{-1}) y^{s-1/2}\frac{dy}{y}=- \int_0^\infty\psi'(y) F_x(y)dy,  $$
  where $F_x$ is the antiderivative of $y\mapsto \phi_x((c_x^2y)^{-1})y^{s-3/2}$, which is of the form
  $$\begin{cases} \frac{A}{s-1/2 -s_\phi} y^{s-1/2 -s_\phi}+\frac{B}{s-3/2+s_\phi} y^{s-3/2+s_\phi},& s_\phi\neq 1/2\\ \frac{A'}{s-1} y^{s-1} + \frac{B}{(s-1)^2}(1 - (s-1)\log y)y^{s-1}, & s_\phi= 1/2, \end{cases}$$
  where $A = A_x \abs{c_x}^{-2s_\phi}$, $B = B_x \abs{c_x}^{2s_\phi - 2}$, and $A' = (A_x - 2B_x \log\abs{c_x})/\abs{c_x}$. This defines a meromorphic continuation to the entire complex plane with possible poles only at $s\in \{1/2+s_\phi, 3/2-s_\phi \}$. When~$\Re(s_\phi)\geq 1$, then~$\phi$ is associated with a holomorphic form, and we have already noted in this case that necessarily~$B = 0$, and therefore there is no pole at~$3/2-s_\phi$. Carrying out a similar computation for
  $$ -\int_0^\infty \psi(y) \phi_\infty(y) y^{s-1/2} \frac{\df y}{y}, $$
  we arrive to the conclusion that~$s\mapsto \E(\phi, x, s)$ has at most simple poles at~$\{1/2+s_\phi, 3/2-s_\phi, 1/2-s_\phi, s_\phi-1/2 \}$ when~$s_\phi \neq 1/2$, with residues given by
  \begin{align*}
    &\Res_{s=1/2+s_\phi} = (i\sgn(c_x))^k \frac{A_x}{\abs{c_x}^{2s_\phi}}, && \Res_{s=3/2-s_\phi} = (i \sgn(c_x))^k \frac{B_x}{\abs{c_x}^{2(1-s_\phi)}}, \\
    &\Res_{s=1/2-s_\phi} = -A_\infty, && \Res_{s=s_\phi-1/2} = -B_\infty;
  \end{align*}
  while for~$s_\phi = 1/2$, we have
  \begin{align*}
      & \E(\phi, x, s) = (i \sgn(c_x))^k \Big( \frac{B_x}{\abs{c_x} (s-1)^2}+  \frac{A_x - 2B_x \log\abs{c_x}}{\abs{c_x}(s-1)} \Big) + O(1) & (s\to 1), \\
      & \E(\phi, x, s) = \frac{B_\infty}{s^2} - \frac{A_\infty}{s} + O(1) & (s\to 0).
  \end{align*}
\end{proof}

\section{Quantum modularity for the Eichler integral}\label{sec:qmod-Eichler}

For all~$x\in C(\Gamma)$, define
$$ \E_\phi(x) := \E(\phi, x, 1/2). $$

\begin{thm}\label{th:hgamma=qmf}
  The map $\E_\phi$ is a quantum modular form for $\Gamma$ with multiplier $u_\gamma$, in the sense that for all~$\gamma\in \Gamma$, the map
  \begin{equation}
    h^\E_\gamma(x) := \E_\phi(\gamma x) - j_\gamma(x)^k \chi(\gamma) \E_\phi(x),\label{eq:def-hP-E}
  \end{equation}
  initially defined for $x\in C(\Gamma)\setminus \{\infty, \gamma^{-1}\infty\}$ extends to a $(1/2-\eps)$-H\"{o}lder continuous function on $\R\setminus \{\gamma^{-1}\infty\}$. More precisely, for~$x, x' \not \in [\gamma^{-1}\infty-\eps, \gamma^{-1}\infty+\eps]$, we have
  \begin{equation}
    \abs{h^\E_\gamma(x) - h^\E_\gamma(x')} \ll_{\eps, \phi, \gamma} \abs{x-x'}^{1/2-\eps} (1+\abs{x}+\abs{x'})^{O_\phi(1)}.\label{eq:holder-bound-hE}
  \end{equation}

  If additionally~$\phi$ is cuspidal, then~$h_\gamma^\E$ is actually~$(1-\eps)$-Hölder continuous, and the bound~\eqref{eq:holder-bound-hE} holds with exponent~$1-\eps$ instead of~$1/2-\eps$.
\end{thm}

Note that for~$\gamma$ and~$x$ fixed, $j_\gamma(x)^k$ depends only on the parity of~$k$.
Note also that~$h_\gamma^\E = 0$ for all~$\gamma \in \Gamma_\infty$, so that Theorem~\ref{th:hgamma=qmf} is trivial in this case.

Throughout the rest of this section, we let~$\gamma\in \Gamma\setminus\Gamma_\infty$,
$$ x_0 = \gamma^{-1}\infty, $$
and we let~$I \subset \R\setminus\{x_0\}$ be a closed interval, not necessarily bounded.

\subsection{A geometric proof of quantum modularity}

We will start by considering the special case where $\phi=iR_0 \varphi$ is of weight $2$, cuspidal, with trivial nebentypus, meaning that $\varphi$ is a Maa{\ss} cusp form of weight $0$ for $\Gamma$ and $R_0$ is the weight $0$ raising operator. In this case one can give a pleasant geometric interpretation of the discrepancy function $h_\gamma=h_\gamma^\E$. This argument should furthermore serve to give some intuition before reading the somewhat technical proof in the general case, in Section~\ref{sec:E-qmf-gencase} below.

The starting point is the following alternative representation of the Eichler integral in the case of weight~$2$:
$$ \E_\phi(x)=  2i \int^{\infty}_{x} \frac{\partial}{\partial z} \varphi(z)dz, $$
which follows directly from the definition of the raising operator $R_0$.
By a change of variable $z\leftrightarrow \gamma z$ we get the following expression for the discrepancy function
$$h_\gamma(x)=2i\int^{x_0}_{x} \frac{\partial}{\partial z} \varphi(z)dz-2i\int_{x}^{ \infty} \frac{\partial}{\partial z} \varphi(z)dz$$
for $\gamma \notin \Gamma_\infty$ (the stabilizer of $\infty$) and $x\in C(\Gamma)\setminus \{\infty, x_0 \}$. The key idea is now that one can apply Stoke's Theorem to express the difference of these two line integrals (over infinite geodesics) to a surface integral over a surface of finite hyperbolic volume. More precisely,  for $x,x_0\in \R$ we denote by $\mathcal{F}_{x,x_0}$ the hyperbolic triangle 
with vertices $\infty,x,x_0$. Recall that the geodesic between two points on the boundary is exactly a Euclidean semi-circle through the two endpoint (if one of the points is $\infty$ this is vertical line), see Figure \ref{triangle}. We apply \cite[Lemma 2]{DukeImamogluToth} to the 1-form $\frac{\partial}{\partial z}\varphi(z)dz$ on the hyperbolic surface $\mathcal{F}_{ \gamma x,\gamma\infty}$ observing that \cite[Lemma 2]{DukeImamogluToth} easily extends to general discrete and cofinite subgroups $\Gamma$ as our 1-form is sufficiently regular at the cusps of $\Gamma$. This gives
$$  \int^{x_0}_{x} i\frac{\partial}{\partial z} \varphi(z)dz+\int^{x}_{ \infty} i\frac{\partial}{\partial z} \varphi(z)dz+\int_{x_0}^{\infty} i \frac{\partial}{\partial z} \varphi(z)dz=\frac{\lambda_\varphi}{2}\int_{\mathcal{F}_{ x,x_0}}\varphi(z)d\mu_0(z),  $$
where $d\mu(x+iy)=\frac{dxdy}{y^2}$ is the hyperbolic measure and~$\lambda_\varphi = \lambda_\phi$ is the Laplace eigenvalue.
This yields the following geometric expression
\begin{equation}\label{eq:hgammaspecialcase} h_\gamma(x)=-2i\int_{x_0}^{\infty} \frac{\partial}{\partial z} \varphi(z)dz+  \lambda_\varphi \int_{\mathcal{F}_{ x,x_0}}\varphi(z)d\mu_0(z).\end{equation}
Notice that the first term is independent of $x$. Now we see that for $x,x'\in I$ (i.e. bounded away from $x_0=\gamma^{-1}\infty $)
\begin{align*}
  h_\gamma(x)-h_\gamma (x')=\lambda_\varphi\left(\int_{\mathcal{F}_{x,x_0}} \varphi(z) d\mu_0(z)-\int_{\mathcal{F}_{x',x_0}} \varphi(z) d\mu_0(z)\right)\leq \abs{\lambda_\varphi} \|\varphi\|_\infty \cdot\mathrm{area}(\mathcal{F}_{ x,x_0}\triangle \mathcal{F}_{x',x_0}),
\end{align*}
where $A\triangle B$ denotes the symmetric difference between the sets $A$ and $B$, $\mathrm{area}$ denotes the hyperbolic area, and $|\!|\varphi|\!|_\infty$ is the sup norm of $\varphi$ (which is finite by cuspidality). The key point is now the following elementary hyperbolic geometric estimate, see Figure \ref{complement}.
\begin{lemma}
  For $x,x'\in I$, we have 
  $$
  \mathrm{area}(\mathcal{F}_{x,x_0}\triangle\mathcal{F}_{x',x_0})\ll_{\gamma, I} |x-x'|^{1/2}.
  $$
  More precisely, for~$x, x' \not \in [x_0-\eps, x_0+\eps]$, we have
  \begin{equation}\label{eq:geometricbound}
    \mathrm{area}(\mathcal{F}_{x,x_0}\triangle\mathcal{F}_{x',x_0}) \ll_{\eps, \gamma} \abs{x-x'}^{1/2} (1+\abs{x}+\abs{x'})^{-1/2}.
  \end{equation}  
\end{lemma}

\begin{proof} 
  First of all we observe that by the Gau{\ss} Defect  Theorem \cite[Theorem 1.3]{Iw} we have $\mathrm{area}(\mathcal{F}_{x,x_0})=\mathrm{area}(\mathcal{F}_{x',x_0})=\pi$. The symmetric difference
  $\mathcal{F}_{x,x_0}\triangle \mathcal{F}_{x',x_0}$  is the union of the two sets 
  $$A_1:= \mathcal{F}_{x',x_0}\setminus \mathcal{F}_{x,x_0}\cap \mathcal{F}_{x',x_0},\quad A_2:=\mathcal{F}_{x,x_0}\setminus \mathcal{F}_{x',x_0}\cap \mathcal{F}_{x,x_0},$$
  as illustrated in Figure \ref{complement} (and similar for the other configurations of $x,x',x_0$). Clearly we may assume that $x_0$ is not between~$x$ and~$x'$, for otherwise~$\abs{x-x'} \gg_\eps 1$ and the claimed bound is trivial. Thus we can restrict to the case $x'\geq x> x_0+\eps$. The three angles of the hyperbolic triangles $A_i$ are $0,0,\pi-\theta$ with  $0\leq \theta=\theta_{x,x'}\leq \tfrac{\pi}{2}$ such that $\cos \theta=\frac{2|x-x_0|}{|x'-x_0|}-1$, and thus by the Gau{\ss} Defect Theorem
  $$\mathrm{area}(A_1)=\mathrm{area}(A_2)=\theta. $$ 
  By the classical fact that $\tfrac{\sin x}{x} \in [\,\tfrac{2}{\pi} , 1) $ for $0\leq x\leq \tfrac{\pi}{2}$, we conclude 
  $$ \frac{4\theta^2}{\pi^2}\leq 1-(\cos \theta)^2= (\sin \theta)^2\leq \theta^2. $$
  This implies by the definition of $\theta$ and since $x'>x>x_0$ that
  $$ \theta^2
  \ll \frac{ \left||x'-x_0|-|x-x_0|\right| |x-x_0|}{|x'-x_0|^2}
  \ll \frac{|x-x'| |x'-x_0|}{|x-x_0|^2}\ll \frac{|x-x'| }{|x'-x_0|}.$$
  Now~\eqref{eq:geometricbound} follows by taking square roots since $x'-x_0\gg_{\eps,\gamma} 1+|x|+|x'|$.
\end{proof}

Thus we conclude that for~$x, x' \not \in [x_0-\eps, x_0+\eps]$, we have   
$$ h_\gamma(x)-h_\gamma(x_0)\ll_{\eps,\gamma} \abs{\lambda_\varphi} \|\varphi\|_\infty |x-x'|^{1/2} (1+|x|+|x'|)^{-1/2},$$
which is a more precise version of Theorem~\ref{th:hgamma=qmf} in this special case.

We end this section with a philosophical remark. As noted the vertical geodesics from $x$ to $i\infty$ have infinite hyperbolic length, which is responsible for the fact that the Eichler integral~$\E_\phi(x)$ itself is not continuous in $x$. However when considering the discrepancy $h_\gamma^\E$ we can transform this using Stoke's Theorem into an integral over a hyperbolic triangle bounded by such infinite geodesics, which has \emph{finite} hyperbolic area. This is a geometrical version of the idea of {\lq\lq}going up in regularity{\rq\rq} underlying the notion of quantum modular forms~\cite{Zagier2010}.

\begin{figure}[t]\centering
  \begin{tikzpicture}
    \begin{scope}
      \clip (-2,0) rectangle (2,3);
      \draw (0,0) circle(1.5);
      \draw (-2,0) -- (2,0);
      \draw (-1.5,0) --(-1.5,3);
      \draw (1.5,0) --(1.5,3);
      \fill[black, opacity=0.2] (-1.5,0) arc[start angle = 180, end angle = 0, radius = 1.5] --(1.5,3)--(-1.5,3)--cycle;
    \end{scope}
    \node at (0,2) [above] {$\mathcal{F}_{x,x_0}$};
    \node at (-1.5,0) [below] {$x$};
    \node at (1.5,0) [below] {$x_0$};

  \end{tikzpicture}
  \caption{Hyperbolic triangle with vertices $x,x_0,\infty$}\label{triangle}
  \begin{tikzpicture}
    \begin{scope}
      \clip (2,0) rectangle (-2,4);
      \draw (0,0) circle(1.5);
      \draw (-0.2,0) circle(1.3);
      \draw (2,0) -- (-2,0);
      \draw (-1.5,0) --(-1.5,4);
      \draw (1.5,0) --(1.5,4);
      \draw (1.1,0) --(1.1,4);
      \fill[black, opacity=0.2] (1.5,0) arc[start angle = 0, end angle = (42.8), radius = 1.5] --(1.1,4)--(1.5,4)--cycle;
      \fill[black, opacity=0.2] (1.1,0) arc[start angle = 0, end angle = 180, radius = 1.3] arc[start angle = 180, end angle = (42.8), radius = 1.5] --cycle; 
    \end{scope}
    \node at (1.5,0) [below] {$x'$};
    \node at (1.1,-0.121) [below] {$x$};
    \node at (-1.5,-0.121) [below] {$x_0$};
    \node at (0,1.5) [above] {$A_2$};
    \node at (1.1,3) [left] {$A_1$};
  \end{tikzpicture}
  \caption{Symmetric difference of $\mathcal{F}_{x,x_0}$ and $\mathcal{F}_{x',x_0}$}\label{complement}
\end{figure}

\subsection{Quantum modularity in the general case}\label{sec:E-qmf-gencase}

In this section we will prove the Hölder continuity of~$h_\gamma^\E$ on~$I$.
Let~$\delta \in (0, \frac12 d(x_0, I))$, where~$d(x_0, I)$ denotes the distance from~$x_0$ to~$I$.
For all~$x\in I$ and~$0<y\leq \delta$, we define
\begin{equation}
  \label{eq:def-mu-nu}
  \begin{aligned}
    \mu_x(y) := {}& \big(1 - \tfrac{y^2}{((x-x_0)/2)^2}\big)^{-1/2}, \\
    \nu_x(y) := {}& \tfrac{x-x_0}2 - \operatorname{sign}(x-x_0)\sqrt{(\tfrac{|x-x_0|}2)^2-y^2}.
  \end{aligned}
\end{equation}
With this definition, the points~$x_0+\nu_x(y)+iy$ and~$x-\nu_x(y)+iy$ both lie on the geodesic half-circle connecting~$x_0$ and~$x$.
We also recall the definition of~$u_\gamma$ in~\eqref{eq:def-ugamma}

The following simple bounds will be used repeatedly: for~$x, x'\in I$, $\delta\geq \eps$ and~$0<y\leq \delta$, we have
\begin{align}
  |\mu_x(y) - 1| &{}\ll_\eps y^2, & |\mu_x(y) - \mu_{x'}(y)| &{}\ll_\eps |x-x'| y^2, \label{eq:bounds-mu} \\
  |\nu_x(y)| &{}\ll_\eps y^2,  &|\nu_x(y) - \nu_{x'}(y)| &{}\ll_\eps |x-x'|y^2. \label{eq:bounds-nu}
\end{align}

The following lemma gives a good bound for an integral near a cusp which, as we will show later, corresponds to the least regular contribution to~$h_\gamma^\E$.
\begin{lemma}\label{lem:intphicusp-holder}
  For~$x \in I$ and~$0<y\leq \delta$, we have
  \begin{equation}
    |u_\gamma(x) \phi(x+iy) - u_\gamma(x-\nu_x(y)+iy) \phi(x-\nu_x(y)+iy) \mu_x(y)| = O_{\eps, \phi}(y^{1/2-\eps}),\label{eq:bound-phiint-diff}
  \end{equation}
  and the map
  \begin{equation}
    x\mapsto H(x) := \int_0^\delta \big(u_\gamma(x) \phi(x+iy) - u_\gamma(x-\nu_x(y)+iy) \phi(x-\nu_x(y)+iy) \mu_x(y)\big)\frac{dy}y\label{eq:int-diffE-cusppart}
  \end{equation}
  defines a~$(1/2-\eps)$-Hölder continuous function on~$I$. More precisely, the bound
  \begin{equation}\label{eq:bd-diff-H}
    H(x) - H(x') = O_{I, \eps, \phi, \gamma}\big(\abs{x-x'}^{1/2-\eps} + \abs{x-x'}\big)
  \end{equation}
  holds uniformly for~$x, x' \in I$.

  If~$\phi$ is cuspidal, then the improved bound
  \begin{equation}\label{eq:bd-diff-H-cusp}
    H(x) - H(x') = O_{I, \eps, \phi, \gamma} \big(\abs{x-x'}^{1-\eps} + \abs{x-x'}\big) \qquad (x, x' \in I)
  \end{equation}
  holds.
\end{lemma}
\begin{proof}
  In the following proof, we allow all implicit constants to depend on~$I, \eps, \phi$ and~$\gamma$.
  Note first that~$u_\gamma$ is constant in a some neighborhood~$I^+$ of~$I$, say~$u_\gamma(x) = u_\gamma^0 \in \R$ for all~$x\in I$. Moreover, or any~$x\in I$ and~$0<y\ll 1$, we have
  \begin{equation}\label{eq:ugamma-diff}
    \abs{u_\gamma(x+iy)-u_\gamma^0} \ll y.
  \end{equation}
  We complement this with the bounds~\eqref{eq:bound-phi}, \eqref{eq:bound-phi-diff}, \eqref{eq:bounds-mu} and \eqref{eq:bounds-nu}. Along with the triangle inequality, this yields the bounds~\eqref{eq:bound-phiint-diff}.
  
  We deduce that the integrand in~\eqref{eq:int-diffE-cusppart} is~$O(y^{-1/2-\eps})$ and the integral in~\eqref{eq:int-diffE-cusppart} is well-defined.
  Let~$\delta'\in (0, \delta)$ be such that~$x + \nu_x(y) \in I^+$ for all~$(x, y) \in I \times (0, \delta')$, which implies that~$u_\gamma(x+\nu_x(y)) = u_\gamma^0$ for~$0<y<\delta'$.
  In the following computations, let us abbreviate
  $$ \nu = \nu_x(y),\quad \nu' = \nu_{x'}(y), \quad \mu = \mu_x(y), \quad \mu' = \mu_{x'}(y). $$
  By the triangle inequality, for~$x, x' \in I$, we have
  $$ \abs{H(x)-H(x')} \leq \int_0^\delta \abs{F_0(y; x, x')} \frac{\df y}{y}, $$
  where
  $$ F_0(y; x, x') := u_\gamma^0(\phi(x+iy) - \phi(x'+iy)) - u_\gamma(x-\nu+iy)\phi(x-\nu+iy)\mu + u_\gamma(x'-\nu'+iy)\phi(x'-\nu'+iy)\mu'. $$
  By taking successive differences, we write
  $$ F_0(y; x, x') = \sum_{1\leq j \leq 6} F_j(y; x, x'), $$
  where
  \begin{align*}
    F_1(y; x, x') := {}& u_\gamma^0(\phi(x+iy) - \phi(x'+iy) - \phi(x-\nu+iy)+\phi(x'-\nu+iy)), \\
    F_2(y; x, x') := {}& u_\gamma^0(\phi(x'-\nu'+iy)-\phi(x'-\nu+iy)), \\
    F_3(y; x, x') := {}& (u_\gamma(x'-\nu'+iy)-u_\gamma^0)(\phi(x'-\nu'+iy) - \phi(x-\nu+iy)), \\
    F_4(y; x, x') := {}& u_\gamma(x'-\nu'+iy)(\phi(x'-\nu'+iy) - \phi(x-\nu+iy))(1-\mu'), \\
    F_5(y; x, x') := {}& u_\gamma(x'-\nu'+iy)\phi(x-\nu+iy)(\mu'-\mu), \\
    F_6(y; x, x') := {}& (u_\gamma(x'-\nu'+iy)-u_\gamma(x-\nu+iy))\phi(x-\nu+iy) \mu. \\
  \end{align*}
  We bound each integral
  $$ D_j := \int_0^\delta \abs{F_j(y; x, x')}\frac{\df y}y $$
  separately.
  \begin{itemize}
  \item By~\eqref{eq:bound-phi-doublediff} and the first bound in~\eqref{eq:bounds-nu}, we have~$F_1(y; x, x') \ll y^{1/2-\eps}\min(1, \frac{|x-x'|}y)$, and so
    $$ D_1 \ll \abs{x-x'}^{1/2-\eps}. $$
  \item By~\eqref{eq:bound-phi-diff} and the second bound in~\eqref{eq:bounds-nu}, we obtain~$F_2(y; x, x') \ll y^{1/2-\eps} \abs{x-x'}$, so that
    $$ D_2 \ll \abs{x-x'}. $$
  \item By~\eqref{eq:bound-phi-diff} again, and~\eqref{eq:ugamma-diff}, we obtain~$F_3(y; x, x')\ll y^{-1/2-\eps} \min(1, \frac{|x-x'|}y)$, and so
    $$ D_3 \ll \abs{x-x'}^{1/2-\eps}. $$
  \item By~\eqref{eq:bound-phi-diff}, the second bound in~\eqref{eq:bounds-nu} and the first bound in~\eqref{eq:bounds-mu}, we have~$F_4(y; x, x')\ll y^{1/2-\eps}\abs{x-x'}$, therefore
    $$ D_4 \ll \abs{x-x'}. $$
  \item By~\eqref{eq:bound-phi} and the second bound in~\eqref{eq:bounds-mu}, we get~$F_5(y; x, x') \ll y^{3/2-\eps} \abs{x-x'}$, so that
    $$ D_5 \ll \abs{x-x'}. $$
  \item Finally, for all~$a, a'\in I^+$, we have
    $$ \abs{u_\gamma(a+iy) - u_\gamma(a'+iy)} \ll \abs{\arg(a+iy) - \arg(a'+iy)} \ll y\abs{a-a'}, $$
    so that, by~\eqref{eq:bound-phi}, the second bound in~\eqref{eq:bounds-nu} and the bound~$\mu\ll 1$, we obtain~$F_6(y; x, x') \ll y^{1/2-\eps}\abs{x-x'}$, and so
    $$ D_6 \ll \abs{x-x'}. $$
  \end{itemize}
  When~$\phi$ is cuspidal, by Lemma~\ref{lem:bounds-phi}, the above bounds on~$D_1$ and~$D_3$ hold with exponent~$1-\eps$ instead of~$1/2-\eps$.

  We conclude that
  $$ \abs{H(x) - H(x')} \leq \sum_{1\leq j \leq 6} D_j \ll \begin{cases} \abs{x-x'}^{1/2-\eps} + \abs{x-x'} & (\text{in general}), \\ \abs{x-x'}^{1-\eps} + \abs{x-x'} & (\phi\text{ cuspidal}), \end{cases} $$
  as claimed.
\end{proof}

\begin{proof}[Proof of Theorem~\ref{th:hgamma=qmf}]
  For~$\Re s>1$ and~$x\in I \cap C(\Gamma)$, we recall the definition~\eqref{eq:def-eichlerint}, and we abbreviate in this proof~$\E(x, s) = \E(\phi, x, s)$.
  We consider
  \begin{equation}
    \Delta(x, s) := u_\gamma^0 \E(x, s) - \abs{j(\gamma, x)}^{2s-1} \E(\gamma x, s).\label{eq:def-Delta-E}
  \end{equation}
  We recall that~$u_\gamma(x) = u_\gamma^0$ for all~$x\in I$.
  Since~$\delta < \frac12 d(x_0, I)$, for each~$x\in I$, the line~$\Im(z) = \delta$ intersects the geodesic connecting~$x$ and~$x_0$ at two distinct points.
  For all~$x\in I$, we let~$\eta(x)>0$ be the smaller solution to
  \begin{equation}
    \Im( \gamma^{-1}(\gamma x + i \eta(x))) = \delta. \label{eq:def-eta-delta}
  \end{equation}
  Note that~$\eta$ defines a smooth map on~$I$, which is bounded and non-zero. We have explicitely
  $$ \eta(x) = \frac{\abs{\nu_x(\delta)}}{c^2 \delta} = \frac1{2c^2 \delta}\Big(1 - \sqrt{1-(2\delta/\abs{x-x_0})^2}\Big). $$
  We deduce in particular
  $$ \abs{\eta(x)}, \abs{\eta'(x)} \asymp_\eps (x-x_0)^{-2}, \qquad (x \in I). $$
  By inserting the integral~\eqref{eq:def-eichlerint} in the definition of~$\Delta(x, s)$, and splitting the two integrals respectively at~$y=\delta$ or~$y=\eta(x)$, we obtain for~$x\in I \cap C(\Gamma)$ the decomposition
  \begin{equation}
    \Delta(x, s) = \Delta_0(x, s) - \Delta_c(x, s) + \Delta_\infty(x, s),\label{eq:Delta-split}
  \end{equation}
  where
  \begin{align*}
    \Delta_0(x, s) := {}& u_\gamma^0 \int_0^\delta \phi(x+iy) y^{s-1/2} \frac{\df y}y - \abs{j(\gamma, x)}^{2s-1} \int_0^{\eta(x)} \phi(\gamma x + iy) y^{s-1/2} \frac{\df y}y, \\
    \Delta_c(x, s) := {}& u_\gamma^0 \int_0^\delta \phi_\infty(y) y^{s-1/2} \frac{\df y}y - \abs{j(\gamma, x)}^{2s-1} \int_0^{\eta(x)} \phi_\infty(y) y^{s-1/2} \frac{\df y}y, \\
    \Delta_\infty(x, s) := {}& u_\gamma^0 \int_\delta^\infty (\phi(x+iy) - \phi_\infty(y)) y^{s-1/2}\frac{\df y}y - \abs{j(\gamma, x)}^{2s-1} \int_{\eta(x)}^\infty (\phi(\gamma x + iy) - \phi_\infty(y)) y^{s-1/2} \frac{\df y}y.
  \end{align*}
  By the explicit expression~\eqref{eq:def-phi-cuspidalpart}, the map~$\Delta_c(x, \cdot)$ can be analytically continued to~$\C\setminus\{1/2-s_\phi, s_\phi-1/2\}$. For some constants~$c_1, c_2, c_3$ depending on~$\phi$, $\gamma$ and~$I$, we evaluate for all~$x\in I\cap C(\Gamma)$,
  $$ \Delta_c(x, 1/2) = c_1 + c_2 \eta(x)^{s_\phi} + c_3 \eta(x)^{1-s_\phi}. $$
  Since~$\eta$ is smooth and non-zero on~$I$, we deduce that the map~$\Delta_c(\cdot, 1/2)$ extends to a smooth function on~$I$, and
  \begin{equation}\label{eq:bd-diff-Deltac}
    \abs{\Delta_c(x, 1/2) - \Delta_c(x', 1/2)} \ll_{I, \gamma, \phi, \eps} \abs{x-x'} (1+\abs{x}+\abs{x'})^{\kappa+\eps}
  \end{equation}
  with~$\kappa = \max\{0, 2(\Re(s_\phi)-1)\}\in\R_{\geq 0}$.

  The integrals in~$\Delta_\infty$ are uniformly convergent for bounded~$s\in \C$, since the maps~$y\mapsto \phi(t+iy)-\phi_\infty(y)$, for~$t\in \{x, \gamma x\}$, have exponential decay at~$\infty$. Since the Fourier expansion~\eqref{fexp} is uniformly convergent for~$\Im(z)>\eps$ for any~$\eps>0$, and~$\eta$ is non-zero on~$I$, we deduce by dominated convergence that~$\Delta_\infty(\cdot, 1/2)$ also extends to a smooth function on~$I$. Moreover, letting~$\phi_1(x+iy) := \frac{\partial}{\partial x} \phi(x+iy)$, we have~$\phi_1(x+iy) \ll_{\phi, \eps} y^{-3/2+\eps}$, and thus
  \begin{align*}
    \frac{d}{dx} \Delta_\infty(x, 1/2) = {}& u_{\gamma}^0 \int_\delta^\infty \phi_1(x+iy) \frac{\df y}y - j(\gamma, x)^{-2} \int_{\eta(x)}^\infty \phi_1(\gamma x+iy) \frac{\df y}{y} + \frac{\eta'(x)}{\eta(x)} \phi(\gamma x + i\eta(x)) \\
    \ll_{I, \gamma, \phi, \eps} {}& \abs{x-x_0} \ll 1+\abs{x}. \numberthis\label{eq:bd-diff-Deltainf}
  \end{align*}
  We obtain by integration
  $$ \abs{\Delta_\infty(x, 1/2) - \Delta_\infty(x', 1/2)} \ll_{I, \gamma, \phi, \eps} \abs{x-x'} (1+\abs{x}+\abs{x'}). $$

  We focus on~$\Delta_0(x, s)$. For~$\Re(s)>1$, we change variables in the second integral, getting
  \begin{align*}
    \Delta_0(x, s) ={}& u_\gamma^0 \int_x^{x+i\delta} \phi(z) (\Im z)^{s-1/2} \df s(z) - \abs{j(\gamma, x)}^{2s-1} \int_{\gamma x}^{\gamma x + i\eta(x)} \phi(z) (\Im z)^{s-1/2} \df s(z) \\
    = {}& u_\gamma^0 \int_x^{x+i\delta} \phi(z) (\Im z)^{s-1/2} \df s(z) - \abs{j(\gamma, x)}^{2s-1} \int_x^{\gamma^{-1}(\gamma x + i\eta(x))} \phi(\gamma z) (\Im \gamma z)^{s-1/2} \df s(z) \\
    = {}& \int_x^{x+i\delta} u_\gamma(x) \phi(z) (\Im z)^{s-1/2} \df s(z) - \int_x^{\gamma^{-1}(\gamma x + i\eta(x))} u_\gamma(z) \phi(z) (\Im z)^{s-1/2} \abs{\frac{j(\gamma, x)}{j(\gamma, z)}}^{2s-1} \df s(z). \numberthis \label{eq:Delta0-varchange-sz}
  \end{align*}
  Here the integrals with respect to the~$\SL(2, \R)$-invariant Poincaré metric~$\df s(z) = (\df x^2 + \df y^2)^{1/2}/y$ are taken along geodesics. In the penultimate integral, we have written~$u_\gamma(x) = u_\gamma^0$ in anticipation of using of Lemma~\ref{lem:intphicusp-holder}.
  We parametrize both integrals according to~$\Im z$, which runs in both cases over~$(0, \delta)$ by our definition~\eqref{eq:def-eta-delta}. The integral~$\int_x^{x+i\delta}$ is over a straight line, and the second integral~$\int_x^{\gamma^{-1}(\gamma x + i\eta(x))}$ is over a portion of geodesic connecting~$x$ with~$\gamma^{-1}\infty = x_0$. By construction of~$\nu_x$ in~\eqref{eq:def-mu-nu}, this is precisely the set
  $$ \{ x - \nu_x(y) + iy, y\in (0, \delta) \}. $$
  Moreover, a quick computation yields~$\df s(z) = \mu_x(y) \df y / y$. We deduce
  $$ \Delta_0(x, s) = \int_0^\delta \Big(u_\gamma(x) \phi(x+iy) - u_\gamma(x-\nu_x(y)+iy) \phi(x-\nu_x(y)+iy) \mu_x(y) \abs{\frac{j(\gamma, x)}{j(\gamma, z)}}^{2s-1}\Big) y^{s-1/2} \frac{\df y}y. $$
  We write~$\abs{j(\gamma, x)/j(\gamma, z)}^{2s-1} = 1 + \xi_{s, x}(y)$, where
  \begin{equation}\label{eq:bound-jgammax-s}
    \xi_{s, x}(y) \ll y
  \end{equation}
  for fixed~$x, \gamma$ and bounded~$s$. Correspondingly, for all~$x\in I\cap C(\Gamma)$ we obtain
  \begin{align*}
    \Delta_0(x, s) ={}& \int_0^\delta (u_\gamma(x) \phi(x+iy) - u_\gamma(x-\nu_x(y)+iy) \phi(x-\nu_x(y)+iy) \mu_x(y)) y^{s-1/2} \frac{\df y}y \\
    &{} - \int_0^\delta u_\gamma(x-\nu_x(y)+iy) \phi(x-\nu_x(y)+iy) \mu_x(y) \xi_{s,x}(y) y^{s-1/2} \frac{\df y}y.
  \end{align*}
  The bound~\eqref{eq:bound-jgammax-s}, combined with~\eqref{eq:bound-phi}, ensures that the second integral here converges uniformly on compact subsets of~$\{\Re(s)>0\}$. The bound~\eqref{eq:bound-phiint-diff} from Lemma~\ref{lem:intphicusp-holder} yields the same conclusion for the first integral. This gives the analytic continuation of~$\Delta_0(x, \cdot)$ on~$\{\Re(s)>0\}$, and the expression
  $$ \Delta_0(x, 1/2) = \int_0^\delta (u_\gamma(x) \phi(x+iy) - u_\gamma(x-\nu_x(y)+iy) \phi(x-\nu_x(y)+iy) \mu_x(y)) \frac{\df y}y, $$
  which we recognize as the quantity~$H(x)$ defined in~\eqref{eq:int-diffE-cusppart}. Lemma~\ref{lem:intphicusp-holder} shows that~$\Delta_0(\cdot, 1/2)$ extends to a function on~$I$ which is~$(1/2-\eps)$-Hölder continuous for any~$\eps>0$, in general, and~$(1-\eps)$-Hölder continuous if~$\phi$ is cuspidal. Since we had shown earlier that~$\Delta_\infty(\cdot, 1/2)$ and~$\Delta_c(\cdot, 1/2)$ are smooth on~$I$, we deduce by~\eqref{eq:Delta-split} that~$\Delta(\cdot, 1/2)$ extends to a~$(1/2-\eps)$-Hölder continuous function on~$I$, and~$(1-\eps)$-Hölder if~$\phi$ is cuspidal. Moreover, the bounds~\eqref{eq:bd-diff-Deltac}, \eqref{eq:bd-diff-Deltainf} and~\eqref{eq:bd-diff-H}, we get
  $$ \abs{\Delta(x, 1/2) - \Delta(x', 1/2)} \ll_{I, \gamma, \phi, \eps} \abs{x-x'}^{1/2-\eps} + \abs{x-x'} (1+\abs{x}+\abs{x'})^{\max(1, 2(\Re(s_\phi)-1)) + \eps}. $$

  However, we see from the definition~\eqref{eq:def-Delta-E} and Proposition~\ref{prop:def-twisted-l} that for~$x\in I \cap C(\Gamma)$,
  $$ \Delta(x, 1/2) = u_\gamma(x) \E(x, 1/2) - \E(\gamma x, 1/2) = - h_\gamma^\E(x). $$
  This gives the desired extension of~$h_\gamma^\E$ to a Hölder continuous function on~$I$, and the claimed bound~\eqref{eq:holder-bound-hE}.
\end{proof}

\subsection{Generalized quantum modularity with conjugation}\label{sec:GQM}

For the applications to reciprocity formulas which we will consider in Section~\ref{sec:arith-applications} below in the case of the Hecke groups~$\Gamma = \Gamma_0(q)$, we will need the action of Fricke involutions, and in the case of non-real nebentypus, it will be handy to generalize slightly the definition of quantum modularity.
Given a matrix $\gamma\in \GL(2, \R)$ we consider the following operation on functions $f:\Hb\rightarrow \C$ 
\begin{equation}\label{eq:extendedQM}
  (\overline{\gamma}f)(z):= \overline{f( \gamma( -\overline{z}))}.
\end{equation}
We extend this  definition ``to the boundary'' as follows: let $\Gamma$ be a Fuchsian group as in our setting (Section~\ref{sec:setting}), and $\gamma\in \Gamma$ be such that we have $\gamma (-x)\in C(\Gamma)$ for all $x\in C(\Gamma)$. Then given a map $f:C(\Gamma)\backslash\{\infty\} \rightarrow \C$, we define 
$$ (\overline{\gamma}f)(x):= \overline{f(\gamma (-x))}, $$
for $x\in C(\Gamma)\backslash\{-\gamma^{-1}\infty\} $.
\begin{thm}\label{th:hgamma-fricke=qmf}
  Let $\gamma\in \SL(2, \R)$ and let $\phi\in \A(\Gamma, \chi, k)$ satisfy hypotheses~$(H_1)$-$(H_6)$ from Section~\ref{sec:setting}, and assume additionally that $(\overline{\gamma}\phi)(z)=\eta j_{\gamma}(z)^k\phi(z)$ for some $\eta\in \C$. Then $\mathcal{E}(\phi,x)$ is quantum modular for $\overline{\gamma}$ in the sense that
  \begin{equation}\label{eq:def-hgamma-E-ext}
    h^\mathcal{E}_{\overline{\gamma}}(x)=\overline{\mathcal{E}_\phi(\gamma (-x))}-\eta j_\gamma(x)^k\mathcal{E}_\phi(x),
  \end{equation}
  initially defined for $x\in C(\Gamma)\backslash\{ \infty, -\gamma^{-1}\infty\}$, extends to a $(1/2-\eps)$-H\"{o}lder continuous map $\R\backslash\{-\gamma^{-1}\infty\}\rightarrow \C$.
\end{thm}
\begin{proof}
  The proof is identical to that of Theorem~\ref{th:hgamma=qmf}, starting from the difference
  $$ \Delta(x, s) := \eta j_\gamma(x)^k \E(x, s) - \abs{j(\gamma, x)}^{2s-1} \overline{\E(\gamma(-x), \bar s)}. $$
\end{proof}

\begin{rem}
  The above notion of quantum modularity can be put into the following general framework. Consider a structured space $Y$  (e.g. a topological space, manifold), a group $G$ and a representation $\rho:G\rightarrow \mathrm{Aut}(\C[Y])$ where $\C[Y]$ denotes the ring of (set theoretic) maps $Y\rightarrow \C$ with pointwise addition and multiplication. Let $u:G\rightarrow  \C[Y]$ be a cocycle for the pair $(G,\rho)$ in the sense that
  $$\rho_g(u(g',\cdot))\cdot u(g,\cdot)= u(gg',\cdot), $$
  for all $g,g'\in G$. Then we say that $f\in \C[Y]$ is \emph{modular} for the pair $(G, \rho)$ with multiplier $u$ if 
  $$ \rho_g(f)(y)=u(g,y)f(y), $$
  for all $g\in G$ and $y\in Y$. Notice that the cocycle condition is forced from the previous equation since $\rho$ preserves the ring structure. Furthermore, let $Y_0$ be some (possibly discrete) subset of $Y$. Then we say that a map $f\in \C[Y_0]$ is \emph{quantum modular} for $(G,\rho)$ with multiplier $u_g$ if
  $$ h_g(y):=\rho_g(f)(y)-u(g,y)f(y),  $$
  is {\lq\lq}more regular{\rq\rq} than $f$ itself, i.e. $h_g$ can be extended to a continuous (or smooth, or analytic) function on some intermediate space $Y_0\subset Y'\subset Y$.

  Let $\overline{\Hb}=\Hb\cup \R \cup \{\infty\}$ and consider the group  $G=\mathrm{PGL}_2(\R)\times \Z/2\Z$ acting on $\C[\overline{\Hb}]$ where $\PSL_2(\R)$ acts by precomposition with the inverse of the associated linear fractional transformation (in order to insure that it is a left action), $\begin{psmallmatrix} -1 & 0 \\ 0 & 1\end{psmallmatrix}$ acts by precomposition with $z\mapsto -\overline{z}$ and the generator $c$ of $\Z/2\Z$ acts by composition with conjugation. One can check that this is indeed a well-defined group action. Now let $G_q\subset G$ be the subgroup generated by $\Gamma_0(q)$ together with the element 
  $$(\begin{psmallmatrix} 0 & 1 \\ q & 0\end{psmallmatrix}, c)\in G,$$
  which acts as the operator $\overline{W_q}$  using the notation~\eqref{eq:extendedQM} with $W_q=\begin{psmallmatrix} 0 & -1 \\ q & 0\end{psmallmatrix}$ the Fricke matrix of level $q$. One can now check that the following defines a cocycle for $G_q$ with the representation $\rho: G_q \rightarrow \mathrm{End}_\mathrm{Ring}(\C[\overline{\Hb}])$ 
  \begin{align}\label{eq:multiplier}u(\overline{W_q},z):= \eta j_{W_q}(z)^k , \quad u(\gamma,z):= \chi(\gamma)j_{\gamma^{-1}} (z)^k,\, \text{for }\gamma\in \Gamma_0(q),\end{align} 
  where $\eta\in \C$ satisfies $|\eta|=1$  and $\chi:\Gamma_0(q)\rightarrow \C^\times$ is a character (here the only relation one has to check is $\overline{W_q} \gamma= \gamma' \overline{W_q}$ for $\gamma \in \SL_2(\R)$ and $\overline{W_q}^2=1$). In particular, we will see in Section \ref{sec:HeckeMaass} that if $\phi$ is a Hecke--Maa{\ss} newform of level $q$, weight $k$ and nebentypus $\chi_\phi$, then we automatically have that $\overline{W_q}\phi= \eta_\phi (j_{W_q})^k \phi$ for some $\eta_\phi\in \C$ of absolute norm $1$. Thus in the just described formalism, Theorem~\ref{th:hgamma=qmf} and Theorem~\ref{th:hgamma-fricke=qmf} can be interpreted as saying that $$\mathcal{E}(\phi,\cdot):C(\Gamma_0(q))\backslash \{\infty\} =\Q\rightarrow \C,$$ is quantum modular for the group $G_q$ together with its obvious representation where the multiplier  is as in \eqref{eq:multiplier} (with $\eta=\eta_\phi$ and $\chi=\chi_\phi$).

\end{rem}

\section{Behaviour of~\texorpdfstring{$h^\mathcal{E}_\gamma$}{h-E-gamma} at \texorpdfstring{$\infty$}{infinity}}\label{sec:hgamma-at-infinity}

For our prospective applications, it will be of importance to have some information of the behaviour of~$h_\gamma^\E$ near~$\infty$ and~$\gamma^{-1}\infty$.

\begin{thm}\label{th:hgamma-asymptotic}
  There exist numbers~$A'_\pm, B'_\pm$ depending on~$\phi$ such that the following holds.
  Let $\gamma\in \Gamma\setminus\Gamma_\infty$, and~$x_0 := \gamma^{-1}\infty$. If~$s_\phi\neq 1/2$, then
  \begin{equation}
    h_\gamma^\E (x)= \chi(\gamma) \big(A'_\pm \abs{x-x_0}^{s_\phi} + B'_\pm \abs{x-x_0}^{1-s_\phi} - i^k\E_\phi(x_0)\big) + O_{\phi, \gamma, \eps}(|x|^{-1+\eps}) \qquad \text{as} \quad x\to \pm \infty.\label{eq:hgamma-E-asymptotic}
  \end{equation}
  If~$s_\phi = 1/2$, the asymptotic formula holds upon replacing the term involving~$B'_\pm$ by~$B'_\pm \abs{x-x_0}^{1/2} \log \abs{x-x_0}$.

  The numbers~$A'_\pm, B'_\pm$ are obtained from the coefficients~$A, B$ in~\eqref{eq:def-phi-cuspidalpart} through
  \begin{equation}
    A'_\pm = A \Upsilon_{k, \pm 1}(s_\phi), \qquad B'_\pm = B \Upsilon_{k, \pm 1}(1-s_\phi), \qquad (s_\phi \neq 1/2).\label{eq:expr-A'-B'-not12}
  \end{equation}
  where~$\Upsilon_{k, \pm}(s_\phi)$ satisfy equation~\eqref{eq:expr-Upsilon} below. When~$s_\phi=1/2$, we have
  \begin{equation}
    A'_\pm = A \Upsilon_{k, \pm 1}(\tfrac12) + B \Upsilon_{k, \pm 1}'(\tfrac12), \qquad B'_\pm = B\Upsilon_{k, \pm 1}(\tfrac12), \qquad (s_\phi = 1/2).\label{eq:expr-A'-B'-12}
  \end{equation}
\end{thm}

We insist that in this statement, it is understood that the bottom left coefficient of~$\gamma$ is positive, or otherwise the ill-defined factor~$\chi(\gamma)$ must be replaced by~$\chi(\gamma) \sgn(c)^k$.

We deduce the following asymptotic expansion close to~$\gamma^{-1}\infty$.

\begin{cor}\label{cor:hgamma-asymptotic-finite}
  Let the notations be as in Theorem~\ref{th:hgamma-asymptotic}, and let~$c>0$ denote the bottom left coefficient of~$\gamma$. Then if~$s_\phi \neq 1/2$, we have
  $$ h_\gamma^\E(\gamma^{-1} \infty - c^{-2}\delta) = -(\mp 1)^k\big(A'_\pm \abs{\delta}^{-s_\phi} + B'_\pm \abs{\delta}^{s_\phi-1} - i^k \E(\gamma \infty)\big) + O_{\phi, \gamma, \eps}(\abs{\delta}^{1-\eps}) \qquad \text{as} \quad \delta \to 0^\pm. $$
  If~$s_\phi=1/2$, the analogous formula holds with the term involving~$B'_\pm$ replaced by~$B'_\pm \abs{\delta}^{-1/2}\log(1/\abs{\delta})$.
\end{cor}

\subsection{A geometric proof in a special case}

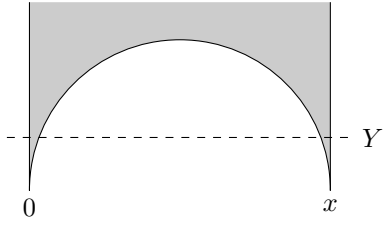
\begin{figure}[t]\centering
  \begin{tikzpicture}
    \begin{scope}
      \clip (-2.5,0) rectangle (2.5,2.5);
      \draw (0,0) circle(2);
      \draw (-2,0) -- (-2 ,2.5);
      \draw (2,0) --(2,2.5);
      \draw[dashed] (-2.3,0.7) --(2.3,0.7);
      \fill[black, opacity=0.2] (-2,0) arc[start angle = 180, end angle = 0, radius = 2] --(2,2.5)--(-2,2.5)--cycle;
    \end{scope}
    \node at (-2,0) [below] {$0$};
    \node at (2,0) [below] {$x$};
    \node at (2.3,0.7) [right] {$Y$};
    
  \end{tikzpicture}
  \caption{The hyperbolic triangle $\mathcal{F}_{x, 0}$ for large $x$}\label{highincusp}
\end{figure}

As for quantum modularity we will begin with a geometric proof in the special case where $\phi=i R_0 \varphi$ with $\varphi$ a Maa{\ss} cusp form of weight $0$, trivial multiplier and now assumed to be of level $1$ (\emph{i.e.} $\Gamma=\PSL(2,\Z)$). We will consider the behavior of $h_S=h_S^\mathcal{E}$ where $S=\begin{psmallmatrix} 0 & -1 \\ 1 & 0 \end{psmallmatrix}$. In this case the proof boils down to an elementary geometric lemma, which in words states that as $|x|\rightarrow \infty$, the majority of the area of $\mathcal{F}_{x,0}$ is near the cusp at $\infty$ (see Figure \ref{highincusp}).

\begin{lemma}\label{lem:smalltriangle}
  We have for any $Y\geq 1$: 
  $$ \mathrm{area}(\mathcal{F}_{x, 0} \cap \{z\in \Hb : \Im z\leq Y \})\leq 2 \min(Y/|x|, 1/2) .$$
\end{lemma}
\begin{proof}
  Clearly we can assume $x>0$. Now the result follows from the following calculation
  \begin{align*}
    \mathrm{area}(\mathcal{F}_{ x, 0}\cap\{z\in \Hb : \Im z\leq Y \}) &=2 \int_0^{\min(x/2,Y)} \int_0^{x-\sqrt{x^2-v^2}} \frac{dudv}{v^2}   \\
    &=  2 \int_0^{\min(x/2,Y)} \frac{1}{x+\sqrt{x^2-v^2}} dv \leq 2 \min(x/2,Y) x^{-1}.
  \end{align*} 
\end{proof}
Using the expression \eqref{eq:hgammaspecialcase} for $h_S$ in this case, we write
\begin{align*}
  h_S(x)=-2i\int_{0}^{\infty} \frac{\partial}{\partial z} \varphi(z)dz+  \lambda_\varphi\int_{\mathcal{F}_{ x, 0}\cap\{z\in \Hb : \Im z\leq Y \}} \varphi(z) d\mu_0(z)+\lambda_\varphi\int_{\mathcal{F}_{ x, 0}\cap\{z\in \Hb : \Im z> Y \}} \varphi(z) d\mu_0(z),
\end{align*}
with $Y=\frac{2}{\pi}\log |x|$. For $\Im z>\frac{2}{\pi}\log |x|$, we have 
$$ \varphi(z)\ll_\varphi e^{-\frac{\pi}{2}\Im z }\ll e^{-\log |x|}= |x|^{-1},$$ 
using the exponential decay of the Whittaker function. Now using that $\mathrm{area}(\mathcal{F}_{ x, 0})=\pi$ as well as Lemma \ref{lem:smalltriangle}, we conclude that
$$h_{S}(x) = -2i\int_{0}^{\infty} \frac{\partial}{\partial z} \varphi(z)dz+O_{\varphi}\left( \frac{\log |x|}{|x|}\right)$$ 
as $x\rightarrow \pm \infty$. This is exactly Theorem \ref{th:hgamma-asymptotic} in our case (even with an improved error-term). 

Finally we notice that by Fourier expanding it can easily be seen shown that 
\begin{align*}
  -2i\int_{0}^{\infty} \frac{\partial}{\partial z} \varphi(z)dz=\begin{cases} \gamma_\varphi(1/2) L(\varphi ,1/2),& \epsilon_\varphi=-1\\ 0,& \epsilon_\varphi=1,\end{cases},
\end{align*}
where $\gamma_\varphi(s) = 2^{3/2}(2\pi)^{-s} \Gamma\left(\frac{s+1+it_\varphi}{2}\right)\Gamma\left(\frac{s+1-it_\varphi}{2}\right)$, $L(\varphi, s) = \sum_{n>0} a(n) n^{1/2-s}$ is the standard $L$-series associated to $\varphi$, and $\epsilon_\varphi$ is the eigenvalue of $\varphi$ under the involution $Q_{s_\phi, 0}$. We will see in Section \ref{sec:relperiod} this fact in greater generality.

\subsection{Proof in the general case}

The rest of this section is devoted to the proof of Theorem~\ref{th:hgamma-asymptotic}. We allow all implicit constants to depend on~$\phi, \gamma$ and~$\eps$, possibly in addition to other parameters indicated in subscript.
Let~$x\in C(\Gamma)$ be given, and
$$ x_0 := \gamma^{-1}\infty, \qquad Y := \abs{x}^\eps. $$
We may assume that~$x\neq x_0$. It suffices to prove the theorem for~$x\in C(\Gamma)$, since we have shown that~$h_\gamma^\E$ is continuous.
We take a representative~$\gamma = \begin{psmallmatrix}a&b\\c&d\end{psmallmatrix}$ with~$c>0$.
For~$r, r' \in C(\Gamma)$, let~$G(r, r')$ denote the geodesic going from~$r$ to~$r'$ in~$\H$.
Let
$$ Z = G(x, x_0), \qquad L_r := G(r, \infty), $$
so that~$L_r$ is follows a vertical half-line upwards. Finally, we denote the truncations
$$ Z^{\pm}, \qquad L_r^{\pm} $$
where~$Z^+$ is the portion of~$Z$ located in~$\{\Im z \geq Y\}$, and~$Z^-$ consists of the two portions of~$Z$ located in~$\{\Im z \leq Y\}$.
Similarly, $L_r^+$ is the portion of~$L_r$ starting from~$r+iY$ upwards, and~$L_r^-$ is the segment~$[r, r+iY]$.
For~$\Re(s)>1$, we consider the integrals
\begin{align}
  \cG(x, s) := {}& \int_{Z} u_\gamma(z) \phi(z) (\Im z)^{s-1/2} \df s(z), \label{eq:def-Gxs} \\
  V_\phi(t) := {}& \int_{G(\pm 1, 0)} (u/\abs{u})^k \phi_\infty(\abs{t}\Im u) \df s(u), & (t\in\R\setminus\{0\}, \pm = \sgn(t)). \label{eq:def-Vintegral}
\end{align}

The main point of the definition~\eqref{eq:def-Gxs}, as we shall see, is that it is both easy to relate to~$h^\E_\gamma(x)$, and to estimate asymptotically.

\begin{prop}
  \begin{enumerate}
  \item The map~$\cG(x, \cdot)$ extends to a meromorphic function in~$\{\Re(s)>0\}$, which is analytic at~$1/2$.
  \item We have~$\cG(x, 1/2) = \E_\phi(\gamma x)$.
  \item As~$\abs{x}\to \infty$, we have
    \begin{equation}
      \cG(x, 1/2) = \chi(\gamma) (V_\phi(x-x_0) - i^k \E(x_0)) + u_\gamma(x) \E(x) + O(\abs{x}^{-1+\eps}).\label{eq:G-asympt-partial}
    \end{equation}
  \end{enumerate}
\end{prop}
\begin{proof}
  To prove the first assertion, we change variables, and write
  \begin{align*}
    \cG(x, s) = {}& \int_x^{x_0} \phi(\gamma z) (\Im z)^{s-1/2} \df s(z) \\
    = {}& \int_{\gamma x}^{\infty} \phi(z) (\Im \gamma^{-1} z)^{s-1/2} \df s(z) \\
    = {}& \int_{\gamma x}^\infty \phi_\infty(\Im z) (\Im \gamma^{-1} z)^{s-1/2} \df s(z) \numberthis\label{eq:G-anacont-1} \\
    &{} + \int_{\gamma x}^\infty (\phi(z) - \phi_\infty(\Im z))\Big((\Im \gamma^{-1} z)^{s-1/2} - \frac{(\Im z)^{s-1/2}}{\abs{j(\gamma^{-1}, \gamma x)}^{2s-1}}\Big) \df s(z) \numberthis\label{eq:G-anacont-2} \\
    &{} + \frac{\E(\gamma x, s)}{\abs{j(\gamma^{-1}, \gamma x)}^{2s-1}}. \numberthis\label{eq:G-anacont-3}
  \end{align*}
  We note that~$j(\gamma^{-1}, \gamma x) = j(\gamma, x)^{-1}$, but we will not use this at this point.
  The meromorphic extension and regularity at~$1/2$ of the third summand~\eqref{eq:G-anacont-3} is a consequence of Proposition~\ref{prop:def-twisted-l}.
  Next, by the exponential decay of~$\phi(z) - \phi_\infty(\Im z)$ for large~$\Im z$, and since
  $$ (\Im \gamma^{-1} z)^{s-1/2} - \frac{(\Im z)^{s-1/2}}{\abs{j(\gamma^{-1}, \gamma x)}^{2s-1}} = (\abs{j(\gamma^{-1}, z)}^{1-2s} - \abs{j(\gamma^{-1}, \gamma x)}^{1-2s}) (\Im z)^{s-1/2} = O_x((\Im z)^{s+1/2}) $$
  as~$\Im z \to 0$, we deduce that the second integral~\eqref{eq:G-anacont-2} converges uniformly with respect to~$s$ on compacts of~$\{\Re(s)>0\}$, and is therefore an analytic function of~$s$ in the same region. It obviously vanishes at~$s=1/2$.

  Finally we focus on the first integral~\eqref{eq:G-anacont-1}. Writing~$z=\gamma x + iy$, we have
  $$ \Im \gamma^{-1} z = \frac{y}{\abs{j(\gamma^{-1}, \gamma x + iy)}^2} = \frac{y}{j(\gamma, x)^{-2} + (cy)^2}. $$
  Since~$j(\gamma, x) = c(x-x_0)$, we deduce
  \begin{align*}
    \int_{\gamma x}^\infty \phi_\infty(\Im z) (\Im \gamma^{-1} z)^{s-1/2} \df s(z)
    ={}& \int_0^\infty \phi_\infty(y) \Big(\frac{y}{j(\gamma, x)^{-2}+(cy)^2}\Big)^{s-1/2} \frac{\df y}y \\
    ={}& \abs{x-x_0}^{s-1/2} \int_0^\infty \phi_\infty\Big(\frac{y}{c^2\abs{x-x_0}}\Big) \Big(\frac{y}{1+y^2}\Big)^{s-1/2} \frac{\df y}y.
  \end{align*}
  Given the shape of~$\phi_\infty$ in~\eqref{eq:def-phi-cuspidalpart}, we deduce that for~$\Re(s)>1$, the integral~\eqref{eq:G-anacont-1} is a linear combination of
  $$ g_\tau(s) := \int_0^\infty \frac{y^{\tau+s-1/2}}{(1+y^2)^{s-1/2}}\frac{\df y}y $$
  for~$\tau \in \{s_\phi, 1-s_\phi\}$, and its derivative with respect to~$\tau$.
  But for~$\Re(s)>\frac12+\abs{\Re(\tau)}$, we have by~\cite[8.380.3]{GradshteynRyzhik2007}
  $$ g_\tau(s) = \frac{\Gamma(\frac{s+\tau-1/2}2) \Gamma(\frac{s-\tau-1/2}2)}{2\Gamma(s-1/2)}. $$
  Thus, both~$g_\tau$ and~$\frac{\df}{\df\tau} g_\tau$ extend meromorphically to~$\C$, are regular at~$s=1/2$ and vanish there.
  This finishes the proof of the first item.
  Evaluating at~$s=1/2$, the terms~\eqref{eq:G-anacont-1} and \eqref{eq:G-anacont-2} vanish and we are left with~$\cG(x, 1/2) = \E(\gamma x)$ as announced in the second item.

  To prove the asymptotic estimate~\eqref{eq:G-asympt-partial}, we split~$\int_Z = \int_{Z^+} + \int_{Z^-}$, and write accordingly
  $$ \cG(x, s) = \cG^+(x, s) + \cG^-(x, s). $$
  It is clear that~$\cG^+(x, \cdot)$ is defined and holomorphic on~$\C$, since the integration path~$Z^+$ is compact in~$\H$.
  We write
  $$ \cG^+(x, s) = \int_{Z^+} u_\gamma(z)(\phi(z) - \phi_\infty(\Im z)) (\Im z)^{s-1/2} \df s(z) + \int_{Z^+} u_\gamma(z) \phi_\infty(\Im z) (\Im z)^{s-1/2} \df s(z). $$
  We recall that~$Z^+$ is contained in~$\{\Im z\geq Y\}$ and~$Y\to \infty$.
  By the Fourier expansion~\eqref{fexp} at~$\infty$ and the exponential decay of the Whittaker function, we deduce that~$\abs{\phi(z) - \phi_\infty(\Im z)} \ll Y^{-A}$ for any fixed~$A>0$.
  Using the triangle inequality and a rough estimate of the hyperbolic length of~$Z^+$, it follows that for~$\abs{s}\ll 1$, we have
  $$ \int_{Z^+} u_\gamma(z)(\phi(z) - \phi_\infty(\Im z)) (\Im z)^{s-1/2} \df s(z) = O_A(Y^{-A} \abs{x}) = O(\abs{x}^{-2}). $$
  by picking~$A$ large enough in terms of~$\eps$.
  Evaluating at~$s=1/2$, we deduce
  \begin{equation}
    \label{eq:Gplus-estim}
    \cG^+(x, 1/2) = O(\abs{x}^{-2}) + \int_{Z^+} u_\gamma(z) \phi_\infty(\Im z) (\Im z)^{s-1/2} \df s(z).
  \end{equation}
  
  On the other hand, we let~$\mu, \nu$ be shorthands for~$\mu_x(y)$ and~$\nu_x(y)$. We insist that these depend on~$y$.
  By parametrizing~$Z^-$, we have
  \begin{align*}
    \cG^-(x, s) = {}& \int_0^Y u_\gamma(x - \nu + iy) \phi(x - \nu + iy) y^{s-1/2} \mu \frac{\df y}y \\
    &{} - \int_0^Y u_\gamma(x_0 + \nu + iy) \phi(x + \nu + iy) y^{s-1/2} \mu \frac{\df y}y \\
    = {}& \cG_1^-(x, s) - \cG_2^-(x, s),
  \end{align*}
  say. We split
  $$ \cG_1^-(x, s) = \sum_{i=1}^5 G_i(x, s), $$
  where
  \begin{align*}
    G_1(x, s) := {}& \int_0^Y (u_\gamma(x-\nu+iy) \phi(x-\nu+iy) \mu - u_\gamma(x) \phi(x+iy)) y^{s-1/2}\frac{\df y}y, \\
    G_2(x, s) := {}& u_\gamma(x) \E(x, s) = u_\gamma(x) \int_0^\infty (\phi(x+iy) - \phi_\infty(y))y^{s-1/2}\frac{\df y}y, \\
    G_3(x, s) := {}& -u_\gamma(x) \int_Y^\infty (\phi(x+iy) - \phi_\infty(y)) y^{s-1/2} \frac{\df y}y, \\
    G_4(x, s) := {}& \int_0^Y (u_\gamma(x) - u_\gamma(x-\nu+iy) \mu)\phi_\infty(y) y^{s-1/2} \frac{\df y}y, \\
    G_5(x, s) := {}& \int_0^Y u_\gamma(x-\nu+iy) \phi_\infty(y) y^{s-1/2} \frac{\df y}y.
  \end{align*}
  All five terms here are meromorphic on~$\{\Re(s)>0\}$ and analytic at~$s=1/2$:
  For~$G_1$, this follows from the bounds~\eqref{eq:bound-phi}, \eqref{eq:bounds-mu} and \eqref{eq:bounds-nu}.
  For~$G_2$, this follows from Proposition~\ref{prop:def-twisted-l}.
  For~$G_3$, this follows from uniform convergence.
  For~$G_4$ and~$G_5$, this follows from the bound~$\phi_\infty(y) \ll y^{\sigma+\eps}$ as~$y\to 0$, where~$\sigma = 1-\Re(s_\phi)>0$ if~$\Re(s_\phi)>1$ or~$\sigma=\Re(s_\phi)>0$ if~$s_\phi\equiv \frac k2\pmod 1$. Here we used the fact that in the latter case we have~$B_\infty=0$ in~\eqref{eq:def-phi-cuspidalpart}.
  
  We now use the more precise bounds, valid for~$y\leq Y$,
  \begin{equation}
    \abs{\mu - 1} \ll x^{-2} y^2, \qquad \abs{\nu} \ll \abs{x}^{-1} y^2,\label{eq:bounds-mu-nu-finer}
  \end{equation}
  and
  \begin{equation}
    \abs{u_\gamma(x-\nu+iy) - u_\gamma(x)} \ll \abs{\e^{ik \arctan(y/(x-x_0-\nu))} - 1} \ll \abs{x}^{-1} y.\label{eq:bound-ugamma-finer}
  \end{equation}
  Using these bounds along with~\eqref{eq:bound-phi} and \eqref{eq:bound-phi-diff}, we obtain
  $$ (u_\gamma(x-\nu+iy) \phi(x-\nu+iy) \mu - u_\gamma(x) \phi(x+iy)) \ll y^{1/2+\eps} \abs{x}^{-1}, \qquad (0\leq y\leq Y), $$
  and therefore
  $$ G_1(x, \tfrac12) = O(\abs{x}^{-1+\eps}). $$
  Concerning~$G_3$, using the bound~$\phi(x+iy) - \phi_\infty(y) = O_A(y^{-A})$ for any fixed~$A>0$, we deduce
  $$ G_3(x, \tfrac12) = O_A(Y^{-A}) = O(\abs{x}^{-1}) $$
  by picking~$A$ large enough in terms of~$\eps$.
  We bound~$G_4$ by using again~\eqref{eq:bounds-mu-nu-finer} and \eqref{eq:bound-ugamma-finer}, obtaining
  $$ G_4(x, \tfrac12) \ll \int_0^Y y^{1-\Re(s_\phi)} \abs{x}^{-1} \df y \ll \abs{x}^{-1+\eps}. $$
  Finally, we have~$G_2(x, \tfrac12) = u_\gamma(x) \E(x)$, and leaving~$G_5(x, \tfrac12)$ unevaluated for now, we conclude the analytic continuation of~$\cG_1^-(x, \cdot)$ and the expression
  \begin{equation}
    \cG^-_1(x, \tfrac12) = u_\gamma(x) \E(x) + \int_0^Y u_\gamma(x-\nu+iy) \phi_\infty(y) \mu \frac{\df y}y + O(\abs{x}^{-1+\eps}).\label{eq:estim-Gminus-1}
  \end{equation}

  We evaluate~$\cG^-_2(x, s)$ in a similar way, except for the terms involving~$u_\gamma$. First we notice that for~$y>0$, the quantity
  $$ u_\gamma(x_0+iy) = \chi(\gamma) (i \sgn(c))^k $$
  is independent of~$y$. Then we bound
  $$ \abs{u_\gamma(x_0+\nu+iy) - u_\gamma(x_0+iy)} \ll \abs{\e^{ik \arctan(\nu/y)} - 1} \ll \abs{\nu}/y \ll y \abs{x}^{-1}. $$
  Using this bound in place of~\eqref{eq:bound-ugamma-finer}, we may reproduce the above computations, to the effect that~$\cG^-_2(x, \cdot)$ extends to a meromorphic function on~$\{\Re(s)>0\}$ which is analytic at~$1/2$, and
  \begin{equation}
    \cG^-_2(x, \tfrac12) = \chi(\gamma) i^k \E(x_0) + \int_0^Y u_\gamma(x_0+\nu + iy) \phi_\infty(y) \mu \frac{\df y}y + O(\abs{x}^{-1+\eps}).\label{eq:estim-Gminus-2}
  \end{equation}
  
  The map~$\cG^-(x, \cdot)$ therefore extends meromorphically to~$\{\Re(s)>0\}$ since~$\cG_1^-$ and~$\cG_2^-$ do.
  Summing our two estimates~\eqref{eq:estim-Gminus-1} and \eqref{eq:estim-Gminus-2}, and by parametrizing again~$Z^-$, we get
  \begin{equation}
    \cG^-(x, \tfrac12) = u_\gamma(x) \E(x) - \chi(\gamma) i^k \E(x_0) + \int_{Z^-} u_\gamma(z) \phi_\infty(\Im z) \df s(z) + O(\abs{x}^{-1+\eps}). \label{eq:Gminus-estim}
  \end{equation}
  
  By summing the estimates~\eqref{eq:Gplus-estim} and \eqref{eq:Gminus-estim}, we deduce
  $$ \cG(x, \tfrac12) = u_\gamma(x) \E(x) - \chi(\gamma) i^k \E(x_0) + \int_{Z} u_\gamma(z) \phi_\infty(\Im z) \df s(z) + O(\abs{x}^{-1+\eps}). $$
  There remains to evaluate the last integral. To do this we change variables $z = \abs{x-x_0} u + x_0$.
  This sends the geodesic circle~$Z = G(x, x_0)$ to~$G(\eta, 0)$, where~$\eta = \sgn(x-x_0)$, and matches the corresponding geodesic lengths.
  We also have~$\Im z = \abs{x-x_0} \Im u$, and since we assumed~$c>0$, we get~$u_\gamma(z) = \chi(\gamma) (u / \abs{u})^k$. We deduce as claimed
  $$ \int_{Z} u_\gamma(z) \phi_\infty(\Im z) \df s(z) = \chi(\gamma) \int_{G(\eta, 0)} \Big(\frac{u}{\abs{u}}\Big)^k \phi_\infty(\abs{x-x_0} \Im u) \df s(u). $$
\end{proof}

Comparing the second and third items of the previous proposition, we have
\begin{align*}
  h_\gamma^\E(x) = {}& \E(\gamma x) - u_\gamma(x) \E(x) \\
  ={}& \chi(\gamma) (V_\phi(x-x_0) - i^k \E(x_0)) + O(\abs{x}^{-1+\eps}).
\end{align*}
The asymptotic formula~\eqref{eq:hgamma-E-asymptotic} follows using the expression~\eqref{eq:def-Vintegral} for~$V_\phi(x-x_0)$ and~\eqref{eq:def-phi-cuspidalpart} for~$\phi_\infty$.

For~$k\in \Z$, $\eta \in \{\pm 1\}$ and~$\tau\in\C$, $\Re(\tau)>0$, define
\begin{equation}
  \Upsilon_{k, \eta}(\tau) := \int_{G(\eta, 0)} (u / \abs{u})^k (\Im u)^\tau \df s(u).\label{eq:def-Upsilon}
\end{equation}

The expressions~\eqref{eq:expr-A'-B'-not12} and \eqref{eq:expr-A'-B'-12} are clear from the definition of~$V_\infty$ and~$\Upsilon_{k, \eta}$.
To finish the proof of Theorem~\ref{th:hgamma-asymptotic}, we describe~$\Upsilon_{k, \eta}(\tau)$ in terms of~$\Gamma$-functions.

\begin{lemma}\label{lem:Upsilon}
  We have
  \begin{equation}\label{eq:expr-Upsilon}
    \Upsilon_{k, \eta}(\tau) = \e^{\pi i (2-\eta) k/4} \frac{2\pi \Gamma(\tau)}{4^\tau\Gamma(\frac{\tau+1+k/2}2)\Gamma(\frac{\tau+1-k/2}2)}.
  \end{equation}
\end{lemma}

\begin{proof}[Proof of Lemma~\ref{lem:Upsilon}]
  Assume first~$\eta=1$.
  In the integral~\eqref{eq:def-Upsilon}, we let~$u = \e^{i\theta/2} \cos(\theta/2)$ with~$\theta \in [0, \pi]$, so that~$\Im u = \frac12\sin(\theta)$, $u/\abs{u} = \e^{i\theta/2}$ and~$(\Im u)\df s(u) = \frac12\df \theta$.
  This gives
  $$ \Upsilon_{k, 1}(\tau) = 2^{-\tau} \int_0^\pi \e^{ik\theta/2} (\sin \theta)^{\tau-1}\df \theta. $$
  By equations~(3.631.1) and~(3.631.8) of\cite{GradshteynRyzhik2007}, the stated result follows for~$\eta=1$.

  When~$\eta=-1$, then we change variables~$u \gets -\overline{u}$ in~\eqref{eq:def-Upsilon}. This sends~$G(-1, 0)$ to~$G(1, 0)$, changes~$(u/\abs{u})^k$ to~$(-1)^k (u/\abs{u})^{-k}$ and leaves~$\Im u$, $\df s(u)$ unchanged.
  We deduce~$\Upsilon_{k, -1}(\tau) = (-1)^k \Upsilon_{-k, 1}(\tau)$, and by using the formula above for~$\eta=1$ we get the stated result.
\end{proof}

\begin{proof}[Proof of Corollary~\ref{cor:hgamma-asymptotic-finite}]
  Let~$x_1 := \gamma \infty$ and
  $$ x = \delta^{-1} + x_1. $$
  This is defined so that
  $$ \gamma^{-1}\infty - c^{-2}\delta = \gamma^{-1}x. $$
  We have also~$j(\gamma^{-1}, x) = -cx+a = -c \delta^{-1}$, and therefore
  $$ u_\gamma(\gamma^{-1}x) = \overline{u_{\gamma^{-1}}(x)} = \chi(\gamma) (-\sgn(\delta))^k. $$
  We deduce
  \begin{align*}
    h_\gamma^\E(\gamma^{-1}\infty - c^{-2}\delta) = {}& \E(x) - \overline{u_{\gamma^{-1}}(x)} \E(\gamma^{-1}x) \\
    = {}& - \chi(\gamma) (-\sgn(\delta))^k h_{\gamma^{-1}}^\E(x).
  \end{align*}
  As~$\delta \to 0$ with~$\sgn(\delta) = \pm$, we have~$x\to \pm \infty$, therefore Theorem~\ref{th:hgamma-asymptotic} applies and yields the announced estimate.
\end{proof}

For future reference, we list the special cases
\begin{equation}\label{eq:Upsilon-0-2}
  \begin{aligned}
    \Upsilon_{0, \pm}(\tfrac12) ={}& \frac{\pi^{3/2}}{\Gamma(3/4)^2}, &
    \Upsilon_{1, \pm}(\tfrac12) ={}& \tfrac{\eta+i}{\sqrt{2}} \pi, &
    \Upsilon_{2,\pm}(\tfrac12) ={}& i\eta \frac{4\pi^{3/2}}{\Gamma(1/4)^2}, \\
    \tfrac{\Upsilon'_{0,\pm}}{\Upsilon_{0,\pm}}(\tfrac12) ={}& -\tfrac\pi2-\log 2, & & &
    \tfrac{\Upsilon'_{2,\pm}}{\Upsilon_{2,\pm}}(\tfrac12) ={}& \tfrac\pi2-\log 2-2.
  \end{aligned}
\end{equation}
They are obtained using the values listed in~\cite[8.366]{GradshteynRyzhik2007}.

\subsection{Behaviour at infinity for generalized quantum modular forms}

Theorem~\ref{th:hgamma-asymptotic} admits an extension to the period function $h_{\overline{\gamma}}$ which was introduced in Section~\ref{sec:GQM}.
We recall the notation~\eqref{eq:extendedQM} and the definition~\eqref{eq:def-hgamma-E-ext}.

\begin{thm}\label{th:hgamma-asymptotic-ext}
  Assume that~$\overline{\gamma} \phi = \eta j_\gamma^k \phi$ for some~$\eta\in\C$. Then with the same notation as in Theorem~\ref{th:hgamma-asymptotic}, we have
  $$ h_{\overline{\gamma}}^\E(x) = \eta\big(A_\pm' \abs{x+x_0}^{s_\phi} + B_\pm' \abs{x+x_0}^{1-s_\phi} - i^k \E_\phi(-x_0)\big) + O_{\phi, \gamma, \eps}(\abs{x}^{-1+\eps}). $$
\end{thm}

We omit the proof, since it is identical to Theorem~\ref{th:hgamma-asymptotic}, starting instead from the integral
$$ \cG(x, s) := \int_{G(x, -x_0)} \eta j_\gamma^k(z) \phi(z) (\Im z)^{s - 1/2} \df s(z). $$
Theorem~\ref{th:hgamma-asymptotic-ext} will be useful when studying reciprocity laws in Section~\ref{sec:arith-applications}.
This will involve the action under the Fricke involution.

\section{Quantum modularity for the central value of the twisted \texorpdfstring{$L$}{L}-series}\label{sec:qmod-L}

In this section, we prove an analogue of Theorem~\ref{th:hgamma=qmf} for the following twisted Dirichlet series~:
\begin{itemize}
\item If~$s_\phi \not\equiv \frac k2 \pmod{1}$, then we let
  \begin{equation}\label{eq:def-Lxs-maass}
    L^{\pm}(\phi, x, s) := \begin{cases} \sum_{n>0} a(n) \cos(2\pi n x) \abs{n}^{1/2-s}, & (\pm = +), \\ i \sum_{n>0} a(n) \sin(2\pi n x) \abs{n}^{1/2-s}, & (\pm = -). \end{cases}
  \end{equation}
\item If~$s_\phi \equiv \frac k2 \pmod{1}$, then we let
  \begin{equation}\label{eq:def-Lxs-holo}
    L(\phi, x, s) := \sum_{n>0} a(n) \e(nx).
  \end{equation}
\end{itemize}
This corresponds to the Maa\ss{}, respectively holomorphic cases studied in~\cite[Appendices~A.3 and~A.4]{KoMiVa02}.
Using the assumption~\eqref{eq:bound-an-rough}, the right-hand sides in both definitions are defined and analytic with respect to~$s$ in~$\{\Re(s)>1\}$.

\subsection{Properties of a Mellin integral of the Whittaker function}

For any~$\alpha \in \C$ and~$\beta\in\C$, $0\leq \Re(\beta)<1/2$, we define
\begin{equation}\label{eq:def-Wmellin}
  \Omega_\beta(\alpha, s) := \int_0^\infty W_{\alpha, \beta}(y) y^{s-1/2} \frac{\df y}y.
\end{equation}
Since~$W_{\alpha, \beta}(y) \ll y^{1/2-\Re(\beta)+\eps}$ (see \cite[eq.~(4.19)]{DukeEtAl2002}), the integral certainly converges absolutely for~$\Re(s)>\abs{\Re(\beta)}$.
In particular, it is always regular at~$s=1/2$.
At this particular point, we have by~\cite[(7.611.1)]{GradshteynRyzhik2007}
\begin{equation}\label{eq:Wmellin-gamma-half}
  \Omega_\beta(\alpha, \tfrac12) = \frac{\pi^{3/2}2^\alpha}{\cos(\pi \beta)\Gamma(\frac34-\frac12\alpha+\frac12\beta)\Gamma(\frac34-\frac12 \alpha - \frac12\beta)}.
\end{equation}

In the rest of this section,~$\beta$ is fixed.
Mellin integrals of~$W_{\alpha, \beta}$ have been studied in~\cite[section~8]{DukeEtAl2002} (see also~\cite[section~12]{Young2019}), and we will return to these works below in Remark~\ref{rmk:Delta}.
We are interested in the quantity
\begin{equation}
  \Delta_\beta(\alpha, s) := (\tfrac12+\alpha+\beta)(\tfrac12+\alpha-\beta)\Omega_\beta(-\alpha-1, s)\Omega_\beta(\alpha, s) + \Omega_\beta(-\alpha, s)\Omega_\beta(\alpha+1, s), \label{eq:def-Delta}
\end{equation}
which will appear later as a certain determinant.

\begin{prop}\label{prop:props-Delta}
  The map~$\Delta_\beta(\alpha, \cdot)$ extends to a meromorphic function on~$\C$. It is independent of~$\alpha$, and
  $$ \Delta_\beta(\alpha, s) = 4^s \Gamma(s+\beta) \Gamma(s-\beta). $$
  In particular,~$\Delta_\beta(\alpha, s) \neq 0$.
\end{prop}

The special case~$s=1/2$ of Proposition~\ref{prop:props-Delta} can be checked directly from the explicit expression~\eqref{eq:Wmellin-gamma-half}.

\begin{proof}
  By analyticity, we may assume that~$\Re(s)>\abs{\Re(\beta)}$, and also that~$s \not\in \alpha+\frac12+\Z$.
  By equations~(7.621.3) and (9.131.1) of\cite{GradshteynRyzhik2007}, we have for~$\Re(s)>\abs{\Re(\beta)}$ the equality
  \begin{equation}\label{eq:rel-Wmellin-hypergeom}
    \Omega_\beta(\alpha, s) = \frac{\Gamma(s+\beta)\Gamma(s-\beta)}{\Gamma(s-\alpha+\tfrac12)} F(s-\beta, s+\beta; s-\alpha+\tfrac12; \tfrac12).
  \end{equation}
  Letting
  $$ (a, b, c) := (s-\beta, s+\beta, s-\alpha+\tfrac12), $$
  we then have~$s+\alpha+\frac12 = a+b-c$, $\frac12+\alpha+\beta = b-c$, $\frac12+\alpha-\beta = a-c$, and we eventually arrive at
  \begin{align*}
    \Delta_\beta(\alpha, s) = {}& \frac{\Gamma(a)^2\Gamma(b)^2}{\Gamma(a+b-c+1)\Gamma(c+1)}\big((c-a)(c-b)F(a,b;c+1;\tfrac12) F(a,b;a+b-c+1;\tfrac12) \\
    & \hspace{12em} + c(a+b-c)F(a,b;c;\tfrac12) F(a,b;a+b-c;\tfrac12)\big).
  \end{align*}
  An identity between hypergeometric functions, which we have stated and proved in Lemma~\ref{lem:hypergeom-2} in the appendix, allows us to deduce~$\Delta_\beta(\alpha, s) = 2^{a+b}\Gamma(a)\Gamma(b)$, which is the claimed equality.
\end{proof}

\begin{rem}\label{rmk:Delta}
  In~\cite{DukeEtAl2002}, section~8, integrals closely related to~$\Omega_\beta(\alpha, s)$ are studied.
  A minor mistake in the computations there was recently corrected in~\cite{Young2019}, section 12.
  With the notation from~\cite[section~12]{Young2019}, we find for~$k\in \Z_{\geq 0}$ that
  \begin{align*}
    \Omega_\beta(k/2, s) ={}& \pi^{-1/2} 4^{s-1} (\Phi_k^+(s, \beta) + \Phi_k^-(s, \beta)), \\
    \Omega_\beta(-k/2, s) = {}& \pi^{-1/2} 4^{s-1} \frac{\Gamma(\beta+\frac{1-k}2)}{\Gamma(\beta+\frac{1+k}2)} (\Phi_k^+(s, \beta) - \Phi_k^-(s, \beta)).
  \end{align*}
  Then a quick computation shows that Proposition~\ref{prop:props-Delta} at~$\alpha=k/2$ is essentially equivalent to the equality
  $$ p_k^+(s, \beta) p_{k+2}^-(s, \beta) - p_{k+2}^+(s, \beta) p_k^-(s, \beta) = 2\frac{\Gamma(\beta+\frac{1+k}2)}{\Gamma(\beta+\frac{1-k}2)}. $$
  Notice the right-hand side is independent of~$s$.
  The sequences of polynomials~$(p_k^\pm)$ satisfy a recurrence relation, see~\cite[eq.~(12.2)]{Young2019}. It would be interesting to know if one can show the above relation, and therefore Proposition~\ref{prop:props-Delta} for~$\alpha=k/2$, by induction on~$k$ instead of the arguments presented here.
  
  As we will see, the quantity~$\Delta_\beta(\alpha, s)$ arises from the computation of the discriminant in a~$2\times 2$ linear system~\eqref{eq:system-Ephipsi-L} below, which is somewhat analogous to a Wronskian: the second row arises from Mellin integrals of Whittaker functions similar to the first row, but to which a level-raising operator was applied.
  It would be interesting to know if a more direct argument could be used.
  Our early attempts were unsuccessful.
\end{rem}

\subsection{Relation with the period integral}\label{sec:relperiod}

In the setting described in Section~\ref{sec:setting}, we let $\phi$ be a Maa{\ss} form for $\Gamma$ of weight $k$ and eigenvalue~$s_\phi$. In this section, we relate the twisted $L$-series~\eqref{eq:def-Lxs-maass}-\eqref{eq:def-Lxs-holo} to the Eichler integral~\eqref{eq:def-eichlerint}.
This can be seen as a generalization of computations done in Sections A.3 and A.4 of~\cite{KoMiVa02}.

Recall the definition \eqref{eq:levelraising} of the raising and lowering operators.
We will require two different forms related to~$\phi$. When~$\frac12 \leq \Re(s_\phi) < 1$ and~$s_\phi\neq 1/2$, let
\begin{equation}
  \psi := R_k \phi. \label{eq:def-psi}
\end{equation}
This is a weight~$k+2$ Maa\ss{} form for~$\Gamma$ with the same eigenvalue~$s_\phi$ and nebentypus~$\chi$ as~$\phi$.
In particular, it follows from Proposition~\ref{prop:def-twisted-l} and Theorem~\ref{th:hgamma=qmf} that for all~$x\in C(\Gamma)$, the map~$\E(\psi, x, \cdot)$ extends to a meromorphic function on~$\C$ which is regular at~$1/2$, and that the map~$\E_\psi(x) = \E(\psi, x, 1/2)$ is such that for all~$\gamma\in \Gamma$, the difference
$$ h_\gamma^\E(\psi; x) := \E_\psi(\gamma x) - j_\gamma(x)^{k+2} \chi(\gamma) \E_\psi(x) $$
extends to a Hölder continuous function of~$x$ on~$\R\setminus\{\gamma^{-1}\infty\}$, with~Hölder exponent~$1/2-\eps$ (and~$1-\eps$ if~$\phi$ is cuspidal). Note that since~$j_\gamma(x) \in \{\pm 1\}$, we obviously have
\begin{equation}
  h_\gamma^\E(\psi; x) = \E_\psi(\gamma x) - j_\gamma(x)^k \chi(\gamma) \E_\psi(x). \label{eq:def2-hpsi}
\end{equation}

When~$\phi$ is associated to a holomorphic form, say~$\Re(s_\phi) = \ell/2$ with~$1\leq \ell \leq k$, $\ell\equiv k\pmod{2}$, then we let
\begin{equation}\label{eq:def-fholom}
  f := \Lambda_{\ell+2} \dotsb \Lambda_{k-2} \Lambda_k \phi
\end{equation}
be the underlying Maa{\ss} form of weight~$\ell$, same eigenvalue~$s_\phi = \ell/2$ and same nebentypus as~$\phi$. The map~$z\mapsto y^{-\ell/2} f(z)$ is holomorphic of weight~$\ell$, see~\cite[p.~507]{DukeEtAl2002}.
Similarly as above, it follows from Proposition~\ref{prop:def-twisted-l} and Theorem~\ref{th:hgamma=qmf} that the map~$\E_f(x) = \E(f, x, 1/2)$ is such that for all~$\gamma\in\Gamma$, the difference
$$ h_\gamma^\E(f; x) = \E_f(\gamma x) - j_\gamma(x)^k \chi(\gamma) \E_f(x) $$
extends to a~$(1/2-\eps)$-Hölder continuous function on~$\R\setminus\{\gamma^{-1}\infty\}$ (and~$(1-\eps)$-Hölder if~$\phi$ is cuspidal). Here we have used that~$k$ and~$\ell$ have the same parity.

We will pass through the intermediate object
$$ U^\pm(\phi, x, s) := \sum_{\substack{n\neq 0 \\ \sgn(n) = \pm}} a(n) \e(nx) \abs{n}^{1/2-s}, $$
defined and analytic in~$\{\Re(s)>1\}$.
This is for convenience only, and we will switch soon thereafter to~$L^\pm$ itself.

\begin{lemma}\label{lem:system-Ephipsi-L}
  For~$\Re(s)>1$, the following holds:
  \begin{itemize}
  \item If~$s_\phi \not\equiv \frac k2\pmod{1}$, then
    \begin{equation}\label{eq:system-Ephipsi-L}
      \left\{
        \begin{aligned}
          (4\pi)^{s-1/2} \E(\phi, x, s) ={}& \Omega_{it_\phi}(\tfrac k2, s) U^+(\phi, x, s) + \Omega_{it_\phi}(-\tfrac k2, s) U^-(\phi, x, s) , \\
          (4\pi)^{s-1/2} \E(\psi, x, s) ={}& -\Omega_{it_\phi}(\tfrac k2+1, s) U^+(\phi, x, s) + (s_\phi+\tfrac k2)(1-s_\phi+\tfrac k2) \Omega_{it_\phi}(-\tfrac k2-1, s) U^-(\phi, x, s).
        \end{aligned}
      \right.
    \end{equation}
  \item If~$s_\phi = k/2$, then
    \begin{equation}\label{eq:rel-E-U-holo}
      (2\pi)^{s-1/2} \E(\phi, x, s) = 2^{k/2}\Gamma(\tfrac k2+s-\tfrac12) U^+(\phi, x, s).
    \end{equation}
  \end{itemize}
\end{lemma}

\begin{proof}
  Consider first the case~$s_\phi\not\equiv \frac k2\pmod{1}$.
  For~$\Re(s)>1$ we have
  \begin{align*}
    \E(\phi, x, s) = {}& \sum_{n\neq 0} a(n) \e(nx) \int_0^\infty W_{k/2\sgn(n), it_\phi}(4\pi\abs{n} y) y^{s-1/2} \frac{\df y}y \\
    = {}& \sum_{n\neq 0} a(n) \e(nx) (4\pi \abs{n})^{1/2-s} \int_0^\infty W_{k/2\sgn(n), it_\phi}(y) y^{s-1/2} \frac{\df y}y \\
    = {}& (4\pi)^{1/2-s} \sum_{n\neq 0} a(n) \e(nx) \abs{n}^{1/2-s} \Omega_{it_\phi}(\tfrac k2\sgn(n), s).
  \end{align*}
  This gives the first claimed equation.

  To prove the second claimed equation, we recall the action of~$R_k$ on the Fourier expansion described by~\cite[eqs.~(4.25)-(4.26)]{DukeEtAl2002}, frow which we obtain for all~$z\in\C$
  $$ \psi(z) = \psi_\infty(y) + \sum_{n\neq 0} a_\psi(n) \e(nx) W_{(\frac k2+1)\sgn(n), it_\phi}(4\pi\abs{n}y), $$
  where
  \[ a_\psi(n) = \begin{cases}
      -a(n) & (n>0), \\
      (s_\phi+\tfrac k2)(1-s_\phi+\tfrac k2) a(n) & (n<0).  
    \end{cases} \]
  We deduce that
  \begin{equation}
    \begin{aligned}
      U^+(\psi, x, s) ={}& - U^+(\phi, x, s), \\
      U^-(\psi, x, s) ={}& (s_\phi+\tfrac k2)(1-s_\phi+\tfrac k2) U^-(\phi, x, s).
    \end{aligned}\label{eq:Lpsi-rel-Lphi}
  \end{equation}
  On the other hand, by the above computations applied to~$\psi$ instead of~$\phi$, we deduce
  \begin{equation}
    (4\pi)^{s-1/2} \E(\psi, x, s) = \Omega_{it_\phi}(\tfrac k2+1, s) U^+(\psi, x, s) + \Omega_{it_\phi}(-\tfrac k2-1, s) U^-(\psi, x, s).\label{eq:Ppsi-rel-Lpsi}
  \end{equation}
  Grouping~\eqref{eq:Lpsi-rel-Lphi} and \eqref{eq:Ppsi-rel-Lpsi} proves the second claimed equation and completes the proof when~$\phi$ is not associated with a holomorphic form.
  
  Assume that~$s_\phi = k/2$, with~$k\geq 1$. Then we have the explicit expression~\cite[eq.~(4.21)]{DukeEtAl2002}
  $$ W_{\frac k2, \frac{k-1}2}(y) = y^{k/2} \e^{-y/2}, $$
  from which we deduce in the larger region of absolute convergence~$\Re(s)>0$ that
  $$ \E(\phi, x, s) = (4\pi)^{s-1/2} 2^{k/2+s-1/2} \Gamma(\tfrac k2 + s - \tfrac12) \sum_{n>0} a(n) \e(nx) n^{1/2-s}. $$
  This gives our second claim.
\end{proof}

\begin{lemma}\label{lem:Lxs-extension}
  The maps~$L^{\pm}(\phi, x, \cdot)$ extend to meromorphic functions on the domain~$\{\Re(s)>0\}$, which are regular at~$1/2$.
  
  More precisely, the following holds:
  \begin{itemize}
  \item If~$s_\phi \not\equiv\frac k2\pmod 1$, then for some constants~$c_\phi^\pm$ and~$c_\psi^\pm$, we have
    \begin{equation}
      L^{\pm}(\phi, x, 1/2) = c_\phi^\pm \E(\phi, x, 1/2) + c_\psi^\pm \E(\psi, x, 1/2),\label{eq:rel-L-E-1/2-maass}
    \end{equation}
    where~$\psi$ was defined in~\eqref{eq:def-psi}.
  \item If~$s_\phi = \ell/2$ for~$1\leq \ell \leq k$, $\ell \equiv k\pmod{2}$, then for some constant~$c_\phi^\pm$, we have
    \begin{equation}\label{eq:rel-L-E-holo}
      L(\phi, x, 1/2) = c_f \E(f, x, 1/2),
    \end{equation}
    where~$f$ was defined in~\eqref{eq:def-fholom}.
  \end{itemize}
\end{lemma}

The constants~$c_\phi^\pm, c_\psi^\pm, c_f^\pm$ depend only on~$s_\phi$ and~$k$. They are given in~\eqref{eq:expr-cpm} and~\eqref{eq:expr-cpm-holo} below.

\begin{proof}
  First we assume~$s_\phi\not\equiv\frac k2\pmod1$.
  Let
  $$ \beta = it_\phi, $$
  which satisfies~$\Re(\beta)\in[0, 1/2)$. Note that
  $$ (s_\phi+\tfrac k2)(1-s_\phi+\tfrac k2) = (\tfrac{1+k}2+\beta)(\tfrac{1+k}2-\beta). $$
  With this notation, the system of equations from Lemma~\ref{lem:system-Ephipsi-L} has discriminant~$\Delta_\beta(\frac12, s)$, where~$\Delta_\beta$ was defined in~\eqref{eq:def-Delta}. Since it does not vanish by Proposition~\ref{prop:props-Delta}, we have the existence of functions~$g_\phi^{\pm}$ holomorphic on~$\{\Re(s)>0\}$ such that for~$\Re(s)>1$,
  \begin{equation}
    U^\pm(\phi, x, s) = \Delta(\tfrac k2, s)^{-1}( g_\phi^{\pm}(s) \E(\phi, x, s) + g_\psi^\pm(s) \E(\psi, x, s)),\label{eq:rel-U-P}
  \end{equation}
  and this identity between meromorphic functions holds for~$\Re(s)>0$ by analytic continuation. More precisely, we have by Cramer's rule
  \begin{align*}
    g_\phi^+(s) = {}& (4\pi)^{s-1/2} \Omega_{it_\phi}(-\tfrac k2-1, s) (s_\phi+\tfrac k2)(1-s_\phi+\tfrac k2), \\
    g_\psi^+(s) = {}& -(4\pi)^{s-1/2} \Omega_{it_\phi}(-\tfrac k2, s), \\
    g_\phi^-(s) = {}& (4\pi)^{s-1/2} \Omega_{it_\phi}(\tfrac k2+1, s), \\
    g_\psi^-(s) = {}& (4\pi)^{s-1/2} \Omega_{it_\phi}(\tfrac k2, s).
  \end{align*}

  Now, by the hypothesis~$(H_4)$ and computing the action of~$Q_{sk}$ on the Fourier expansion~\cite[equation~(4.70)]{DukeEtAl2002}, we have for all~$n>0$ that
  $$ a(n) = \epsilon_\phi \frac{\Gamma(s_\phi-\frac k2)}{\Gamma(s_\phi+\frac k2)} a(-n). $$
  We deduce, for all~$x\in C(\Gamma)$,
  \begin{equation}\label{eq:rel-Uplusminus} U^+(\phi, -x, s) = \epsilon_{\phi} \frac{\Gamma(s_\phi-\frac k2)}{\Gamma(s_\phi+\frac k2)} U^-(\phi, x, s). \end{equation}
  Therefore, we have
  $$ L^\pm(\phi, x, s) = \frac12\Big(U^+(\phi, x, s) \pm \epsilon_\phi \frac{\Gamma(s_\phi-\frac k2)}{\Gamma(s_\phi+\frac k2)} U^-(\phi, x, s)\Big), $$
  we deduce that for any~$x\in C(\Gamma)$, the map~$L^\pm(\phi, x, \cdot)$ also extends to~$\Re(s)>0$.
  By Proposition~\ref{prop:props-Delta}, we have~$\Delta(\tfrac k2, 1/2)\neq 0$, therefore both sides of~\eqref{eq:rel-U-P} are regular at~$1/2$. The relation~\eqref{eq:rel-L-E-1/2-maass} follows with
  $$ c_\ast^\pm = (2\Delta(\tfrac k2, 1/2))^{-1} \Big(g_\ast^+(1/2) \pm \epsilon_\phi \frac{\Gamma(s_\phi-\frac k2)}{\Gamma(s_\phi+\frac k2)}  g_\ast^-(1/2)\Big), \qquad (\ast \in \{\phi, \psi\}). $$
  The coefficients here are expressed using~\eqref{eq:Wmellin-gamma-half} and the functional equation of the~$\Gamma$ function as
  \begin{align*}
    g_\phi^\pm(\tfrac12) = {}& \frac{\pi^{3/2} 2^{1\mp k/2}}{\cosh(\pi t_\phi) \Gamma(\frac{1\pm k}4+i\frac{t_\phi}2) \Gamma(\frac{1\pm k}4-i\frac{t_\phi}2)}, \\
    g_\psi^\pm(\tfrac12) = {}& \frac{\mp \pi^{3/2} 2^{1\mp k/2}}{\cosh(\pi t_\phi) \Gamma(\frac{3\pm k}4+i\frac{t_\phi}2) \Gamma(\frac{3\pm k}4-i\frac{t_\phi}2)}.
  \end{align*}
  Along with the explicit expression for~$\Delta(\frac k2, 1/2)$ from Proposition~\ref{prop:props-Delta}, and upon using the complement and duplication formula for the Euler~$\Gamma$ function, this yields
  \begin{equation}
    \begin{aligned}
      c_\phi^{\pm} = {}& \frac{\sqrt{\pi}2^{-1-k/2}}{\Gamma(\frac{1+k}4+i\frac{t_\phi}2) \Gamma(\frac{1+k}4-i\frac{t_\phi}2)} \Big(1 \pm \epsilon_\phi \frac{\cos(\frac\pi2(\frac{1+k}2+it_\phi))}{\sin(\frac\pi2(\frac{1+k}2-it_\phi))} \Big), \\
      c_\psi^{\pm} = {}& \frac{\sqrt{\pi}2^{-1-k/2}}{\Gamma(\frac{3+k}4+i\frac{t_\phi}2) \Gamma(\frac{3+k}4-i\frac{t_\phi}2)} \Big(-1 \pm \epsilon_\phi \frac{\sin(\frac\pi2(\frac{1+k}2+it_\phi))}{\cos(\frac\pi2(\frac{1+k}2-it_\phi))} \Big).
    \end{aligned}\label{eq:expr-cpm}
  \end{equation}
  Note that the quotients of~$\cos, \sin$ on the right-hand sides simplify somehow when~$k\in2\Z_{\geq 0}$, for then we have
  \begin{equation}\label{eq:cphipsi-cossin-keven}
    \frac{\cos(\frac\pi2(\frac{1+k}2+it_\phi))}{\sin(\frac\pi2(\frac{1+k}2-it_\phi))} = \frac{\sin(\frac\pi2(\frac{1+k}2+it_\phi))}{\cos(\frac\pi2(\frac{1+k}2-it_\phi))} = (-1)^{k/2}, \qquad (k\text{ even}).
  \end{equation}

  Next, we assume~$s_\phi = \ell / 2$ for some integer~$1\leq \ell \leq k$ of the same parity as~$k$.
  Then, by definition, we have~$L(\phi, x, s) = U^+(\phi, x, s)$.
  By the action of the level-lowering operator~$\Lambda_k$ on the Fourier expansion~\cite[eqs.~(4.27), (4.28) and~(4.30)]{DukeEtAl2002}, we deduce that for~$n>0$,
  $$ a_f(n) = (-1)^{(k-\ell)/2} \frac{\Gamma(\frac{k-\ell}2+1)\Gamma(\frac{k+\ell}2)}{\Gamma(\ell)} a_\phi(n). $$
  Hence, for~$\Re(s)>1$
  $$ L(\phi, x, s) = U^+(\phi, x, s) = (-1)^{(k-\ell)/2} \frac{\Gamma(\ell)}{\Gamma(\frac{k-\ell}2+1)\Gamma(\frac{k+\ell}2)} U^+(f, x, s). $$
  Since~$f$ has eigenvalue~$s_f = s_\phi = \ell/2$ and weight~$\ell$, we can apply~\eqref{eq:rel-E-U-holo} to obtain
  $$ L(\phi, x, s) = c_f(s) \E(f, x, s), $$
  where
  \begin{equation}
    c_f(s) := \frac{(-1)^{(k-\ell)/2} (2\pi)^{s-1/2} \Gamma(\ell)}{2^{\ell/2} \Gamma(\frac\ell2+s-\frac12) \Gamma(\frac{k-\ell}2+1)\Gamma(\frac{k+\ell}2)}.
    \label{eq:def-cfs}
  \end{equation}
  This shows the claimed analytic extension of~$L(\phi, x, \cdot)$, and the equation~\eqref{eq:rel-L-E-holo} with
  \begin{equation}\label{eq:expr-cpm-holo}
    c_f = c_f(1/2) = \frac{(-1)^{(k-\ell)/2} \Gamma(\ell)}{2^{\ell/2} \Gamma(\frac\ell2) \Gamma(\frac{k-\ell}2+1)\Gamma(\frac{k+\ell}2)}.
  \end{equation}
\end{proof}

\begin{rem}
  In the holomorphic case~$s_\phi = \ell/2$, another option would be to compute directly the Eichler integral~$\E(\phi, x, s)$ from the Fourier expansion of~$\phi$, instead of passing through the original form~$f$.
  This raises two difficulties. The first is that the explicit expression for the Whittaker function~$W_{\frac k2, \frac{\ell-1}2}$ is more involved, see~\cite[(9.237.3)]{GradshteynRyzhik2007}, although still elementary.
  The second difficulty is that~$\E_\phi(x)$ vanishes completely if~$k-\ell \equiv 2\pmod{4}$. Indeed one way to see this is the following: use the formulas~(9.237.3) and~(8.970.1) of~\cite{GradshteynRyzhik2007} to reduce to proving that~$F_{n,\ell}(v) := \sum_{m=0}^n \binom{n}{m}\frac{\Gamma(\ell/2+m)}{\Gamma(\ell+m)}v^m$ vanishes at~$v=-2$ for odd~$n$; by~(8.380.1) \emph{ibid.}, we have~$F_{n,\ell}(v) = \Gamma(\ell/2)^{-1} \int_0^1(1+vt)^n(t(1-t))^{\ell/2-1}\df t$, which obviously vanishes at~$v=-2$ for odd~$n$ by changing~$t$ to~$1-t$. Therefore, we cannot avoid having to switch to a lowered or raised form in that case.
  
  In practical situations, the formulas~\cite[(4.25)-(4.30)]{DukeEtAl2002} provide all the information one needs to translate data from~$\phi$ to~$f$.
\end{rem}

\subsection{Quantum modularity for the twisted central \texorpdfstring{$L$}{L} values}

By Lemma~\ref{lem:Lxs-extension}, we may now define
\begin{equation}
  L_\phi(x) := L(\phi, x, 1/2) \qquad \text{or} \qquad L_\phi^\pm(x) := L^\pm(\phi, x, 1/2), \label{eq:def-Lx}
\end{equation}
depending on whether~$s_\phi\equiv\frac k2\pmod 1$ or not.
We recall that~$\psi$ or~$f$ were defined in terms of~$\phi$ in~\eqref{eq:def-psi} or~\eqref{eq:def-fholom}.

\begin{thm}\label{th:Lx=qmf}
  When~$s_\phi \not\equiv\frac k2\pmod1$, the maps~$L_\phi^\pm(\cdot)$ are quantum modular forms for $\Gamma$ with multiplier $u_\gamma$ in the same sense as in Theorem~\ref{th:hgamma=qmf}, meaning that for all~$\gamma\in \Gamma$, the map
  \begin{equation}
    h^\pm_\gamma(x) := L^\pm_\phi(\gamma x) - j_\gamma(x)^k \chi(\gamma) L^\pm_\phi(x),\label{eq:def-hP}
  \end{equation}
  initially defined for $x\in C(\Gamma)\setminus \{\infty, \gamma^{-1}\infty\}$ extend to a $(1/2-\eps)$-H\"{o}lder continuous function on $\R\setminus \{\gamma^{-1}\infty\}$. More precisely, for~$x,x' \not \in [\gamma^{-1}\infty-\eps, \gamma^{-1}\infty+\eps]$, we have
  $$ \abs{h_\gamma^\pm(x) - h_\gamma^\pm(x')} \ll_{I, \gamma, \phi, \eps} \abs{x-x'}^{1/2-\eps} (1+\abs{x}+\abs{x'})^{O_\phi(1)}. $$
  
  The analogous statement holds for the map~$L_\phi(\cdot)$ when~$s_\phi \equiv \frac k2 \pmod{1}$. If $\phi$ is cuspidal the same statements hold with H\"{o}lder exponent $(1-\eps)$ in place of $(1/2-\eps)$.
\end{thm}

\begin{proof}
  Assume first that~$s_\phi\not\equiv\frac k2\pmod1$.
  By the relation~\eqref{eq:rel-L-E-1/2-maass} and linearity, we deduce
  \begin{equation}\label{eq:rel-hE-hL}
    h^\pm_\gamma(x) = c_\phi^\pm h_\gamma^\E(\phi; x) + c_\psi^\pm h_\gamma^\E(\psi; x),
  \end{equation}
  where~$h_\gamma^\E(\phi; x)$ satisfies~\eqref{eq:def-hP-E} and~$h_\gamma^\E(\psi; x)$ satisfies~\eqref{eq:def2-hpsi}. By Theorem~\ref{th:hgamma=qmf}, both extend to~$(1/2-\eps)$-Hölder continuous functions of~$x$ to~$\R\setminus\{\gamma^{-1}\infty\}$  ($(1-\eps)$-Hölder if~$\phi$ is cuspidal), which proves our claim in the case where~$\phi$ is not associated to a holomorphic form.
  
  The proof in the case~$s_\phi = \ell/2\equiv k/2\pmod 1$ is similar, invoking~\eqref{eq:rel-L-E-holo} in place of~\eqref{eq:rel-L-E-1/2-maass}.
\end{proof}

From Theorems~\ref{th:Lx=qmf} and \ref{th:hgamma-asymptotic}, it is straightforward to obtain the asymptotic behaviour of~$h_\gamma^+$ and~$h_\gamma^-$ as~$x\to \infty$ or~$x\to \gamma^{-1}\infty$.
However, the constants involved do not seem to admit a particularly simple expression in terms of~$\phi_\infty$ or~$L^\pm(x_0)$, so we refrain from carrying this out here.

Similarly as for Theorems~\ref{th:hgamma-fricke=qmf} and~\ref{th:hgamma-asymptotic-ext}, we obtain the following asymptotic behaviour on the generalized period function~$h_{\overline{\gamma}}^\E$ define in~\eqref{eq:def-hgamma-E-ext}.

\begin{thm}\label{th:Lx=gqmf}
  Let~$\gamma \in \SL(2, \R)$ and let~$\phi$ be as above. Assume that for some~$\eta\in \C$, we have~$(\overline\gamma \phi)(z) = \eta j_\gamma(z)^k \phi(z)$ for all~$z\in \H$. Then the analogue of Theorem~\ref{th:Lx=qmf} holds for the period functions defined for~$x\in C(\Gamma)\setminus \{\infty, -\gamma^{-1}\infty\}$ by
  \begin{align}
    \label{eq:hgammapm-gqmf}
    h_{\overline{\gamma}}^\pm(x) := {}& \overline{L_\phi^\pm(\gamma(-x))} - \eta j_\gamma(x)^k L_\phi^{\pm(-1)^k}(x), & (s_\phi\not\equiv\tfrac k2\pmod1), \\
    \notag
    h_{\overline{\gamma}}(x) := {}& \overline{L_\phi(\gamma(-x))} - \eta j_\gamma(x)^k L_\phi(x), & (s_\phi \equiv \tfrac k2 \pmod{1}).
  \end{align}
\end{thm}

Note the change of sign~$\pm$ in the first case, which, as we will see, is ultimately due to the conjugation.

\begin{proof}
  We first note that~$c_f\in\R$. Therefore the second equality immediately follows from~\eqref{eq:rel-L-E-holo}.
  
  Next, when $k$ is even, we see from~\eqref{eq:cphipsi-cossin-keven} and from the fact that~$t_\phi\in \R \cup i\R$, that~$c_\phi^\pm, c_\psi^\pm$ are real numbers. Therefore, the first equality follows by linearity from~\eqref{eq:rel-L-E-1/2-maass} and  Theorem~\ref{th:hgamma-fricke=qmf} when~$k$ is even.
  
  When~$k$ is odd, then as remarked in~\cite[p.~508]{DukeEtAl2002}, we have~$t_\phi\in \R$ always. Then a simple computation shows that
  $$ \overline{c_\phi^\pm} = c_\phi^{\mp}, \qquad \overline{c_\psi^\pm} = c_\psi^\mp. $$
  By~\eqref{eq:rel-L-E-1/2-maass} we deduce
  $$ h_{\bar\gamma}^\pm(x) = c_\phi^\mp h_{\bar\gamma}^\E(\phi, x) + c_\psi^\mp h_{\bar\gamma}^\E(\psi, x), $$
  and our claim follows from Theorem~\ref{th:hgamma-fricke=qmf}.
\end{proof}

\subsection{Functional equation for the additive twist \emph{L}-series}\label{sec:fe}

In this section only, we assume that~$\phi$ is cuspidal at $\infty$, in other words
$$ \phi_\infty \equiv 0. $$
We state a functional equation relating, for instance in the non-holomorphic case~$s_\phi \not\equiv \frac k2\pmod{1}$,
\begin{equation}\label{eq:prospective-FE}
  L^\pm(\phi, x, s) \qquad \text{ and } \qquad L^{\pm(-1)^k}(\phi, -\bar x, s),
\end{equation}
where~$x\in C(\Gamma)\setminus\{\infty\}$ is any cusp equivalent to~$\infty$, and given~$\gamma \in \Gamma$ with~$x=\gamma \infty$, $-\bar x$ is by definition
$$ -\bar x = \gamma^{-1}\infty \qquad (x = \gamma\infty). $$
This depends modulo~$1$ only the class of~$\gamma$ in~$\Gamma_\infty\backslash\Gamma/\Gamma_\infty$.
Note, as in Theorem~\ref{th:Lx=gqmf}, the change of sign in~\eqref{eq:prospective-FE} when~$k$ is odd.
The content of this section is not strictly related to quantum modularity, but it is convenient to include it at this point.

Define
\begin{equation}
  \begin{aligned}
    \eta := {}& -\chi(\gamma) i^k, \\
    \nu_\pm := {}& (-1)^{\floor{\frac k2 - \frac{1\mp \epsilon_\phi}4}}, & (s_\phi \not\equiv \tfrac k2\pmod{1}), \\
    \nu_f := {}& (-1)^{(k-\ell)/2}, & (s_\phi = \tfrac\ell2 \equiv \tfrac k2\pmod{1}).
  \end{aligned}\label{eq:FE-def-roots}
\end{equation}
If~$c\neq 0$ denotes the bottom-left coefficient of a matrix~$\gamma\in \Gamma$ such that~$x=\gamma \infty$, we let
\begin{align*}
  \Lambda^\pm(\phi, x, s) ={}& \Big(\frac{\abs{c}}\pi\Big)^s \Gamma\Big(\frac{s+it_\phi}2+\frac{1\mp\epsilon_\phi(-1)^k}4\Big) \Gamma\Big(\frac{s-it_\phi}2+\frac{1\mp \epsilon_\phi}4\Big) L^\pm(\phi, x, s),
  & (s_\phi \not\equiv \tfrac k2 \pmod{1}), \\
  \Lambda^f(\phi, x, s) ={}& \Big(\frac{\abs{c}}{2\pi}\Big)^s \Gamma(s_\phi-\tfrac12+s) L(\phi, x, s), & (s_\phi = \tfrac \ell2 \equiv \tfrac k2\pmod{1}).
\end{align*}

\begin{prop}\label{prop:Lpm-FE}
  With the above notations, for any cusp~$x\in C(\Gamma)\setminus\{\infty\}$ equivalent to~$\infty$, we have
  $$ \Lambda^\pm(\phi, x, s) = \eta \nu_\pm \Lambda^{\pm(-1)^k}(\phi, -\bar x, 1-s) $$
  if~$s_\phi \not\equiv \frac k2\pmod{1}$, and otherwise
  $$ \Lambda^f(\phi, x, s) = \eta \nu_f \Lambda^f(\phi, -\bar x, 1-s). $$
\end{prop}

The special case~$\chi=1$ is Proposition~3.3 of\cite{Nordentoft2021}, see also Lemma~1.2.(iv) of~\cite{MazurRubin} for an expression in terms of modular symbols.
The special case when~$\phi$ is a certain Eisenstein series of weight~$0$ is the functional equation of the Estermann function~\cite{Estermann1930}, which we mention below in Section~\ref{sec:examples-eisenstein}.
We focus on the non-holomorphic case~$s_\phi\not\equiv \frac k2\pmod{1}$, the complementary case being similar and comparatively simpler.
The proof is based on the argument of Hecke~\cite{Hecke36} in the case of holomorphic forms. The functional equation is deduced from the corresponding one for Eichler integrals, but thereafter one needs to prove a similar functional equation for integrals of Whittaker functions of the shape~\eqref{eq:def-Wmellin}.

The case~$k=0$ of Proposition~\ref{prop:Lpm-FE} was proved in~\cite[Appendix~A]{KoMiVa02}, see equations~(A.12) and~(A.13).
This uses explicit expressions for the integrals~$\Omega_{it_\phi}(0, s)$ in~\eqref{eq:Wmellin-gamma-half}, which we do not have for more general~$k$.

The computations we need for general~$k$ are done in~\cite[Section~8]{DukeEtAl2002}.
They correspond to the local functional equation at a real infinite place for the~$L$ function of an automorphic representation of~$\GL(2)$, which was worked out in~\cite[Chapter~5]{JacquetLanglands1970}; see also~\cite[Sections~2.7-2.8]{Godement1970}.

Here we sketch a more classical argument which passes through properties of hypergeometric functions. This circumvents the induction over~$k$ carried out in~\cite{DukeEtAl2002}, and highlights the relevance of $k$ being an integer in this context.
We start by the functional equation for Eichler integrals.
\begin{lemma}\label{lem:eichler-FE}
  For any~$x\in C(\Gamma)$ equivalent to~$\infty$, with denominator~$c$, we have
  $$ \E(\phi, x, s) = \eta c^{1-2s} \E(\phi, -\bar x, 1-s), $$
  where
  $$ \eta = \chi(\gamma) i^k. $$
\end{lemma}
\begin{proof}
  Recall the expression~\eqref{eq:def-eichlerint} and the fast decay of~$\phi$ near the cusps.
  Let~$x \in C(\Gamma)$ and~$\gamma \in \Gamma$ be such that~$x = \gamma \infty$. We write~$c = c_\gamma>0$.
  For~$z = \gamma^{-1}\infty + iy$, $y>0$, we have~$j(\gamma, z) = icy$. For those values of~$z$, we deduce~$u_\gamma(z) = \chi(\gamma) i^k$.
  Hence, for any~$s\in \C$,
  \begin{align*}
    \E(\phi, x, s) = {}& \int_{\gamma\infty}^\infty \phi(z) (\Im z)^{s-1/2} \df s(z) \\
    = {}& \int_{\gamma^{-1}\infty}^\infty \phi(\gamma z) (\Im \gamma z)^{s-1/2} \df s(z) \\
    = {}& \chi(\gamma) i^k \int_{\gamma^{-1}\infty}^\infty \phi(z) \abs{j(\gamma, z)}^{1-2s} (\Im z)^{s-1/2} \df s(z) \\
    = {}& \chi(\gamma) i^k c^{1-2s} \E(\phi, \bar x, 1-s).
  \end{align*}
\end{proof}

The previous Lemma can be applied to~$\psi = R_k \phi$, formally replacing~$k$ with~$k+2$. This yields
$$ \E(\psi, x, s) = -\eta c^{1-2s} \E(\phi, -\bar x, 1-s) $$
with the same value~$\eta = \chi(\gamma) i^k$ as in Lemma~\ref{lem:eichler-FE}.

We write, as in Section~\ref{sec:relperiod},
$$ L^\pm(\phi, x, s) = c_\phi^\pm(s) \E(\phi, x, s) + c_\psi^\pm(s) \E(\psi, x, s). $$
By Lemma~\ref{lem:eichler-FE}, using the value of~$\eta$ defined there, we find
$$ L^\pm(\phi, x, s) = \eta c^{1-2s} \big(c_\phi^\pm(s) \E(\phi, -\bar x, 1-s) - c_\psi^\pm(s) \E(\psi, -\bar x, 1-s)\big). $$
Proposition~\ref{prop:Lpm-FE} is therefore an immediate consequence of the following lemma.
\begin{lemma}\label{lem:cphipsi-FE}
  Define
  \begin{equation}
    \Psi_\pm(s) = \nu_\pm \pi^{2s-1} \frac{\Gamma(\frac{1-s-it_\phi}2+\frac{1\mp \epsilon_\phi(-1)^k}4) \Gamma(\frac{1-s+it_\phi}2+\frac{1\mp \epsilon_\phi}4)} {\Gamma(\frac{s+it_\phi}2+\frac{1\mp \epsilon_\phi(-1)^k}4) \Gamma(\frac{s-it_\phi}2+\frac{1\mp \epsilon_\phi}4)}\label{eq:FE-expr-W-symm}
  \end{equation}
  where~$\nu_\pm$ was defined in~\eqref{eq:FE-def-roots}.
  Then we have
  \begin{align*}
    c_\phi^\pm(s) ={}& \Psi_{\pm, \phi}(s) c_\phi^{\pm(-1)^k}(1-s), \\
    c_\psi^\pm(s) ={}& -\Psi_{\pm, \phi}(s) c_\psi^{\pm(-1)^k}(1-s).
  \end{align*}
\end{lemma}
\begin{proof}
  This statement is precisely the local functional equation at the real infinite place for the representation of~$\Gamma\backslash\GL(2, \A_\Q)$ associated with~$\phi$~\cite[Theorem~5.15]{JacquetLanglands1970}.
  We give here a proof in classical terms which relies on computations involving hypergeometric functions. These computations are carried out in the appendix.

  Recall the definition~\eqref{eq:def-Wmellin}. For~$s, \beta \in \C$ with~$0\leq \Re(\beta) <\frac12$ and~$\Re(s)> \abs{\Re(\beta)}$, and~$k\in\Z$, let
  $$ F_k^\pm(s, \beta) := \frac1{\Gamma(s-\beta)\Gamma(s+\beta)}\Big(-\Omega_\beta(-\tfrac k2, s) \pm \frac{\Gamma(\frac{1-k}2+\beta)}{\Gamma(\frac{1+k}2+\beta)} \Omega_\beta(\tfrac k2, s)\Big). $$
  \begin{lemma}\label{lem:EF-Fk}
    The following identity between meromorphic functions of~$s$ holds:
    $$ F_k^{\pm}(s, \beta) = -(-1)^{\floor{k/2 - (1\mp 1)/4)}} \frac{\Gamma(\frac{1-s+\beta}2+\frac{1\mp 1}4) \Gamma(\frac{1-s-\beta}2 + \frac{1\mp(-1)^k}4)}{\Gamma(\frac{s+\beta}2+\frac{1\mp(-1)^k}4)\Gamma(\frac{s-\beta}2+\frac{1\mp 1}4)} F_k^{\pm(-1)^k}(1-s, \beta). $$
  \end{lemma}
  \begin{proof}
    First note that
    $$ F_k^\pm(s, -\beta) = F_k^{\pm (-1)^k}(s, \beta). $$
    This is a restatement of~\cite[eq.~(8.34)]{DukeEtAl2002}; it follows from the invariance of~$\Omega_\beta(\alpha, s)$ by~$\beta \gets -\beta$ and by the complement formula.
    We now use Lemma~\ref{lem:hypergeom-FE} and the functions~$G_\pm$, $Q_\pm$ defined therein.
    By~\eqref{eq:rel-Wmellin-hypergeom}, we have
    \begin{align*}
      F_k^\pm(s, \beta) ={}& -\frac{1}{\Gamma(s+\frac{1-k}2)} \Big(\frac{\Gamma(s+\frac{1-k}2)}{\Gamma(s+\frac{1+k}2)} F(s-\beta, s+\beta; s+\tfrac{1+k}2; \tfrac12) \mp \frac{\Gamma(\frac{1-k}2+\beta)}{\Gamma(\frac{1+k}2+\beta)} F(s-\beta, s+\beta; s+\tfrac{1-k}2; \tfrac12) \Big) \\
      ={}& -\frac{1}{\Gamma(s+\frac{1-k}2)} G_\mp(s-\beta, s+\beta, s+\tfrac{1-k}2).
    \end{align*}
    We then compute
    \begin{align*}
      \frac{F_k^\pm(s, \beta)}{F_k^\pm(1-s, -\beta)} = {}& \frac{\Gamma(1-s+\frac{1-k}2)}{\Gamma(s+\frac{1-k}2)} \frac{G_\mp(s-\beta, s+\beta, s+\frac{1-k}2)}{G_\mp(1-s+\beta, 1+s-\beta, 1-s+\frac{1-k}2)} \\
      = {}& \frac{\Gamma(1-s+\frac{1-k}2)}{\Gamma(s+\frac{1-k}2)} Q_\mp(s-\beta, s+\beta, s+\tfrac{1-k}2).
    \end{align*}
    Using the explicit expression from Lemma~\ref{lem:hypergeom-FE} with~$n=k-1$, we get
    $$ Q_\mp(s-\beta, s+\beta, s+\tfrac{1-k}2) = (-1)^k 2^{2s-1} \frac{\Gamma(1-s+\beta)\Gamma(1-s-\beta)}{\Gamma(\frac{1+k}2-s)\Gamma(\frac{3-k}2-s)}
    \Big(1 \mp \frac{\cos(\pi(\beta+\frac{k}2))}{\cos(\pi(s+\frac{k}2))}\Big), $$
    and therefore
    \begin{align*}
      \frac{F_k^\pm(s, \beta)}{F_k^\pm(1-s, -\beta)} ={}& (-1)^k 2^{2s-1} \frac{\Gamma(1-s-\beta)\Gamma(1-s+\beta)}{\Gamma(\frac12-s+\frac k2)\Gamma(\frac12+s-\frac k2)} \Big(1 \mp \frac{\cos(\pi(\beta+\frac{k}2))}{\cos(\pi(s+\frac{k}2))}\Big)
    \end{align*}
    and the claimed formula follows by the complement and duplication formulas for the~$\Gamma$ function.
  \end{proof}

  We turn to the proof of Lemma~\ref{lem:cphipsi-FE}.
  Note that~$\Psi_{\pm, \psi}(s) = - \Psi_{\pm, \phi}(s)$, therefore it suffices to show either of the two formulas. We prove the second.
  From the definition of~$c_\psi^\pm(s)$, we have
  $$ c_\psi^\pm(s) = \frac{\pi^{s-1/2}}{4\Gamma(s+it_\phi)\Gamma(s-it_\phi)} \Big( -\Omega_{it_\phi}(-\tfrac k2, s) \pm \epsilon_\phi \frac{\Gamma(\frac{1-k}2+it_\phi)}{\Gamma(\frac{1+k}2+it_\phi)} \Omega_{it_\phi}(\tfrac k2, s)\Big). $$
  We express this as~$c_\psi^\pm(s) = \frac14 \pi^{s-1/2} F_k^{\pm\epsilon_\phi}(s, it_\phi)$, which yields
  $$ \frac{c_\psi^\pm(s)}{c_\psi^{\pm(-1)^k}(1-s)} = \pi^{2s-1} \frac{F_k^{\pm \epsilon_\phi}(s, it_\phi)}{F_k^{\pm(-1)^k \epsilon_\phi}(1-s, it_\phi)}. $$
  Using the previous lemma concludes the proof.
  
  When~$s_\phi = \ell/2 \equiv k/2\pmod{1}$, we do not reproduce the proof since it is similar and much simpler, since we have in that case the explicit expression~\eqref{eq:def-cfs}.
\end{proof}

\section{Examples}\label{sec:examples}

In this section and the following ones, we will be interested in applications in which the group~$\Gamma$ is an arithmetic group, and more precisely a Hecke congruence subgroup~$\Gamma_0(q)$ (see~\cite[p.~44]{Iw}).
This is our primary motivation for the above results.

Let~$q\in\Z_{>0}$, and denote by $\Gamma_0(q)$ the Hecke congruence subgroup of level $q$.
The associated set of cusps is given by~$C(\Gamma_0(q)) = \Q \cup \{\infty\}$.

Let~$\phi \in \A(\Gamma_0(q), \chi, k)$. If~$s_\phi\not\equiv\frac k2\pmod1$, then we have already defined in~\eqref{eq:def-Lxs-maass} and~\eqref{eq:def-Lx} the central~$L$-value~$L_\phi^\pm(x)$.
If~$s_\phi \equiv \frac k2\pmod{1}$ with~$k\geq 2$, on the other hand, we have merely defined by~\eqref{eq:def-Lxs-holo} and~\eqref{eq:def-Lx} the single value~$L_\phi(x) = U^+(\phi, x, 1/2)$. We now take advantage of the fact that the symmetry~$\varpi := \begin{psmallmatrix} -1 & 0 \\ 0 & 1\end{psmallmatrix}$ normalizes~$\Gamma_0(q)$, and more precisely, if~$\gamma = \begin{psmallmatrix} a & b \\ c & d\end{psmallmatrix}\in\Gamma_0(q)$, then~$\varpi \gamma \varpi = \begin{psmallmatrix} a & -b \\ -c & d \end{psmallmatrix} \in \Gamma_0(q)$. This easily implies that the map
$$ x\mapsto L_\phi(-x) $$
satisfies the same quantum modularity relation as~$L_\phi$ in Theorem~\ref{th:Lx=qmf}.
The same is therefore true for the maps
\begin{equation}\label{eq:def-Lxpm-holo}
  L^\pm_\phi(x) := \tfrac12(L_\phi(x) \pm L_\phi(-x))
\end{equation}
which are the even, resp.\ the odd part of~$L_\phi(\cdot)$. Thus~$L_\phi^\pm(x)$ is now defined in all cases, and clearly satisfies Theorem~\ref{th:Lx=qmf}. Moreover, we check that, setting
$$ h_{\bar\gamma}^\pm(x) := \overline{L_\phi^\pm(\gamma(-x))} - \eta j_\gamma(x)^k L_\phi^{\pm(-1)^k}(x), $$
we have
$$ h_{\bar\gamma}^\pm(x) = \frac12\big( h_{\bar\gamma}(x) + h_{\overline{\varpi\gamma\varpi}}(-x)\big). $$
Here we recall our convention that we pick the representative of~$\gamma$ in~$\PSL(2, \R)$ with~$c_\gamma\geq 0$ in the notation~$j_\gamma(x)$, and in particular~$j_{\varpi\gamma\varpi}(z) = cz - d$ (with~$c\geq 0$).
Therefore, under the hypotheses of Theorem~\ref{th:Lx=gqmf}, the equation~\eqref{eq:hgammapm-gqmf} also holds when~$s_\phi\equiv k/2\pmod{1}$ with the definition~\eqref{eq:def-Lxpm-holo}.

\subsection{Hecke--Maa\ss{} cuspidal newforms}\label{sec:HeckeMaass}
We will now review the theory of newforms. We refer to \cite[Section 8.5]{Iw}, \cite[Section 6.6]{Iw2} for a more detailed account. Let $\chi$ be a Dirichlet character modulo $q$ and define a character of $\Gamma_0(q)$ by 
$$\begin{psmallmatrix} a & b \\ c & d \end{psmallmatrix}\mapsto \chi(d), $$
which we also denote by $\chi$ (by slight abuse of notation). In this arithmetic setup we have a huge family of commuting linear operators acting on the space of automorphic functions $\A(\Gamma_0(q), \chi, k)$ for each~$k\in \Z_{\geq 0}$. These are the Hecke operators~\cite[eq.~(6.1)]{DukeEtAl2002} defined for $n\geq 1$ as
$$ T_n \phi(z)= \frac{1}{n^{1/2}}\sum_{ad=n}\chi(a) \sum_{0\leq b <d} \phi\left(\frac{az+b}{d}\right). $$
A Maa{\ss} cuspform $\phi\in \A(\Gamma_0(q), \chi, k)$  is called a \emph{Hecke--Maa{\ss} cuspform} if it is an eigenvector for all Hecke operators $T_n$ with $(n,q)=1$ and normalized so that $a(1)=1$. We furthermore, say that $\phi\in \A(\Gamma_0(q), \chi, k)$ with $k\geq 0$ is a \emph{Hecke--Maa{\ss} cuspidal newform} if it is not of the form 
$$ R_{k-2}\cdots R_\ell u(dz),$$ 
for a Hecke--Maa{\ss} form $u\in \A(\Gamma_0(q'), \chi, \ell)$ of level $q'<q$ with $dq'|q$ and weight $\ell\geq 0$ congruent to $ k$ modulo $2$. Notice in particular that if $\chi$ is a primitive Dirichlet character modulo $q$ then all Hecke--Maa{\ss} forms $\phi\in \A(\Gamma_0(q), \chi, k)$ with $k=0,1$ are new. Furthermore, if $\phi$ is Hecke--Maa{\ss} cuspidal newform not of weight $0$ nor $1$ then $\phi$ is holomorphic, meaning that $\phi(z)=y^k f(z)$ for some cuspidal holomorphic Hecke newform $f\in \mathcal{S}_k(\Gamma_0(q),\chi)$ of level $q$, weight $k$ and nebentypus $\chi$.

It is a consequence of {\lq}multiplicity one{\rq}~\cite[p.~520]{DukeEtAl2002} that a Hecke--Maa{\ss} cuspidal newform $\phi$ is automatically an eigenfunction for \emph{all} Hecke operators $T_n$ with $n\geq 1$. Let $\lambda_\phi(n)$ denote the Hecke eigenvalues of $\phi$ meaning that $T_n\phi =\lambda_\phi(n) \phi$ for $n\geq 1$. It follows from the properties of the Hecke operators that we have the Hecke relation 
\cite[(6.24)]{Iw}
\begin{equation}\label{eq:Heckerel}
  \lambda_\phi(mn)= \sum_{d|(m,n)} \mu(d) \chi_\phi(d) \lambda_\phi(m/d)\lambda_\phi(n/d),
\end{equation}
for all $m,n\geq 1$ (here it is crucial that $\phi$ is assumed to be a newform).
In particular $n\mapsto\lambda_\phi(n)$ is (weakly) multiplicative.
The Fourier coeffcients at $\infty$ of $\phi$ can be written as
$$ a_\phi(n)= \lambda_\phi(n)n^{-1/2}, $$
which implies that for $x\in C(\Gamma_0(q))\backslash \{\infty\}=\Q$, we have
$$ L^\pm (\phi,x,s) = \sum_{n\geq 1} \lambda_\phi(n) \left\{\begin{array}{c}\cos(2\pi n x)\\ i\sin(2\pi n x)\end{array}\right\} n^{-s}. $$
Here the definition is given by~\eqref{eq:def-Lxpm-holo} or by~\eqref{eq:def-Lx} depending on whether~$s\equiv \frac k2\pmod{1}$ or not.

One has the trivial point-wise bound $\lambda_\phi(n)\ll n^{1/2}$, which arises from the bound~$\abs{a(n)}\ll 1$. This last bound holds in the general setting for $\Gamma$ as discussed in Section~\ref{sec:setting} (see~\cite[Theorem~3.2]{Iw}).
The Ramanujan--Petersson conjecture predicts that $\lambda_\phi(n) \ll_\eps n^{\eps}$ for any $\eps>0$. This is a theorem due to Deligne~\cite{DeligneW1} in the case where~$s_\phi = \frac k2$, which means that~$\phi(z)=y^{k/2} f(z)$ with $f\in \mathcal{S}_k(\Gamma_0(q),\chi)$ a holomorphic Hecke cuspform. In general it is known by work of Kim and Sarnak~\cite{KimSarnak2003} that $|\lambda_\phi(n)|\leq d(n) n^{7/64}$.

Consider the involution 
$$\overline{W_{q,k}}:\A(\Gamma_0(q), \chi, k)\rightarrow \A(\Gamma_0(q), \chi, k),$$ 
defined by  
$$(\overline{W_{q,k}}\phi)(z) := \left(\frac{|z|}{-z}\right)^k\overline{\phi(1/(q\overline{z}))}. $$
Notice that $\overline{W_{q,k}}$ is not linear but skew-linear. It can be shown that $\overline{W_{q,k}}$ commutes with the Hecke operators and satisfies
\begin{equation}
  \label{eq:Wkraising} \overline{W_{q,k+2}} K_k=K_k\overline{W_{q,k}}, 
\end{equation}
By multiplicity one we conclude that for any Hecke--Maa{\ss} form $\phi$, we have 
$$\overline{W_{q,k}}\phi= \eta_\phi \phi,$$
for some $\eta_\phi$ of absolute norm $1$. In the terminology of Section \ref{sec:GQM} this means that a Hecke--Maa{\ss} form is automorphic of weight $k$ for the group $G_q$ generated by $\Gamma_0(q)$ and $\overline{W_q}$ where
$$W_q=\begin{psmallmatrix} 0 & -1 \\ q & 0 \end{psmallmatrix}, $$ 
is the Fricke matrix of level $q$. Thus we conclude from Theorem \ref{th:Lx=qmf} and Theorem \ref{th:Lx=gqmf} that the central values $L^\pm_\phi(x)=L^\pm(\phi,x,1/2)$ define (generalized) quantum modular forms, meaning that 
$$h^\pm_\gamma(x):=L^\pm_\phi(\gamma x)- j_\gamma(x)^k \chi(\gamma) L^\pm_\phi(x), $$
for $\gamma\in \Gamma_0(q)$, as well as 
$$h^\pm_{\overline{W_q}}(x):=\overline{L^\pm_\phi(1/qx)}- \eta_\phi (-\sgn(x))^k  L^{\pm(-1)^k}_\phi(x),$$
initially defined for $x\in \Q/ \{\gamma^{-1}\infty\}$ (resp.\ $x\in \Q/ \{0\}$) extends to a $(1/2-\eps)$-H\"{o}lder continuous function on $ \R/ \{\gamma^{-1}\infty\}$ (resp.\ $x\in \R/ \{0\}$). For applications to reciprocity formulae, we will need precise information on the discrepancy function $h^\pm_{\overline{W_q}}$ as follows.
\begin{prop}\label{prop:cuspidalhWq}
  For $\phi$ a Hecke--Maa{\ss} cuspidal newform of weight $k$ with $s_\phi\not\equiv k/2\pmod{1}$, we have as $x\rightarrow \infty$ 
  $$h^\pm_{\overline{W_q}}(x)= \begin{cases} -\eta_\phi \epsilon_\phi L(\phi,1/2)+O_{\phi, \eps}(|x|^{-1+\eps}), & \pm=+\\O_{\phi, \eps}(|x|^{-1+\eps}),& \pm =-,  \end{cases}$$
  for $k=0$, and
  $$h^\pm_{\overline{W_q}}(x)= \begin{cases}  -\eta_\phi\epsilon_\phi \frac{ \sinh (\pi  t_\phi)+i\epsilon_\phi}{\cosh (\pi  t_\phi)} L(\phi,1/2)+O_{\phi, \eps}(|x|^{-1+\eps}),& \pm =+,\\
  O_{\phi, \eps}(|x|^{-1+\eps}), & \pm=-\end{cases}$$
  for $k=1$.
  For $\phi$ a Hecke--Maa{\ss} cuspidal newform of weight $k$ with $s_\phi=k/2$ (i.e. holomorphic), we have
  $$h^\pm_{\overline{W_q}}(x)= \begin{cases} -i^k \eta_\phi L(\phi,1/2)+O_{\phi, \eps}(|x|^{-1+\eps}), & \pm=+\\O_{\phi, \eps}(|x|^{-1+\eps}),& \pm =-. \end{cases}$$
\end{prop}
\begin{proof}
  We proceed by using Lemma~\ref{lem:Lxs-extension} combined with Theorem \ref{th:hgamma-asymptotic-ext}. Implied constants are allowed to depend on~$\phi$ and~$\eps$. If $s_\phi=k/2$ (meaning that $\phi(z)=y^{k/2}f(z)$ for $f\in \mathcal{S}_k(\Gamma_0(q),\chi)$ a holomophic cupsidal Hecke newform) then we have
  $$ h^\pm_{\overline{W_q}}(x)=-i^k \eta_\phi  \frac{c_f (\E_\phi(0)\pm  \E_\phi(0))}{2}+O(|x|^{-1+\eps}). $$
  Using Lemma \ref{lem:Lxs-extension} in reverse, we get the wanted in the holomorphic case.
  Similarly we have for $k\in\{0, 1\}$ and~$s_\phi\not\equiv k/2\pmod{1}$
  \begin{align}
    h^\pm_{\overline{W_q}}(x)&=-i^k \eta_\phi  c_\phi^{\pm(-1)^k} \E_\phi(0)-i^{k+2} \eta_\psi  c_\psi^{\pm(-1)^k} \E_\psi(0)+O(|x|^{-1+\eps})\\
    &=-i^k \eta_\phi  \left(c_\phi^{\pm(-1)^k}\E_\phi(0)- c_\psi^{\pm(-1)^k}\E_\psi(0)\right)+O(|x|^{-1+\eps}),
  \end{align}
  using that $\eta_\phi=\eta_\psi$ which follows from \eqref{eq:Wkraising}. Recall that  
  $$ L^\pm_\phi (x)= c^\pm_\phi \E_\phi(x)+c^\pm_\psi \E_\psi(x), $$
  with $c^\pm_\phi, c^\pm_\psi$ defined as in \eqref{eq:expr-cpm}.
  Consider the matrix 
  $$C_\phi:=\begin{pmatrix} c_\phi^+ & c_\psi^+\\ c_\phi^- & c_\psi^- \end{pmatrix},$$
  which has determinant 
  \begin{align} 
    c_\phi^+c_\psi^- -c_\phi^-c_\psi^+
    = \sigma\Big(1 + \epsilon_\phi T_1 \Big)\Big(-1 - \epsilon_\phi T_2 \Big) 
    -\sigma\Big(1 - \epsilon_\phi T_1 \Big)\Big(-1 + \epsilon_\phi T_2 \Big)= -2\sigma\epsilon_\phi (T_1+T_2),
  \end{align}
  where 
  $$\sigma=\frac{\pi 2^{-2-k}}{\Gamma(\frac{1+k}4+i\frac{t_\phi}2) \Gamma(\frac{1+k}4-i\frac{t_\phi}2)\Gamma(\frac{3+k}4+i\frac{t_\phi}2) \Gamma(\frac{3+k}4-i\frac{t_\phi}2)},$$
  and 
  $$T_1=\frac{\cos(\frac\pi2(\frac{1+k}2+it_\phi))}{\sin(\frac\pi2(\frac{1+k}2-it_\phi))},\quad T_2=\frac{\sin(\frac\pi2(\frac{1+k}2+it_\phi))}{\cos(\frac\pi2(\frac{1+k}2-it_\phi))} .$$
  Since
  $$T_1+T_2= \frac{2\cos( \pi i t_\phi)}{\sin( \pi(\frac{k+1}{2} -i t_\phi))},$$
  we conclude that $C_\phi$ is invertible (since $it_\phi\notin 1/2+\Z$). 
  Now by simple linear algebra we get that
  $$c_\phi^\pm\E_\phi(x)- c_\psi^\pm\E_\psi(x)= \frac{c_\phi^\pm c_\psi^- +c_\phi^-c_\psi^\pm}{c_\phi^+c_\psi^- -c_\phi^-c_\psi^+}L^+_\phi(x) -
  \frac{c_\phi^+ c_\psi^\pm +c_\phi^\pm c_\psi^+}{c_\phi^+c_\psi^- -c_\phi^-c_\psi^+}L^-_\phi(x). $$
  For $x=0$ we have
  $$L^-_\phi(0)=-L^-_\phi(-0)=0,$$
  and
  $$L^+_\phi(0)=L(\phi,1/2), $$
  where $L(\phi,s)=\sum_{n\geq 1}\lambda_\phi(n)n^{-s} $ denotes the (standard) $L$-series of $\phi$. Using that 
  \begin{align} 
    c_\phi^\pm c_\psi^- +c_\phi^-c_\psi^\pm
    &= \sigma\Big(1 \pm \epsilon_\phi T_1 \Big)\Big(-1 - \epsilon_\phi T_2 \Big) 
    +\sigma\Big(1 - \epsilon_\phi T_1 \Big)\Big(-1 \pm \epsilon_\phi T_2 \Big)\\
    &= \sigma( -2(1\pm T_1T_2)+\epsilon_\phi (T_1-T_2)(1\mp 1)),
  \end{align}
  this yields 
  $$c_\phi^\pm\E_\phi(0)- c_\psi^\pm\E_\psi(0)= \frac{c_\phi^\pm c_\psi^- +c_\phi^-c_\psi^\pm}{c_\phi^+c_\psi^- -c_\phi^-c_\psi^+}L(\phi,1/2)=\frac{2(1\pm T_1T_2)-\epsilon_\phi (T_1-T_2)(1\mp 1)}{2\epsilon_\phi (T_1+T_2)}L(\phi,1/2) . $$
  For $k=0$, we have $T_1=T_2=1$ and for $k=1$ 
  $$ T_1=-\tan(\frac{\pi}{2}i t_\phi), \qquad T_2=\cot(\frac{\pi}{2}i t_\phi).$$
  Plugging this in gives the wanted. 
\end{proof}

\subsection{Eisenstein series}\label{sec:examples-eisenstein}

Let~$q_1, q_2 \in \Z_{>0}$ and~$\chi_i\pmod{q_i}$ be primitive characters, and~$k\in\{0, 1\}$ be such that~$(-1)^k = \chi_1\chi_2(-1)$.
For~$\Re(s)>1$ define the twisted Dirichlet series
$$ D_{\chi_1, \chi_2}(x, s) := \sum_{n\geq 1} (\chi_1 \ast \overline{\chi_2})(n) \e(nx) n^{-s}. $$
When~$x\in \Q$, orthogonality of additive characters yields an expression of~$D_{\chi_1, \chi_2}(x, s)$ as a linear combination of the Estermann function~$D_{1, 1}(x, s)$, which is known~\cite{Estermann1930} to have a meromorphic continuation to~$\C$ which is analytic on~$\C\setminus\{1\}$.
We deduce that the map~$s\mapsto D_{\chi_1, \chi_2}(x, s)$ extends to a meromorphic function of~$s$ which is analytic on~$\C\setminus\{1\}$, and also at~$1$ if~$x\not\in\Z$. When~$\chi_1 = \chi_2 = \1$, this is Estermann's function~\cite{Estermann1930}.

We are interested in the central value
$$ D_{\chi_1, \chi_2}(x) := D_{\chi_1, \chi_2}(x, 1/2). $$

\begin{thm}\label{th:Dchi=QM}
  The map~$D_{\chi_1, \chi_2}(\cdot)$ is a quantum modular form of weight~$k$ for the Hecke congruence group~$\Gamma_0(q_1q_2)$ with nebentypus~$\chi(\begin{psmallmatrix}\ast&\ast\\\ast&d\end{psmallmatrix}) := \chi_1\overline{\chi_2}(d)$, in the sense that for all~$\gamma = \begin{psmallmatrix} \ast&\ast\\c& d\end{psmallmatrix}\in \Gamma_0(q_1q_2)$ with~$c \neq 0$, the map
  $$ h_\gamma(x) := D_{\chi_1, \chi_2}(\gamma x) - \sgn(cx + d)^k \chi_1 \overline{\chi_2}(d) D_{\chi_1,\chi_2}(x) $$
  initially defined on~$\Q\setminus\{-d/c\}$, extends to a~$(1/2-\eps)$-continuous function of~$x\in \R\setminus\{-d/c\}$.
\end{thm}

\begin{proof}
  For~$\Re(s)>1$ and~$z\in \H$, define as in~\cite[Section~3.2]{Young2019} the Eisenstein series
  \begin{equation}\label{eq:eisensteindef} E_{\chi_1,\chi_2}(z, s) := \frac12 \sum_{\substack{c, d \in \Z \\ (c, d)=1}} \frac{(q_2y)^s \chi_1(c) \chi_2(d)}{\abs{cq_2 z + d}^{2s}}\Big(\frac{cq_2 z + d}{\abs{cq_2 z + d}}\Big)^k.\end{equation}
  It is proved in Section~3.2 \emph{ibid.} that~$E_{\chi_1, \chi_2}(\cdot, s)$ is a Maa\ss{} form for the Hecke congruence group~$\Gamma_0(q_1q_2)$, of weight~$k$, nebentypus~$\chi(\begin{psmallmatrix}\ast&\ast\\\ast&d\end{psmallmatrix}) := \chi_1\overline{\chi_2}(d)$, and eigenvalue~$s$.
  Moreover, it is proved in~\cite[Proposition~4.2]{Young2019} that the modified series
  \begin{equation}
    E_{\chi_1, \chi_2}^*(z, s) := \frac{(q_2/\pi)^s}{i^{-k} \tau(\chi_2)} \Gamma(s+\tfrac k2) L(2s, \chi_1\chi_2) E_{\chi_1, \chi_2}(z, s)\label{eq:Eisenstein-normalized}
  \end{equation}
  extends to a meromorphic function of~$s\in \C$ which is regular at~$s=1/2$.

  Define~$\phi = E_{\chi_1, \chi_2}^*(z, 1/2)$.
  By~\cite[Proposition~4.1]{Young2019} evaluated at~$s=1/2$ (see the proof of Proposition 4.1 \emph{ibid.}), we have the Fourier expansion~\eqref{fexp} with
  $$ a(n) = \frac{(\chi_1 \ast \overline{\chi_2})(n)}{\sqrt{n}} \qquad (n > 0). $$
  This is~$\lambda_{\chi_1, \chi_2}(n, 1/2)$ in the notations of~\cite[eq.~(4.2)]{Young2019}.
  We deduce that
  $$ D_{\chi_1, \chi_2}(x) = U^+_\phi(x) = \tfrac12(L^+_\phi(x) + L^-_\phi(x)) $$
  in the notation~\eqref{eq:def-Lx}. Theorem~\ref{th:Lx=qmf} yields the desired conclusion.
\end{proof}

The behaviour at infinity of~$h_\gamma$ can be spelled out using the expression~\cite[eq.~(4.1)]{Young2019} for the coefficients~$A_\phi, B_\phi$ in~\eqref{eq:def-phi-cuspidalpart}. We have explicitely
$$ B_\phi = \begin{cases} 1, & (q_1=q_2=1), \\ 0 & (\text{otherwise}), \end{cases} $$
$$ A_\phi = \begin{cases}
  \gamma_0 - \log(4\pi), & (q_1=q_2=1), \\
  q_2 \frac{i^k}{\sqrt{\pi} \tau(\chi_2)} \Gamma(\tfrac{1+k}2) L(1, \chi_2), & (q_1=1<q_2), \\
  0, & (q_1, q_2 > 1).
\end{cases} $$

\begin{prop}\label{prop:hgamma-eisenstein-asymptotic}
  As~$x\to\pm\infty$, we have an asymptotic expansion of the shape
  $$ h_\gamma(x) = \chi(\gamma)\big(A'_\pm \abs{x-x_0}^{1/2} + B'_\pm \abs{x-x_0}^{1/2}\log\abs{x-x_0} + C' \big) + O_{\chi_1, \chi_2, \gamma, \eps}(\abs{x}^{-1+\eps}). $$
  When~$k=0$, the coefficients are given in terms of~$A_\phi, B_\phi$ by
  \begin{align*}
    A'_\pm ={}& \frac{1\mp i}2 A_\phi - \Big(\frac{1\pm i}2 \frac\pi 2 + \frac{1\mp i}2 \log 2\Big) B_\phi, \\
    B'_\pm ={}& \frac{1\mp i}2 B_\phi, \\
    C' ={}& - \chi_2(-1) D_{\chi_1, \chi_2}(-x_0).
  \end{align*}
  When~$k=1$, the coefficients are
  \begin{align*}
    A'_\pm ={}& \sqrt{\pi} \frac{\pm 1 + i}2, \\
    B'_\pm ={}& 0, \\
    C'_\pm ={}& - D_{\chi_1, \chi_2}(x_0).
  \end{align*}
\end{prop}

\begin{proof}
  In the proof, implied constants are allowed to depend on~$\chi_j, \gamma$ and~$\eps$.
  Assume first~$k=0$. Then we use the expression
  $$ U^+(\phi, x, \tfrac12) = (\Delta_0(0, s))^{-1}(g_\phi^+(\tfrac12) \E_\phi(x) + g_\psi^+(\tfrac12)\E_\psi(x)) $$
  given at~\eqref{eq:rel-U-P}. In our case, we find
  \begin{equation*}
    g_\phi^+(\tfrac12) = \frac{2\pi^{3/2}}{\Gamma(1/4)^2}, \qquad g_\psi^+(\tfrac12) = \frac{-\pi^{3/2}}{\Gamma(3/4)^2}.
  \end{equation*}
  In particular, we deduce
  $$ h_\gamma(x) = \sqrt{\pi}\Big(\frac1{\Gamma(1/4)^2}h_\gamma^\E(\phi, x) - \frac1{2\Gamma(3/4)^2} h_\gamma^\E(\psi, x)\Big). $$
  By~\cite[eq.~(4.29)]{DukeEtAl2002}, the coefficients~$A_\psi, B_\psi$ in the cuspidal expansion~\eqref{eq:def-phi-cuspidalpart} for~$\psi=R_0\phi$ satisfy
  $$ A_\psi = \tfrac12 A_\phi + B_\phi, \qquad B_\psi = B_\phi. $$
  By Theorem~\ref{th:hgamma-asymptotic}, we deduce as~$x\to\pm\infty$ the asymptotic formula
  $$ h_\gamma(x) = \chi(\gamma) \Big(A'_\pm \abs{x-x_0}^{1/2} + B'_\pm \abs{x-x_0}^{1/2}\log\abs{x-x_0} + C'\Big) + O(\abs{x}^{-1+\eps}), $$
  where the coefficients are
  \begin{align*}
    C' ={}& -\sqrt{\pi}\Big(\frac1{\Gamma(1/4)^2} \E_\phi(x_0) + \frac1{2\Gamma(3/4)^2}\E_\psi(x_0)\Big), \\
    B'_\pm ={}& \Big(\frac{\Upsilon_{0,\pm}(1/2)}{\Gamma(1/4)^2} - \frac{\Upsilon_{2,\pm}(1/2)}{4\Gamma(3/4)^2}\Big) B_\phi, \\
    A'_\pm ={}& \Big(\frac{\Upsilon_{0,\pm}}{\Gamma(1/4)^2} - \frac{\Upsilon_{2,\pm}(1/2)}{4\Gamma(3/4)^2}\Big)A_\phi + \Big(\frac{\Upsilon_{0,\pm}'(1/2)}{\Gamma(1/4)^2} - \frac{\Upsilon_{2,\pm}(1/2)}{2\Gamma(3/4)^2} - \frac{\Upsilon_{2,\pm}'(1/2)}{4\Gamma(3/4)^2}\Big)B_\phi.
  \end{align*}
  Using the expressions~\eqref{eq:Upsilon-0-2}, we get
  \begin{align*}
    A'_\pm ={}& \frac{1\mp i}2 A_\phi - \Big(\frac{1\pm i}2 \frac\pi2 + \frac{1\mp i}2 \log 2\Big) B_\phi, &
    B'_\pm ={}& \frac{1\mp i}2 B_\phi.
  \end{align*}
  Moreover, using the expressions~\eqref{eq:system-Ephipsi-L}, \eqref{eq:rel-Uplusminus} and~\eqref{eq:Wmellin-gamma-half}, we find
  $$ C' = -\epsilon_\phi U^+(-x_0). $$
  By~\cite[Proposition~4.1, eq.~(4.2)]{Young2019}, we have~$\epsilon_\phi = \chi_2(-1)$, which gives the claimed value of~$C'$.
  
  In the complementary case~$k=1$, we have~$B=0$ necessarily. By~\eqref{eq:rel-L-E-holo}, we express~$U^+(\phi, x, \tfrac12) = L_\phi(x) = \frac1{\sqrt{2\pi}} \E_\phi(x)$, and therefore
  $$ h_\gamma(x) = \frac1{\sqrt{2\pi}} h_\gamma^\E(\phi, x). $$
  By Theorem~\ref{th:hgamma-asymptotic} and the expression~\eqref{eq:Upsilon-0-2}, we obtain
  $$ h_\gamma(x) = \chi(\gamma)\big(\sqrt{\pi} \frac{\pm 1 + i}2 \abs{x-x_0}^{1/2} - L_\phi(x_0)\big) + O(\abs{x}^{-1+\eps}). $$
\end{proof}

Similarly, we can consider the action of the Fricke involution $\overline{W_q}$ as defined in the previous section on the Eisenstein series $\phi=E^\ast_{\chi_1,\chi_2}(z,1/2)$. Consider the operator 
\begin{equation}\label{eq:def-op-X}
  (\overline{X}f)(z)= \overline{f(-\overline{z})},
\end{equation}
for $f:\Hb\rightarrow \C$. Then $\overline{W_q}$ is the composition of $\overline{X}$ with the usual action of the Fricke matrix $\begin{psmallmatrix} 0 & -1 \\q & 0\end{psmallmatrix}$. Using the definition of the Eisenstein series \eqref{eq:eisensteindef} one sees directly that
\begin{equation}
  \overline{X} E_{\chi_1,\chi_2}(z,s)= \chi_2(-1) E_{\overline{\chi_1},\overline{\chi_2}}(z,\overline{s}),   
\end{equation}
where $E_{\chi_1,\chi_2}(z,s)$ denotes the uncompleted newform Eisenstein series. Combining this with the calculations of \cite[Section 9.2]{Young2019} and the functional equation \cite[Proposition 4.2]{Young2019}, one gets for $s=1/2$ 
\begin{align}
  \overline{W_q} E_{\chi_1,\chi_2}(z,1/2)&=\overline{X} \left[  \chi_1(-1)\left(\frac{z}{|z|}\right)^k E_{\chi_2,\chi_1}(z,1/2)\right]\\
  &=(-1)^k \chi_1(-1) \chi_2(-1) \left(\frac{z}{|z|}\right)^k  E_{\overline{\chi_2},\overline{\chi_1}}(z,1/2)\\
  &=\frac{q_2^{1/2}\tau(\overline{\chi_1})L(1,\chi_1\chi_2)}{q_1^{1/2}\tau(\chi_2)L(1,\overline{\chi_1\chi_2})}\left(\frac{z}{|z|}\right)^k E_{\chi_1,\chi_2}(z,1/2).    
\end{align}
By linearity and taking into account the conjugation in~\eqref{eq:def-op-X}, we deduce
$$ \overline{W_q} E_{\chi_1, \chi_2}^*(z, 1/2) = (-1)^k \frac{\tau(\overline{\chi_1}) \tau(\chi_2)}{(q_1q_2)^{1/2}} j_{W_q}(z)^k E_{\chi_1, \chi_2}^*(z, 1/2). $$
This implies that $D_{\chi_1,\chi_2}(\cdot)$ defines a generalized quantum modular form in the sense of Theorem~\ref{th:Lx=gqmf}, with the value~$\eta$ being
$$\eta_{\chi_1,\chi_2} = (-1)^k \frac{\tau(\overline{\chi_1}) \tau(\chi_2)}{(q_1q_2)^{1/2}}, $$
which indeed does satisfy $|\eta_{\chi_1, \chi_2}|=1$ as should be the case. We have the following behaviour at infinity.
\begin{prop}\label{prop:hgammabar-eisenstein-asymptotic}
  The map
  $$ h_{\overline{W_q}}(x) := \overline{D_{\chi_1,\chi_2}(1/(q_1q_2x))} - \eta_{\chi_1,\chi_2} \sgn(x)^k D_{\chi_1,\chi_2}(x) $$
  satisfies, as~$x\to\pm\infty$, the asymptotic expansion
  $$ h_{\overline{W_q}}(x) = \eta_{\chi_1,\chi_2}\big(A'_\pm \abs{x}^{1/2} + B'_\pm \abs{x}^{1/2}\log\abs{x} + C' \big) + O_{\chi_1, \chi_2, \eps}(\abs{x}^{-1+\eps}), $$
  with $A'_\pm,B'_\pm$ as in Proposition \ref{prop:hgamma-eisenstein-asymptotic} and 
  $$
  C' =\begin{cases}
    - \chi_2(-1) L(\chi_1,1/2)L(\chi_2,1/2),& k=0,\\
    - L(\chi_1,1/2)L(\chi_2,1/2), & k=1.
  \end{cases}
  $$
\end{prop}
\begin{proof}
  We proceed by using Theorem \ref{lem:Lxs-extension} combined with Theorem \ref{th:hgamma-asymptotic-ext} noting that 
  $$D_{\chi_1,\chi_2}(0)=L(\chi_1,1/2)L(\chi_2,1/2).$$
\end{proof}

\subsubsection*{The Estermann function}

The Estermann function~\cite{Estermann1930}
$$ D(x) := D_{1, 1}(x), $$
is a particular case, which was studied in~\cite{Bettin16}.
As a consequence of Proposition~1 of~\cite{Bettin16}, the map
$$ h_D(x) := D(-1/x) - D(x) $$
extends to a $(1/2-\eps)$-Hölder-continuous function on~$\R\setminus\{0\}$.
Theorem~\ref{th:Lx=qmf} recovers this statement by a different proof.

Regarding asymptotic formulas, in~\cite[Corollary~9]{Bettin16} (see also~\cite{BettinDrappeau2019}, the formula above~(9.10)), it is shown that for~$x\in \Q$, $x\to 0$ with~$\sgn(x) = \pm 1$,
$$ h_D(x) = A'_\pm \abs{x}^{-1/2} + B'_\pm \abs{x}^{-1/2}\log\abs{x} - D(0) + O(x), $$
with
$$ A'_\pm = -\frac{1\pm i}2(\gamma - \log(8\pi)) - \frac{-1\pm i}2 \frac \pi2, \qquad B'_\pm = \frac{1\pm i}2. $$
It is easy to check that this matches the expansion given by Proposition~\ref{prop:hgamma-eisenstein-asymptotic}.

\subsubsection*{The representation function as sums of two squares}

Let~$r(n) := |\{(a, b)\in Z^2, n = a^2 + b^2\}|$.
It is classically known (see~\cite[eq.~(1.51)]{IwaniecKowalski2004}) that
$$ r(n) = \tfrac14 (\1 \ast \chi_4)(n) $$
where~$\chi_4$ is the non-trivial Dirichlet character modulo~$4$. For~$x\in \Q$ and~$\Re(s)>1$, define
$$ R(x, s) := \sum_{n\geq 1} \frac{r(n) \e(nx)}{n^s}. $$
Then we have
$$ R(x, s) = \tfrac14 D_{1,\chi_4}(x, s) $$
for~$\Re(s)>1$, which gives the holomorphic continuation to~$\C\setminus\{1\}$ of~$R(x, \cdot)$. Let
$$ R(x) := R(x, \tfrac12). $$
Then by Theorem~\ref{th:Dchi=QM}, the map~$R$ is a quantum modular form of weight~$1$ for~$\Gamma_0(4)$.
Since~$\Gamma_0(4)$ is generated by~$\{\begin{psmallmatrix} 1 & 1 \\ & 1\end{psmallmatrix}, \begin{psmallmatrix} 1& \\ 4 & 1\end{psmallmatrix}\}$, this amounts to saying that the map
$$ h_R(x) := R(\tfrac{x}{4x+1}) - \sgn(x+\tfrac14) R(x) $$
extends to a~$(1/2-\eps)$-Hölder continuous map on~$\R\setminus\{-1/4\}$. 

More precisely, in this case, the spectral parameter~$s_\phi = 1/2$ is half the weight~$k=1$, so~$\phi$ is related to a holomorphic form of weight~$1$.
By applying Lemma~\ref{lem:Lxs-extension} in the second case with~$\ell=k=2$, we deduce
$$ R(x) = \frac1{4\sqrt{2\pi}} \E(\phi, x). $$
We have also~$\phi_\infty(y) = (\sqrt{\pi}/2) y^{1/2}$, so that upon applying Theorem~\ref{th:hgamma-asymptotic} and computing~$\Upsilon_{1, \pm}(\frac12)$, we obtain
$$ h_R(x) = \frac{\pi}{16}(\pm 1 + i)\abs{x+1/4}^{1/2} - iR(-1/4) + O_\eps(\abs{x}^{-1+\eps}) $$
as~$x \to \pm \infty$.

Similarly since $R(\cdot)$ has real nebentypus, we get  that  
$$ h_{W_4}(x):= R(\tfrac{-1}{4x})-i\sgn(x) R(x),$$
originally defined for $x\in \Q^\times $ extends to a $(1/2-\eps)$-H\"{o}lder contuous map in $\R^\times$. Here $W_4=\begin{psmallmatrix}0 & -1 \\ 4 & 0 \end{psmallmatrix}$ is the Fricke matrix of level $4$. Furthermore, we have
$$ h_{W_4}(x) =\frac{\pi}{16}(1\mp i)\abs{x}^{1/2} - \frac 14 \zeta(1/2)L(1/2, \chi_4) + O_\eps(\abs{x}^{-1+\eps}), $$
as $x\rightarrow \pm \infty$ using that  $\eta_{1,\chi_4} = -i$.

\section{Normal distribution in the cuspidal case}\label{sec:normal-distrib}

In this section, we work with the full modular group and consider cuspidal forms:
$$ \Gamma = \SL(2, \Z), \qquad \phi_\infty \equiv 0. $$
If we assume that $\phi$ is a Hecke--Maa{\ss} cusp form then the coefficients~$a(n)$ are real numbers, \emph{cf.}~\cite[eq.~(6.6)]{DukeEtAl2002}.
Using Theorem~\ref{th:Lx=qmf}, we can answer completely the question of the asymptotic statistical distribution of the values of~$L_\phi(x)$ as~$x$ runs over the set of rationals of denominators at most~$Q$, and~$Q\to \infty$. This generalizes the case~$G=\SL(2, \Z)$ of~\cite[Theorem~1.1]{Nordentoft2021} and~\cite[Theorem~2.3]{BettinDrappeau2019}, which were concerned with holomorphic forms (see~\cite[Theorem~1.11]{PetridisRisager2018}, \cite[Theorem~C]{LeeSun2019} for related results for congruence groups).

Given~$Q\geq 1$, we let
$$ \Omega_Q := \{ x\in \Q\cap (0, 1), \den(x) \leq Q\}, $$
where~$\den(x)$ denotes the reduced denominator of~$x$, with~$\den(0) = 1$.
Define
$$ \PP_Q, \quad \EE_Q, \quad \VV_Q $$
to be the uniform probability measure over~$\Omega_Q$, and the associated expectation and variance.

\subsection{Distributional result: characteristic function}

Let~$S = \begin{psmallmatrix}&-1\\1&\end{psmallmatrix}$, and for~$x\in \Q\setminus\{0\}$ and~$\phi$ a Maa\ss{} cusp form, define
$$ g_\phi^\pm(x) := -h_S(\phi, x) = L_\phi^\pm(x) - L_\phi^\pm(-1/x). $$
By Theorem~\ref{th:Lx=qmf}, this map extends to a~$(1-\eps)$-Hölder continuous function on~$\R\setminus\{0\}$.
By Theorem~\eqref{eq:hgamma-E-asymptotic} and~\eqref{eq:rel-hE-hL}, the maps~$g_\phi^\pm$ admit limits at~$0$ on both sides, and at~$\pm\infty$, which implies that~$h_\phi^\pm$ is bounded.
By Euclid's algorithm and the~$1$-periodicity of~$L_\phi^\pm$, we deduce
$$ L_\phi^\pm(x) = \sum_{j=1}^r g_\phi^\pm((-1)^{j-1} T^{j-1}x) + L_\phi^\pm(0), $$
where~$T:(0, 1) \to [0, 1)$, $Tx = \{1/x\}$ is the Gauss map, and~$r\geq 0$ is the least integer such that~$T^j x = 0$.

In particular, the boundedness of~$g^\pm$ along with the worst-case estimate for the complexity of the Euclidean algorithm~\cite[Corollary~L, p.~360]{KnuthII} immediately implies the following rough but useful bound.
\begin{lemma}\label{lem:bound-Lpm-pointwise}
  For~$x\in\Q$, we have~$L_\phi^\pm(x) = O(1+\log(\den(x)))$.
\end{lemma}

We now consider~$r\in\N$, $\phi_1, \dotsc, \phi_r$ distinct cuspidal Hecke--Maa\ss{} forms.
Note that~$L^+_{\phi_j}(x) \in \R$ and~$L^-_{\phi_j}(x) \in i\R$. Define the non-normalized vector
\begin{equation}\label{eq:def-cCx}
  \cV(x) = \cV_{\phi_1, \dotsc, \phi_r}(x) := \Big(L_{\phi_1}^+(x), i^{-1}L_{\phi_1}^-(x), \dotsc, L_{\phi_r}^+(x), i^{-1}L_{\phi_r}^-(x) \Big) \in \R^{2r}, \qquad (\den(x)>1)
\end{equation}
where we set~$\cV_\phi^\pm(0) = 0$. Our aim is to show that~$\cV_\phi^\pm(x)$ converges, under a suitable normalization, to an $2r$-dimensional Gaussian random vector. We consider the linear form on~$\C^{2r}$ defined at
$$ \bt = (t_1^+, t_1^-, \dotsc, t_r^+, t_r^-) \in \C^{2r} $$
by the value
$$ \L(\bt; x) := \cV(x) \bt^T = \sum_{1\leq j \leq r} \Big(t_j^+ L_{\phi_j}^+(x) + t_j^- i^{-1} L_{\phi_j}^-(x)\Big), \qquad x \in \Omega_Q, $$

We define
$$ \Psi(\bt) := \EE_Q \Big( \exp\big\{\L(\bt, x)\big\} \Big). $$

\begin{prop}\label{prop:Lpm-charfun}
  For some~$\delta, t_0>0$, and some maps~$U, V$ holomorphic on~$B := \{\bt\in \C^{2r}, \abs{t_j}<t_0\}$, the estimate
  \begin{equation}
    \Psi(\bt) = \exp\Big\{ U(\bt) \log Q + V(\bt) + O(Q^{-\delta}) \Big\}\label{eq:quasi-mom}
  \end{equation}
  holds for~$\bt\in B$. The implied constant may depend on~$(\phi_j)$. Moreover, for some row vector~$\mu \in \C^r$ and some~$d\times d$ symmetric non-negative matrix~$\Sigma$, we have
  \begin{equation}
    \Psi(\bt) = \exp\Big\{ (\bt \mu^T) \log Q + \frac12 \bt \Sigma \bt^T \log Q + O(Q^{-\delta} + O(\|\bt\|^3 \log Q + \|\bt\| + Q^{-\delta}) \Big\}.\label{eq:quasi-mom-taylor}
  \end{equation}
\end{prop}

\begin{proof}
  We apply Theorem~\cite[Theorem~3.1]{BettinDrappeau2019} with~$m\gets2$, $d\gets 2r$, and the maps~${\bm \phi}_j$, $j\in\{1, 2\}$ given by
  $$ {\bm \phi}_j(x) = \big(g_{\phi_1}^+((-1)^{j-1}x), i^{-1} g_{\phi_1}^-((-1)^{j-1}x), \dotsc \big) \in \R^{2r}. $$
  By Theorem~\ref{th:Lx=qmf}, the map~${\bm\phi}_j$ is~$(1-\eps)$-Hölder continuous, and by Theorem~\ref{th:hgamma-asymptotic} and the relation~\eqref{eq:rel-hE-hL}, it is bounded. Moreover, by~\eqref{eq:holder-bound-hE}, for all~$n\geq 1$ and~$x, x' \in [0, 1]$ we have
  $$ \abs{g_{\phi_j}^\pm(1/(n+x)) - g_{\phi_j}^\pm(1/(n+x'))} = \abs{g_{\phi_j}^\pm(-n-x) - g_{\phi_j}^\pm(-n-x')} \ll_\eps n^C \abs{x-x'}^{1/2-\eps} $$
  for some~$C = C(\phi_j)$, and similarly for~$g_{\phi_j}^\pm(-1/(n+x))$.
  Therefore, the hypotheses~(1)-(3) of~\cite[Theorem~3.1]{BettinDrappeau2019} are satisfied with
  $$ \kappa_0 = 1/2-\eps, \qquad \alpha_0 = 3, \qquad \lambda_0 = \tfrac12\min\{1/C(\phi_1), \dotsc, 1/C(\phi_j)\}. $$
  This yields the claimed estimate for~$\bt \in B\cap \R^{2r}$.
  The holomorphicity of~$U, V$ and the fact that~\eqref{eq:quasi-mom} holds throughout~$B$ comes from the boundedness of~${\bm\phi}_j$, as in~\cite{BaladiVallee2005}.
  The estimate~\eqref{eq:quasi-mom-taylor} follows from a Taylor expansion at~$\bt = 0$.
\end{proof}

By reasonning similarly as~\cite[eqs.~(9.7)-(9.9)]{BettinDrappeau2019}, we have
\begin{equation}\label{eq:rel-mu-Sigma-mom}
  \begin{aligned}
    \mu_{(\phi_j)} ={}& \lim_{Q\to \infty} \frac{\EE_Q(\cV(x))}{\log Q}, &
    \Sigma_{(\phi_j)} ={}& \lim_{Q\to\infty} \frac{\EE_Q(\cV(x)^T \cV(x))}{\log Q}.
  \end{aligned}
\end{equation}

Our next task is to evaluation these two quantities.

\subsection{Computation of the first and second moment}

\begin{prop}\label{prop:moments-1-2}
  Let $\phi$ and~$\psi$ be Hecke--Maa{\ss} cusp forms of level $1$, normalized so that $a_\phi(1)=a_\psi(1)=1$, and~$\eps_1, \eps_2 \in\{\pm 1\}$.
  Then for some~$\delta>0$ and some~$b_{\phi, \pm}, c_{\phi, \psi, \pm} \in\C$, we have
  \begin{align}
    &\EE_Q(L_\phi^{\pm}(x)) = O_\phi(Q^{-1/3}), \label{eq:mom1} \\
    &\EE_Q(L_\phi^{\pm}(x)^2) = \pm L(\sym^2\phi, 1) \log Q + b_{\phi,\pm} + O_\phi(Q^{-\delta}), \label{eq:mom2abs}\\
    &\EE_Q(L_\phi^{\pm}(x) L_\psi^{\pm}(x)) = c_{\phi, \psi, \pm} + O_{\phi, \psi}(Q^{-\delta}) \qquad \text{ if $\phi \neq \psi$.} \label{eq:mom2rel}
  \end{align}
  Here~$L(\sym^2 \phi, s)$ is the symmetric square $L$-function of $\phi$~\cite[Chapter~8.2]{Iw}.
\end{prop}
The value~$L(\sym^2 \phi, 1)$ can be expressed in terms of the appropriate Petersson inner product~$\langle \phi, \phi \rangle$, see~\cite[eq.~(5.101)]{IwaniecKowalski2004} for details.
Note that~$L_\phi^\pm(x)^2 = \pm |L_\phi^\pm(x)|^2$, since in our case~$a(n)\in\R$ for all~$n$. Note also that~$\EE_Q(L_\phi^+(x) L_\psi^-(x)) = 0$ because of the symmetry~$x\gets -x$.
The value~$c_{\phi, \psi, +}$ could be expressed in terms of the Rankin-Selberg $L$-value~$L(\phi\times\psi, 1)$ and the constants in the functional equation, but we will not use it.

When~$\phi$ is associated to a holomorphic form, Proposition~\ref{prop:moments-1-2} corresponds, up to the size of the error terms, to the first few cases~$m+n\leq 2$ of Theorem~5.10 of~\cite{Nordentoft2021}.
An analogous estimate, but for \emph{fixed} denominator, is established in~\cite[Chapter~9]{BlFoKoMiMiSa18}. The argument is very different and much more difficult. We cannot, however, quote them in a straightforward way, due to the coprimality condition on the denominator. This is likely to change in the near future~\cite{Wu}.

Assuming Proposition~\ref{prop:moments-1-2} for a moment, we readily deduce the computation of the first and second moment of the random vector~$\cV_{\phi_1, \dotsc, \phi_r}(x)$, and therefore the values of~$\mu$ and~$\Sigma$ in~\eqref{eq:rel-mu-Sigma-mom}.
\begin{cor}\label{cor:value-mu-Sigma}
  Let $\phi_1, \dotsc, \phi_r$ be distinct Hecke--Maa{\ss} cusp forms of level $1$. Then
  \begin{equation}\label{eq:values-mu-Sigma}
    \mu_{(\phi_j)} = (0, \dotsc, 0), \qquad \Sigma_{(\phi_j)} = \operatorname{Diag}(\sigma_1^2, \sigma_1^2, \sigma_2^2, \sigma_2^2, \dotsc, \sigma_r^2, \sigma_r^2),
  \end{equation}
  where~$\sigma_j = \sqrt{L(\sym^2\phi_j, 1)}$.
\end{cor}
\begin{proof}
  This is immediate for~$\mu$ using~\eqref{eq:mom1}. The coefficients of the matrix~$\Sigma$ are indexed by pairs~$(\phi_j, \pm)$. The coefficient of indices~$((\phi, \eps_1), (\psi, \eps_2))$ is given by
  $$ \lim_{Q\to\infty} \frac{\EE_Q(i^{(\eps_1+\eps_2-2)/2} L_\phi^{\eps_1}(x) L_\psi^{\eps_2}(x))}{\log Q}. $$
  By~\eqref{eq:mom2abs}, this is~$L(\sym^2\phi, 1)$ if~$(\phi, \eps_1) = (\psi, \eps_2)$, which corresponds to diagonal elements, and otherwise this is~$0$.
\end{proof}

The proof of Proposition~\ref{prop:moments-1-2} starts with the following approximate functional equations.

\begin{lemma}\label{lem:AFE}
  For some numbers~$\mu_{\phi, \psi}^\pm$ and some smooth functions~$V^\pm:\R_+\to\C$, we have
  \begin{align*}
    L^\pm_\phi(x) L^\pm_\psi(x) ={}& \pm \sum_{m, n\geq 1} \frac{\lambda_\phi(m)\lambda_\psi(n)}{\sqrt{mn}} \big(C^\pm_{m,n}(x) + \mu^\pm_{\phi, \psi} C^\pm_{m,n}(\bar x)\big) V^\pm\Big(\frac{\pi mn}{\den(x)^2}\Big).
  \end{align*}
  Here
  $$ C^+_{m,n}(x)=\cos(2\pi m x)\cos(2\pi n x), \qquad C^-_{m,n}(x)= \sin(2\pi m x)\sin(2\pi n x), $$
  and writing~$x = a/q$ in lowest terms, we have
  $$ \bar x = -d/q \pmod{1} \qquad (ad \equiv 1\pmod{q}). $$
  Moreover the functions~$V^\pm$ can be chosen to satisfy, for any fixed~$a\in\N$ and~$A>0$,
  $$ (V^\pm)^{(a)}(y) \ll_{a, A, \phi} y^{-a} (1+y)^{-A}. $$
  Moreover~$\mu_{\phi, \psi}^\pm = 1$ if~$\phi = \psi$.
\end{lemma}
\begin{proof}
  This is obtained by an argument identical to~\cite[Lemma~9.1]{DukeEtAl2002}, specialized at~$s=1/2$, see also Theorem~5.3 and Proposition 5.4 of~\cite{IwaniecKowalski2004}.
\end{proof}

Taking expectations of the periodic exponential will give rise to Ramanujan sums.
\begin{lemma}\label{lem:ramanujan-sums}
  Given a map~$f:\N\to\C$ and~$m, n\in\Z$, we have
  $$ \EE_Q(C^\pm_{m,n}(x) f(\den(x))) = \frac1{2\abs{\Omega_Q}} \sum_{\eta\in\{-1, 1\}} \eta^{(1\mp 1)/2}\sum_{\substack{q, d \geq 1 \\ qd \leq Q}} d\mu(q) \1_{d\mid m-\eta n}f(qd). $$
\end{lemma}
\begin{proof} This easily follows from~\cite[eq.~(3.2)]{IwaniecKowalski2004}. \end{proof}

\begin{proof}[Proof of Proposition~\ref{prop:moments-1-2}]
  All implied constants in the proof will be allowed to depend on~$\phi, \psi$ and~$\eps$.
  First note that we have~$\EE_Q(L_\phi^-(x))=0$ by symmetry~$x\gets -x$.
  Consider~$s\in\C$ with~$\Re(s)>1$. Writing out the central~$L$-value and using Lemma~\ref{lem:ramanujan-sums}, we have
  $$ \EE_Q(L^+(\phi, x, s)) = \frac1{\abs{\Omega_Q}} \sum_{\substack{d, q\geq 1 \\ dq\leq Q}} d^{1-s} \mu(q) \sum_{n\geq 1} \lambda(nd) n^{-s}. $$
  We expand~$\lambda(nd)$ using the Hecke relations~\eqref{eq:Heckerel}, to get
  \begin{align*}
    \EE_Q(L^+(\phi, x, s)) ={}& \frac1{\abs{\Omega_Q}} \sum_{dq \leq Q} \sum_{r\mid d} d^{1-s} \mu(q) \mu(r) \chi(r) \lambda(d/r) \sum_{\substack{n\geq 1 \\ r\mid n}} \lambda(n/r) n^{-s} \\
    ={}& \frac{L(\phi, s)}{\abs{\Omega_Q}} \sum_{rdq \leq Q} r^{1-2s} d^{1-s} \mu(q) \mu(r) \chi(r) \lambda(d).
  \end{align*}
  At this point we may set~$s=1/2$ and bound trivially
  $$ \sum_{rdq\leq Q} d^{1/2}\mu(q) \mu(r)\chi(r) \lambda(d) \ll Q^{3/2+\theta+\eps}, $$
  by virtue of the bound~$\lambda(d) \ll d^{\theta+\eps}$, where~$\theta\leq 7/64$ is a bound towards the Ramanujan--Petersson conjecture. Since~$\abs{\Omega_Q}\asymp Q^2$, we deduce our first claim~\eqref{eq:mom1} since~$\EE_Q(L_\phi^+(x)) \ll Q^{-1/2+\theta+\eps} \ll Q^{-1/3}$.
  
  We switch to the computation of~$\EE_Q(L_\phi^+(x) L_\psi^+(x))$, the case of~$L_\phi^-(x) L_\psi^-(x)$ being similar.
  We first motivate the upcoming arguments. Using the functional equation and orthogonality, we expect
  $$ \EE_Q(L_\phi^+(x)L_\psi^+(x)) \approx \frac1{Q^2} \sum_{q\leq Q} \underset{\substack{mn\ll Q^2 \\ q\mid m-n}}{\sum_{m} \sum_{n}} \lambda_\phi(m) \lambda_\psi(n)
  \approx \frac1{Q^2} \underset{mn\ll Q^2}{\sum_{m} \sum_{n}} \lambda_\phi(m) \lambda_\psi(n) \tau_Q(m-n), $$
  where~$MN\ll Q^2$ and~$\tau_Q(h) := \1 \ast \1_{[1, Q]}(h)$.
  Considering~$m\leq n$ and~$m$ fixed, the sum over~$(n, q)$ is an instance of the shifted convolution of~$\lambda_\psi$ with a modified divisor function~$\tau_Q$, for which we have very good error terms using Voronoï summation with Dirichlet's hyperbola method, as we will do here. This was used for instance in~\cite[Lemma~8.1]{Pitt}.
  If we wanted an estimate for a single large~$q$, as in the conjectures of Mazur-Rubin~\cite[section~4]{MazurRubin}, then we would be faced with the shifted convolution of~$\lambda_\phi$ with~$\lambda_\psi$ with a large shift, which is much more delicate~\cite[chapter~9]{BlFoKoMiMiSa18}.

  First we smooth out the sharp cutoff on the denominator of~$x$ in the expectation~$\EE_Q$. To do this, we let~$Y = Q^\delta$ for some fixed~$\delta\in(0, 1/10]$ to be set later, and we bound, using Lemma~\ref{lem:bound-Lpm-pointwise},
  $$ \frac1{\abs{\Omega_Q}} \sum_{\substack{x\in \Omega_Q \\ Q(1-Y^{-1}) \leq \den(x) \leq Q}} \big|L_\phi^+(x)L_\psi^+(x)\big| \ll (\log Q)^2 Y^{-1}. $$
  Let~$W_0:\R \to \R_+$ be such that~$\1_{[0, 1-Y^{-1}]} \leq W_0 \leq \1_{[0, 1]}$ and~$\|W_0^{(j)}\|_\infty \ll_j Y^j$. We get
  $$ \EE_Q(L_\phi^+(x) L_\psi^+(x)) = \EE_Q(L_\phi^+(x) L_\psi^+(x) W_0(\den(x)/Q)) + O(Q^\eps Y^{-1}). $$
  Denoting for convenience~$V = V^+$ and~$\mu = \mu_{\phi, \psi}^+$, we have by Lemma~\ref{lem:ramanujan-sums}
  $$ \EE_Q(L_\phi^+(x) L_\psi^+(x) W_0(\den(x)/Q)) = \frac{1+\mu}{2\abs{\Omega_Q}} \sum_\pm \sum_{q, d \geq 1} d \mu(q) W_0\Big(\frac{qd}{Q}\Big) \sum_{\substack{m, n \geq 1 \\ d\mid m\pm n}} \frac{\lambda_\phi(m) \lambda_\psi(n)}{\sqrt{mn}} V\Big(\frac{\pi m n}{(qd)^2}\Big). $$
  The condition~$qd\leq Q$ was dropped due to redundance with the support of~$W_0$.
  The diagonal~$\pm = -1$, $m=n$ contributes
  $$ \cD_Q := \frac{1+\mu}{2\abs{\Omega_Q}} \sum_{q\geq 1} \varphi(q) W_0\Big(\frac{q}{Q}\Big) \sum_{n\geq 1} \frac{\lambda_\phi(n)\lambda_\psi(n)}n V\Big(\frac{\pi n^2}{q^2}\Big), $$
  where~$\varphi(q) = \sum_{d\mid q} d \mu(q/d)$ is the Euler totient function.
  By the complex analytic properties of the Rankin-Selberg convolution~$\phi\times\psi$~\cite[Chapters~5.11-5.12]{IwaniecKowalski2004}, we have, for some sufficiently small~$\delta>0$,
  $$ \sum_{n\geq 1} \frac{\lambda_\phi(n)\lambda_\psi(n)}n V\Big(\frac{\pi n^2}{q^2}\Big) = \Res_{u=0} \bigg( \Big(\frac{q^2}\pi\Big)^u {\check V}(u) L(\phi\times\psi, 1+2u)\bigg) + O(q^{-\delta}) $$
  Here~$\check V(u)$ is the Mellin transform of~$V$, and by~\cite[eq.~(5.13)]{IwaniecKowalski2004}, we have~$\Res_{u=0}(\dotsb) = P(\log q)$, where~$P(X)$ is either:
  \begin{itemize}
  \item a degree~$1$ polynomial with leading coefficient~$L(\sym^2\phi, 1)$, if~$\psi = \phi$,
  \item a constant polynomial with value~$L(\phi\times \psi, 1)$ if~$\psi \neq \phi$.
  \end{itemize}
  Summing this over~$q, d$, we find
  $$ \cD_Q = \frac{1+\mu}2 P(\log Q) - 1 + O(Y Q^{-\delta}). $$

  Introducing a partition of unity as in~\emph{e.g.} \cite[Lemma~2.27]{BlFoKoMiMiSa18}, we have
  $$ \EE_Q(L_\phi^+(x)L_\psi^+(x)) = \cD_Q + O(Q^\eps Y^{-1} + (\log Q)^4 \sup \abs{A_\pm(Q_1, D, M, N)}), $$
  where
  \begin{align*}
    A_\pm(Q_1, D, M, N) := {}& \frac1{\abs{\Omega_Q}} \sum_{q, d \geq 1} W\Big(\frac q{Q_1}\Big) W\Big(\frac dD\Big) W_0\Big(\frac{qd}{Q}\Big) d \mu(q) \times \\
    &\qquad \times \sum_{\substack{m, n \geq 1 \\ d\mid m\pm n \neq 0}} \frac{\lambda_\phi(m) \lambda_\psi(n)}{\sqrt{mn}} W\Big(\frac mM\Big) W\Big(\frac nN\Big) V\Big(\frac{\pi m n}{(qd)^2}\Big),
  \end{align*}
  the function~$W$ is smooth and supported inside~$[1/2, 2]$, and the supremum is over
  \begin{equation}
    \begin{aligned}
      Q_1, D, M, N \gg{}& 1, & M\leq{}& N, \\ Q_1 D \ll{}& Q, & MN \ll{}& Q^2.
    \end{aligned}\label{eq:mom2rel-conditions}
  \end{equation}
  We focus on~$A_+$, the case~$A_-$ being similar.
  We have by the triangle inequality
  $$ A_+(\dotsb) \ll \frac{D}{Q^2M^{1/2}} \sum_{q\asymp Q_1} \sum_{d\asymp D} \sum_{m\asymp M} \abs{\lambda_\phi(m)} \bigg| \sum_{n\equiv -m\pmod{d}} \frac{\lambda_\psi(n)}{\sqrt{n}} W\Big(\frac nN\Big) V\Big(\frac{\pi m n}{(qd)^2}\Big)\bigg|. $$
  We estimate the sum over~$n$ by Voronoï summation~\cite[Lemma~2.21]{BlFoKoMiMiSa18}, getting
  $$ \sum_{n\equiv -m\pmod{d}} \frac{\lambda_\psi(n)}{\sqrt{n}} W\Big(\frac nN\Big) V\Big(\frac{\pi m n}{(qd)^2}\Big) = \frac{\omega}{d\sqrt{N}} \sum_\pm \sum_{r\mid d} \frac1r \sum_{n\geq 1} \lambda_\psi(n) S(m, \mp n; r) \widetilde{W}_{\pm, N}(n/r^2), $$
  where~$S(m, \mp n; r)$ is the Kloosterman sum~\cite[eq.~(1.56)]{IwaniecKowalski2004} and by~\cite[Lemma~2.23]{BlFoKoMiMiSa18} we have
  $$ \abs{\widetilde{W}_{\pm, N}(y)} \ll N Q^\eps \times \begin{cases} Q^{-100}, & (y\leq Q^\eps/N), \\ (Ny)^{-2\vartheta-\eps} & (y\leq Q^\eps/N). \end{cases} $$
  Using the Weil's bound for Kloosterman sums~\cite[Corollary~11.12]{IwaniecKowalski2004}, we get
  \begin{align*}
    \sum_{n\equiv -m\pmod{d}} \frac{\lambda_\psi(n)}{\sqrt{n}} W\Big(\frac nN\Big) V\Big(\frac{\pi m n}{(qd)^2}\Big) 
    \ll Q^{-10} + Q^\eps\frac{\sqrt{N}}{d} \sum_{r\mid d} \frac1{\sqrt{r}} \sum_{n\ll Q^\eps r^2/N} \abs{\lambda_\psi(n)} (n, r)^{1/2} \Big(\frac{r^2}{nN}\Big)^{2\vartheta}.
  \end{align*}
  We use the Rankin-Selberg bound~$\sum_{n\leq X} \abs{\lambda_\psi(n)}^2 \ll X$, the bound~$\vartheta\leq 7/64 < 1/4$ and the elementary bound
  $$ \sum_{n\ll Q^\eps r^2/N} (n, r) \Big(\frac{r^2}{nN}\Big)^{4\vartheta} \ll Q^\eps \frac{r^2}N. $$
  Combining these estimates by Cauchy's inequality, we get
  $$ \sum_{n\equiv -m\pmod{d}} \frac{\lambda_\psi(n)}{\sqrt{n}} W\Big(\frac nN\Big) V\Big(\frac{\pi m n}{(qd)^2}\Big) \ll Q^{-10} + Q^\eps \frac1{d\sqrt{N}} \sum_{r\mid d} r^{3/2} \ll Q^\eps \sqrt{d/N}, $$
  and summing this over~$m$, $d$ and~$q$, we conclude our first bound
  \begin{equation}\label{eq:mom2rel-bound-voronoi}
    A_+(Q_1, D, M, N) \ll Q^{-2+\eps} Q_1 D^{5/2} M^{1/2} N^{-1/2}.
  \end{equation}
  The same bound holds for~$A_-$. This bound will be acceptable for~$D$ smaller than~$N^{2/3}$, however the smooth average over~$d$ allows for a more effective estimate as soon as~$D>\sqrt{N}$ by switching divisors. We arrange
  $$ A_+(D, Q_1, M, N) \ll \frac{D}{Q^2M^{1/2}} \sum_{q\asymp Q_1} \sum_{m\asymp M} \abs{\lambda_\phi(m)} \abs{B(m, q)}, $$
  where
  $$ B(m, q) := \sum_n \frac{\lambda_\psi(n)}{\sqrt{n}} W\Big(\frac nN\Big) \sum_{n\equiv -m\pmod d} W\Big(\frac dD\Big) \frac dD W_0\Big(\frac{qd}Q\Big) V\Big(\frac{\pi m n}{(qd)^2}\Big). $$
  Changing~$d'\gets (m+n)/d$, we get
  $$ B(m, q) = \sum_{d'\geq 1} \sum_{\substack{n\geq 1 \\ n\equiv -m\pmod{d'}}} \frac{\lambda_\psi(n)}{\sqrt{n}} W_2\Big(\frac nN\Big), $$
  where
  $$ W_2(x) := W(x) W\Big(\frac{m+Nx}{d'D}\Big)\frac{m+Nx}{d'D} W_0\Big(\dotsb\Big) V\Big(\dotsb\Big). $$
  Note that here we have used the fact that~$m+n\neq 0$.
  The support condition here implies~$d' \leq 8N/D$ for any non-zero term in the sum.
  The Leibniz differentiation rule yields~$\abs{W_2^{(j)}(x)} \ll_{\eps, j} Q^\eps Y^j$.
  Having this at hand, we may now apply Voronoï summation as above and bound trivially the dual sum. We get
  $$ \abs{B(m, q)} \ll Y^{O(1)} Q^\eps \sum_{d'\ll N/D} \sqrt{d'/N} \ll Y^{O(1)} Q^\eps N D^{-3/2}. $$
  Summing over~$m$ and~$q$ and using~\eqref{eq:bound-an-rough}, we obtain our second bound
  \begin{equation}
    \label{eq:mom2rel-bound-voronoi2}
    A_+(D, Q_1, M, N) \ll Y^{O(1)} Q^{-2+\eps} Q_1 D^{-1/2} M^{1/2} N.
  \end{equation}
  The same bound holds for~$A_-$, since we have removed the diagonal~$m=n$.
  Under the constraints~\eqref{eq:mom2rel-conditions}, we have
  \begin{align*}  
    Q^{-2}\min(Q_1 D^{5/2} M^{1/2} N^{-1/2}, Q_1 D^{-1/2} M^{1/2} N) 
    \leq \min(D^{3/2} N^{-1}, D^{-3/2} N^{1/2}) 
    \leq Q^{-1/4}.
  \end{align*}
  Therefore, we deduce
  $$ \EE_Q(L_\phi^+(x)L_\psi^+(x)) = \cD_Q + O_\eps(Q^\eps(Y^{-1} + Y^{O(1)} Q^{-1/4})), $$
  and by choosing the exponent~$\delta>0$ of~$Y=Q^\delta$ sufficiently small in terms of the implicit constant in~$O(1)$, we get~\eqref{eq:mom2abs} with~$b_{\phi, +} = P(0)-1$ in case~$\psi=\phi$, and~$c_{\phi, \psi, +} = \frac{1+\mu_{\phi, \psi}^+}2$ if~$\psi\neq \phi$. The analogous computation for~$L_\phi^-(x)^2$ yields an estimate with a possibly different constant term~$b_{\phi, -}$ and leading coefficient~$P'(0) = -L(\sym^2\phi, 1)$.
\end{proof}

\subsection{Distributional result: normal distribution}

We can now use the estimates~\eqref{eq:quasi-mom-taylor} and Corollary~\ref{cor:value-mu-Sigma} in the classical Levy continuity theorem, we obtain the convergence of~$\cV^\pm(x)$ to a complex Gaussian law. We recall the values~$\sigma_j := \sqrt{L(\sym^2 \phi_j, 1)} > 0$ and the matrix~$\Sigma$ defined in Corollary~\ref{cor:value-mu-Sigma}.

\begin{cor}\label{cor:joint-CLT}
  Given a measurable set~$R\subset \R^{2r}$ with measure~$0$ boundary, we have
  $$ \PP_Q\Big(\frac{\cV(x)}{\sqrt{\log\den(x)}} \in R \Big) \to \PP(\cN(0, \Sigma) \in R) \qquad (Q\to\infty), $$
  where~$\cN(0, \Sigma)$ represents a random centered Gaussian vector in~$\R^{2r}$ with covariance matrix~$\Sigma$.
\end{cor}
\begin{proof}
  By partial summation, it suffices to prove the same limiting statement for
  $$ \PP_Q\Big(\frac{\cV(x)}{\sqrt{\log Q}} \in R \Big). $$
  This is an immediate consequence of Proposition~\ref{prop:Lpm-charfun} and Levy's continuity theorem, see \emph{e.g.}~\cite[corollary~2.8]{Stromberg}.
\end{proof}

With more work, we could obtain an error term of the shape~$O(1/\sqrt{\log Q})$, with an application of a multi-dimensional version of the Berry-Esseen theorem. We haven't, however, found a ready-to-use statement in the literature, therefore we restrict to a quantitative statement. Note however that when $r=1$ we obtain a quantitative statement as in \cite[Thm. 2.3]{BettinDrappeau2019} using similar arguments.

We believe that this statement holds without the extra average over~$q$:

\begin{conj}
  In the context and notations of Corollary~\ref{cor:joint-CLT}, we have
  $$ \Big\{\frac{\cV(a/q)}{\sqrt{\log q}}: a\in (\Z/q\Z)^\times \Big) \to \PP(\cN(0, \Sigma) \in R) \qquad (q\to\infty). $$
\end{conj}

\subsection{Moment calculations}

Using the complex moments estimate in Proposition~\ref{prop:Lpm-charfun}, and the computation of the first two moments in Corollary~\ref{cor:value-mu-Sigma}, we readily deduce an estimate with power-saving for all moments, see~\cite[eq.~(2.12)]{BaladiVallee2005}.

\begin{prop}\label{prop:moments-level1}
  Let~$\phi_1, \dotsc, \phi_r$ be distinct Hecke-Maa\ss{} cusp forms,~$k_1, \dotsc, k_r, \ell_1, \dotsc, \ell_r \in \N_{>0}$. Then there exists~$\delta>0$ and a polynomial~$P$ such that
  $$ \EE_Q\Big( \prod_{1\leq j \leq r} L_{\phi_j}^+(x)^{k_j} L_{\phi_j}^-(x)^{\ell_j} \Big) = P(\log Q) + O(Q^{-\delta}). $$
  The implied constant may depend on~$(\phi_j), (k_j)$ and~$(\ell_j)$.
  When~$k_j$ and~$\ell_j$ are even for all~$j$, then~$P$ has degree~$\frac12\sum_j(k_j+\ell_j)$ and leading coefficient
  $$ \prod_{1\leq j \leq r} m_{k_j} m_{\ell_j} \sigma_j^{k_j+\ell_j}, \qquad m_k := (k-1)!! = \frac{k!}{2^{k/2}(k/2)!}. $$
  If at least one of the exponents~$k_j, \ell_j$ is odd, then~$P$ has degree strictly less than~$\frac12\sum_j(k_j+\ell_j)$.
\end{prop}
\begin{proof}
  See \cite{Hwang} for a proof when~$r=1$. The general case is analogous.
\end{proof}

\section{Arithmetic applications}\label{sec:arith-applications}

Let $\phi$ be a Hecke--Maa{\ss} newform (not necessarily cuspidal) of level $q$, weight $k$ and nebentypus $\chi_\phi$. Then we will be studying the following twisted $L$-functions; 
\begin{equation}\label{eq:Ldirichlet}L(\phi, \chi, s):=\sum_{n\geq 1} \lambda_\phi(n)\chi(n)n^{-s},\end{equation}
where $\lambda_\phi(n)$ are the Hecke-eigenvalues of $\phi$ (normalized so that the Ramanujan--Petersson conjecture predicts $\lambda_\phi(n)\ll_\eps n^{\eps}$) and $\chi$ is a primitive Dirichlet character mod $c$. This Dirichlet series admits analytic continuation and functional equation relating $s\leftrightarrow 1-s$ (see e.g.  \cite[Section 2.2]{BlFoKoMiMiSa18} in the case of prime conductor $c=p$ not dividing $q$). If $(c,q)=1$, then this is exactly the finite part of the $L$-function of the automorphic representation $\pi_\phi\otimes \chi$, where $\pi_\phi$ is the automorphic representation corresponding to $\phi$. In this case we will use the notation $L(\phi\otimes \chi, s)$. These $L$-functions have been studied in many contexts from the analytic point of view and have many interesting algebraic aspects as well (e.g. if $\phi$ is holomorphic of weight 2 corresponding to an elliptic curve), see \cite[Section 1.2]{BlFoKoMiMiSa18}, \cite[Chapter~14.8]{IwaniecKowalski2004} and the references therein. In the monograph \cite{BlFoKoMiMiSa18} the full power of the approximate functional equation-approach is explored including deep input from spectral theory and algebraic geoemtry. In particular they show how one can use a second moment computation for the family 
$$\{L(\phi, \chi, 1/2): \chi \modulo p\},$$ 
with a power saving as $p\rightarrow \infty$ to obtain non-vanishing for a positive proportion using mollification (among other applications). 

In this section we will show some surprising applications of the above results to the family of  $L$-functions \eqref{eq:Ldirichlet}. On the one hand, we will obtain certain reciprocity formulas for the twisted second moment of~\eqref{eq:Ldirichlet} generalizing \cite{Conrey07}, \cite{Bettin16}, \cite{Nordentoft20}. They can be seen as the simplest incarnations of spectral reciprocity formulas, see \cite{AndersenKiral18},\cite{BlomerKhan19},\cite{BlomerLiMiller19}. The second application are to certain computations of wide moments as have been explored in other contexts \cite{Be17}, \cite[Corollary 1.9]{Nordentoft2021}, \cite{Nordentoft21}, \cite{Nordentoft21.2}. 

The starting point for each of the two applications is the Birch--Stevens formula which relates additive and multiplicative twists.

\subsection{The Birch--Stevens formula}\label{sec:BirchStevens}

The central values of $L$-series of additive twists of holomorphic cusp forms of weight~2 are known as modular symbols, introduced by Birch and Manin~\cite{Birch1971, Manin1973}. Modular symbols have been used extensively in the study of the arithmetics of $L$-functions due to the {\it Birch--Stevens formula}.
\begin{prop}\label{BS}
  Let $\phi$ be a Hecke--Maa{\ss} newform (not neccesarily cuspidal) of level $q$, weight $k\in \Z_{\geq 0}$ and  neben-typus $\chi_\phi$. Then for $\chi$ a Dirichlet character mod $c$, we have
  \begin{align}
    \nu(\phi, \chi^*,c/c(\chi)) L(f, \chi^*,1/2)=\sum_{a\in (\Z/c\Z)^\times} \overline{\chi}(a )(L^+(\phi,a/c,1/2)+L^-(\phi,a/c,1/2)),
  \end{align}
  and 
  \begin{align}
    \label{birchstevens2}L^{\pm}(\phi,a/c, 1/2)=\frac{2}{\varphi(c)} \sum_{\substack{\chi \modulo c\\ \chi(-1)=\pm 1}} \nu(\phi, \chi^*,c/c(\chi)) L(\phi, \chi^*,1/2)\chi(a),
  \end{align}
  where $\chi^*\modulo c(\chi)$ denotes the unique primitive character that induces $\chi$ and the arithmetic weight $\nu$ defined by
  \begin{equation} \label{eq:weightBS}\nu(\phi, \chi, n):= \tau(\overline{\chi})\sum_{n_1n_2n_3=n}\chi_\phi(n_1)\chi(n_1) \mu(n_1)\overline{\chi}(n_2) \mu(n_2) \lambda_\phi(n_3)n_3^{1/2}, \end{equation}
  with $\tau(\overline{\chi})$ the Gau{\ss} sum of $\overline{\chi}$.
  
  If $\phi=E^*_{\chi_1,\chi_2}$ is the newform Eisenstein series from~\eqref{eq:Eisenstein-normalized}, then we have 
  $$\lambda_\phi(n)= (\chi_1\ast \overline{\chi_2})(n)=\sum _{d|n} \chi_1(d) \overline{\chi_2}(n/d),$$
  and we simply write $ \nu(\phi, \chi, n)=\nu(\chi_1\ast \overline{\chi_2}, \chi, n)$.
\end{prop}
\begin{proof}
  The proof is a straightforward adaption of the proof of \cite[Proposition 6.1]{Nordentoft2021} recalling that we have the Hecke relations~\eqref{eq:Heckerel} since $\phi$ is assumed to be a newform. 
\end{proof}
We will make a few comments on the arithmetic weights $\nu$. Note that if $\chi$ is primitive modulo $c$ then the arithmetic weight is simply given by $\nu(\phi, \chi^*,c/c(\chi))=\tau(\overline{\chi})$. In particular the weight is of absolute norm $c^{1/2}$ in this case. In general, if we have $\lambda_\phi(n)\ll_\eps n^{\theta+\eps}$ then one gets the following bound
\begin{equation} \label{eq:boundnu}\nu(\phi, \chi^*,c/c(\chi))\ll_\eps c(\chi)^{1/2} (c/c(\chi))^{1/2+\eps} (c/c(\chi))^\theta \ll_\eps c^{1/2+\theta+\eps},  \end{equation}
which is $O_\eps(c^{1/2+\eps})$ assuming the Ramanujan--Petersson conjecture. Finally we observe that we can express $\nu(\phi, \chi, n)$ in terms of a triple convolution as follows
$$ \tau(\overline{\chi})\cdot [( \chi_\phi\chi\mu)\ast ( \overline{\chi}\mu)\ast (\lambda_\phi |\cdot|^{1/2})](n).$$

If we restrict to prime conductor and level $1$, we get the following more pleasant form.

\begin{cor}\label{cor:BSprimelevel}
  Let $\phi$ be a Hecke--Maa{\ss} newform of level $1$ and let $p$ a prime number. Then 
  \begin{align}
    L^\pm(\phi,a/p,1/2)=\frac{2}{p-1}\sum_{\substack{\chi \modulo p \text{ \rm primitive},\\ \chi(-1)=\pm 1}}\tau(\overline{\chi})L(\phi, \chi,1/2)\chi(a)+O_{\phi, \eps}(p^{\theta-1+\eps}) ,
  \end{align}
  for $\eps>0$, where $\theta= \frac{7}{64}$ is the best bound towards the Ramanujan--Petersson conjecture for Maa{\ss} forms due to Kim and Sarnak \cite{KiSa03}. 
\end{cor}

\subsection{Applications to reciprocity formulae} %%%%%%%%%%%%%%%%%%%%%%%%%%%%%%%%%%

The starting point is the following unpublished paper of Conrey \cite[Theorem 10]{Conrey07}, in which a reciprocity relation satisfied by twisted second moment of Dirichlet $L$-functions was discovered. Here the terminology {\lq\lq}reciprocity{\rq\rq} refers to the cosmetic similarity with quadratic reciprocity; one relates the arithmetics of the seemingly unrelated finite fields $\mathbb{F}_q$ and $\mathbb{F}_\ell$ for  primes $q\neq \ell$. In the case of Conrey, the reciprocity relation relates the following two objects
\begin{align}\label{eq:dirichlettwist} \sum_{\chi \text{ mod }q}|L(\chi,1/2)|^2 \chi(\ell)   \rightsquigarrow \sum_{\chi \text{ mod }\ell}|L(\chi,1/2)|^2 \chi(-q),\end{align}
for primes $q\neq \ell$. The results were later refined by Young \cite{Young11} and Bettin \cite{Bettin16}. This can be seen as the simplest example of a \emph{spectral reciprocity}-relation. Another example being the $\GL_3\times \GL_2$-relation due to Blomer and Khan \cite{BlomerKhan19} taking the following shape;
$$  \sum_{f\text{ level }q} L(F\otimes f,1/2)L(f,1/2)\lambda_f(p) \rightsquigarrow \sum_{f\text{ level }p} L(F\otimes f,1/2)L(f,1/2)\lambda_f(q), $$
where $f$ runs over an orthonormal basis of Hecke--Maa{\ss} forms of level $q$ (resp.\ $p$), $\lambda_f(n)$ is the $n$-th Hecke eigenvalue of $f$ and $F$ is a (fixed) $\GL_3$-automorphic form.

The left-hand side of~\eqref{eq:dirichlettwist} can be seen as the twisted first moment of a twisted Eisenstein series. The second named author \cite{Nordentoft20} extended this result to general cuspidal holomorphic cusp forms of even weight using a connection to quantum modularity. In this paper we have extended the quantum modularity to general $\GL_2$-forms, and this implies the following reciprocity relation in the cuspidal case.
\begin{thm} \label{recipr}
  Let $\phi$ be a Hecke--Maa{\ss} cuspidal newform of level $q$, weight $k$, nebentypus $\chi_\phi$, sign $\epsilon_\phi$, and Fricke eigenvalue $\eta_\phi$. Assume that either~$s_\phi = k/2$, or~$k\in\{0, 1\}$. Then for any pair of integers $0<c_1<c_2$ with $(c_1,c_2)=(c_1c_2,q)=1$ and a sign $\pm$, we have
  \begin{align}
    \label{mainreci}&\frac{2}{\varphi(c_1)} \sum_{\substack{\chi \modulo c_1\\ \chi(-1)=\pm 1}}\nu(\overline{\phi}, \chi^*,c_1/c(\chi)) L(\overline{\phi}, \chi^*,1/2)\chi(c_2)\\
    \nonumber &  \qquad \qquad \mp\frac{2(-1)^k\eta_\phi }{\varphi(qc_2)} \sum_{\substack{\chi \modulo qc_2\\ \chi(-1)=\pm (-1)^k}}\nu(\phi, \chi^*,qc_2/c(\chi)) L(\phi, \chi^*,1/2)\chi(c_1) = M_{\phi,\pm}+O_{\phi,\eps}((c_1/c_2)^{1-\eps}),
  \end{align}
  where $\nu(\cdot,\cdot,\cdot)$ is a finite Euler product defined as in~\eqref{eq:weightBS} and 
  $$M_{\phi,\pm}=
  \begin{cases}  
    -\eta_\phi \epsilon_\phi L(\phi,1/2), & k=0,\pm=+,\\
    -\eta_\phi\epsilon_\phi \frac{ \sinh (\pi  t_\phi)+i\epsilon_\phi}{\cosh (\pi  t_\phi)} L(\phi,1/2),& k=1, s_\phi\neq 1/2, \pm =+,\\
    -i^k \eta_\phi L(\phi,1/2), & s_\phi=k/2, \pm=+\\
    0,&  \text{else}.
  \end{cases}$$
  Here $L(\phi,s)$ denotes the (standard) $L$-function of $\phi$. 
\end{thm}
\begin{proof}
  By Proposition \ref{BS}, the left hand side of~\eqref{mainreci} is exactly $h_{\overline{W_q}}(-c_2/c_1)$. Now the result follows directly from Proposition \ref{prop:cuspidalhWq}.
\end{proof}
Similarly in non-cuspidal case of Eisenstein series we get the following reciprocity relation for products of Dirichlet $L$-functions extending \cite{Bettin16}.  

\begin{thm} \label{thm:reciprnoncuspidal}
  Let $\chi_i \modulo q_i$ be primitive Dirichlet characters with $\chi_1(-1)\chi_2(-1)=(-1)^k,k\in \{0,1\}$. Then for any pair of integers $0<c_1<c_2$ satisfying $(c_1,c_2)=(c_1c_2,q)=1$ where $q=q_1q_2$, we have
  \begin{align}
    \label{eq:mainreciEisenstein}&\frac{1}{\varphi(c_2)} \sum_{\substack{\chi \modulo c_2}} \nu(\overline{\chi_1}\ast \chi_2, \chi^*,c_2/c(\chi)) L(\chi_2\chi^*,1/2)L(\overline{\chi_1}\chi^*,1/2)\chi(c_1)\\
    \nonumber &  \qquad \qquad -\frac{(-1)^k\eta_{\chi_1,\chi_2} }{\varphi(qc_1)} \sum_{\substack{\chi \modulo qc_1}}\nu(\chi_1\ast \overline{\chi_2}, \chi^*,qc_1/c(\chi)) L(\chi_1\chi^*,1/2)L(\overline{\chi_2}\chi^*,1/2) \chi(-c_2) \\
    \nonumber &= \eta_{\chi_1,\chi_2}\big(A'_\pm (c_2/c_1)^{1/2} + B'_\pm (c_2/c_1)^{1/2}\log(c_2/c_1) + C' \big) +O_{\chi_1,\chi_2,\eps}((c_1/c_2)^{1-\eps}),
  \end{align}
  where $\nu(\cdot,\cdot,\cdot)$ is a finite Euler product defined as in~\eqref{eq:weightBS},  $A'_\pm,B'_\pm$ as in Proposition \ref{prop:hgamma-eisenstein-asymptotic} and
  $$
  C' =\begin{cases}
    - \chi_2(-1) L(\chi_1,1/2)L(\chi_2,1/2),& k=0,\\
    - L(\chi_1,1/2)L(\chi_2,1/2), & k=1,
  \end{cases}$$
  and
  $$\eta_{\chi_1,\chi_2} = (-1)^k \frac{\tau(\overline{\chi_1}) \tau(\chi_2)}{(q_1q_2)^{1/2}}. $$
  Here $L(\chi,s)$ denotes the Dirichlet $L$-function of a Dirichlet character $\chi$. 
\end{thm}
\begin{proof}
  Again by Proposition \ref{BS}, the left hand side of~\eqref{eq:mainreciEisenstein} is exactly $h_{\overline{W_q}}(-c_2/(qc_1))$. Now the result follows directly from Proposition \ref{prop:hgammabar-eisenstein-asymptotic} using that 
  $$ L(E_{\chi_1,\chi_2}^*, \chi,s)=L(\chi_1\chi,s)L(\overline{\chi_2}\chi,s) .$$
\end{proof}
In the special case cuspidal Maa{\ss} forms of level $1$ and where $c_1,c_2$ are prime, we get the simplified version Corollary \ref{cor:recilevel1.1} stated in the introduction using Corollary \ref{cor:BSprimelevel}.

\subsection{Wide moments of automorphic \texorpdfstring{$L$}{L}-functions}

We will now use Proposition \ref{BS} to obtain asymptotic calculations of certain \emph{wide moments} of automorphic $L$-functions. These moments calculations are new and go beyond what has been obtained with the approximate functional equation-approach. We note that these moment calculations are derived using quite surprising input: dynamics of the Gau{\ss} map combined with quantum modularity of additive twists (which we have seen is a very general and non-arithmetic phenomena).

These moment evaluations fit into the a general philosophy of wide moments and distribution of automorphic periods as described in \cite{Nordentoft21} and \cite{Nordentoft21.2} (see also \cite{Be17}). The starting point is that for many natural families of automorphic $L$-functions
$$ \{L(\pi \otimes \chi,1/2): \chi\in G\},$$
the (finite) Fourier transforms 
$$ \widehat{L}(a):= \frac{1}{|G|} \sum_{\chi\in G}L(\pi \otimes \chi,1/2) \chi(a),\quad a\in \widehat{G},$$
are {\lq\lq}well behaved{\rq\rq} in some suitable limit. Here $\pi$ is an automorphic representation of $\GL_n(\mathbb{A}_F)$ and $G$ is some finite group of Hecke characters of some number field $F$. In \cite{Nordentoft21} the setting is that $\pi=\mathrm{BC}_{F/\Q}(\pi_0)$ is the base change to $F$ of a (fixed) cuspidal automorphic representation $\pi_0$ of $\GL_2(\mathbb{A}_\Q)$, $G=\widehat{\mathrm{Cl}_F}$ are class group characters of $F$, where $F$ is an imaginary quadratic field of discriminant tending to infinity. In \cite{Nordentoft21} one has $\pi=\psi$ a Dirichlet character and $G=\{\text{Dirichlet characters modulo $p$}\}$ as $p\rightarrow \infty$.  

In this language we can reinterpret Proposition \ref{BS} as follows; let
$$ G=\{\text{Dirichlet characters modulo $c$}\},\quad \widehat{G}\cong (\Z/c\Z)^\times.$$
Then the Fourier transform of 
$$\chi\text{ mod }c\mapsto \nu(\phi, \chi^*,c/c(\chi)) L(\phi, \chi^*,1/2),$$
is equal to the central values of the additive twist $L$-series
$$ (\Z/c\Z)^\times \ni a\mapsto L(\phi,a/c, 1/2)=L^+(\phi,a/c, 1/2)+L^-(\phi,a/c, 1/2).$$
The well behavedness of the Fourier transform in this context is exactly Proposition \ref{prop:moments-level1}; we can calculate all (even, mixed) moments of the additive twists of level $1$ Maa{\ss} forms! Now it is an easy exercise in Fourier theory that given function $L_i:G\rightarrow \C$ with Fourier transforms $\widehat{L}_i$ for $1\leq i\leq m$, we have
\begin{equation}\label{eq:Fourierfact}
  \sum_{a\in \widehat{G}}\,\, \prod_{i=1}^m \widehat{L}_i(a)= \frac{1}{|G|^{m-1}}\sum_{\substack{\chi_1,\ldots, \chi_m\in G:\\ \chi_1\cdots \chi_m=1}}\,\, \prod_{i=1}^m L_i(\chi_i).   
\end{equation}
The righthand side is what we call the \emph{wide moment} of $L_1,\ldots, L_m$.

As a corollary of Proposition \ref{prop:moments-level1} we obtain the following wide moment calculation.  

\begin{cor} \label{cor:widemoments}
  Let $\phi_1,\ldots, \phi_r$ be distinct Hecke--Maa{\ss} cusp forms of level $1$, $k_1, \dotsc, k_r, \ell_1, \dotsc, \ell_r \in 2 \N_{\geq 0}$ and put $n=\sum_j(k_j+\ell_j)$.  Then we have as $Q\rightarrow \infty$
  \begin{align}
    %\nonumber &\sum_{\substack{0<c\leq Q}} \frac{2^{n-1}}{\varphi(c)^{n-1}} \sideset{}{^+}\sum_{\substack{\chi_{j,k} \modulo c,\\ 1\leq j\leq r, 1\leq k\leq k_j }} 
    %\sideset{}{^-}\sum_{\substack{\psi_{j,\ell} \modulo c,\\1\leq j\leq r, 1\leq \ell\leq \ell_j,\\
    %\Pi_{j,k,\ell}\,\,\chi_{j,k}\psi_{j,\ell}=\mathbf{1}}}\, \prod_{j,k,\ell}  \nu_{j}(\chi_{j,k},\psi_{j,\ell}) L(\phi_j, \chi_{j,k}^*,1/2)L(\phi_j,  \psi_{j,\ell}^*,1/2)\\   
     \nonumber &\sum_{\substack{0<c\leq Q}} \frac{2^{n-1}}{\varphi(c)^{n-1}} \sideset{}{^+}\sum_{\substack{\chi_{j,k} \modulo c,\\ 1\leq j\leq r, 1\leq k\leq k_j }} 
    \sideset{}{^-}\sum_{\substack{\psi_{j,\ell} \modulo c,\\1\leq j\leq r, 1\leq \ell\leq \ell_j:\\
    \Pi_j (\Pi_{k}\,\chi_{j,k})(\Pi_{\ell}\,\psi_{j,\ell})=\mathbf{1}}}\,
    \prod_{j=1}^r  \left(\prod_{k=1}^{k_j} \nu_{j}(\chi_{j,k})L(\phi_j, \chi_{j,k}^*,1/2)\right)\left(\prod_{\ell=1}^{\ell_j} \nu_{j}(\psi_{j,\ell}) L(\phi_j,  \psi_{j,\ell}^*,1/2)\right)\\ 
    \label{eq:widemoment}&= P(\log Q) Q^2+O_{\phi_j, k_j, \ell_j}(Q^{2-\delta}),
  \end{align}
  for some $\delta>0$, where $P$ is a degree $n/2$ polynomial with leading coefficient 
  $$ \prod_{1\leq j \leq r} m_{k_j} (m_{\ell_j} L(\sym^2 \phi_j, 1)^{(k_j+\ell_j)/2}), \qquad m_k := (k-1)!! = \frac{k!}{2^{k/2}(k/2)!}. $$ 
  Here the decorations on the sums mean restricting to characters with $\chi(-1)=\pm1$, $\mathbf{1}$ denotes the principal character (of the relevant modulus suppressed in the notation), and the weights are given by
  $$\nu_{j}(\chi)=\nu(\phi_j, \chi^*, c/c(\chi)),$$
  where $\chi^*\modulo c(\chi)$ denotes the primitive character inducing the Dirichlet character $\chi$, and $\nu(\cdot,\cdot,\cdot)$ is the finite Euler product defined in \eqref{eq:weightBS}.
\end{cor}

\begin{proof}
  This follows directly by combining Proposition \ref{prop:moments-level1}, the Birch--Stevens formula \eqref{birchstevens2} and the Fourier theoretic fact \eqref{eq:Fourierfact}.  
\end{proof}

\appendix

\section{Computations with hypergeometric functions}\label{app:hypergeom}

In this appendix, the notation~$F(a,b;c;z) = {}_2F_1(a,b;c;z)$ stands for the hypergeometric function~\cite[Chapter~9.1]{GradshteynRyzhik2007}.

\subsection{Computation of a determinant}

In this section we provide the proof for a certain identity between hypergeometric functions, which was used in the proof of Proposition~\ref{prop:props-Delta}.

\begin{lemma}\label{lem:hypergeom-1}
  For~$a, b, c\in\C$, $c\not\in\Z$ and~$\abs{z}<1$, we have the equality
  \begin{equation}\label{eq:hypergeom-1}
    \begin{aligned}
      &\frac{(c-a)(c-b)}{c(c-1)} z F(a-c+1,b-c+1; 2-c; z) F(a,b;c+1;z) \\
      &\hspace{12em} + F(a-c,b-c;1-c;z)F(a,b;c;z) = (1-z)^{c-a-b}.
    \end{aligned}
  \end{equation}
\end{lemma}
\begin{proof}
  Assume first that~$\abs{z-\frac12}<\frac12$. Let
  $$ A(z) = F(a-c,b-c;1-c;z), \qquad B(z) = z^c F(a,b;c+1;z), $$
  where~$z^c$ is understood to be the principal value.
  By~\cite[9.103.1-3]{GradshteynRyzhik2007}, the left-hand side can be rewritten as
  $$ \frac{z^{1-c}}c W(z), \qquad \text{where } W(z) = A(z) B'(z) - A'(z) B(z). $$
  Both~$U=A$ and~$U=B$ satisfy the equation
  $$ z(1-z) U'' + (1-c-(a+b+1-2c)z)U' - (a-c)(b-c)U = 0, $$
  see \emph{e.g.} (9.153.1) of~\cite{GradshteynRyzhik2007}, so that~$W(z)$ is a Wronskian determinant.
  We deduce
  $$ \frac{W'}W (z) = -\frac{1-c-(a+b+1-2c)z}{z(1-z)} = \frac{c-1}z + \frac{a+b-c}{1-z}. $$
  Therefore, for some~$C\in\C$ independent of~$z$, we have~$W(z) = C z^{c-1} (1-z)^{c-a-b}$, which gives the equality~\eqref{eq:hypergeom-1} for~$\abs{z-1/2}<1/2$, up to a factor~$C/c$.
  This equality holds for all~$\abs{z}<1$ by analytic continuation.
  We then see that~$C=c$ by letting~$z=0$.
\end{proof}

\begin{lemma}\label{lem:hypergeom-2}
  For~$a, b, c \in \C$, $c\not\in\Z$, we have
  \begin{equation}
    \begin{aligned}
      & (c-a)(c-b) F(a,b;a+b-c+1;\tfrac12) F(a,b;c+1;\tfrac12) \\
      &\quad{} + c(a+b-c)F(a,b;a+b-c;\tfrac12)F(a,b;c;\tfrac12) \\
      & \hspace{7em} =  \frac{2^{a+b}}{\Gamma(a)\Gamma(b)} \Gamma(a+b-c+1)\Gamma(c+1).
    \end{aligned}\label{eq:hypergeom-2}
  \end{equation}
\end{lemma}

\begin{proof}
  By~\cite[9.131.2]{GradshteynRyzhik2007} with~$\gamma\gets a+b-c$ and~$z\gets 1/2$, we have
  $$ F(a,b;a+b-c;\tfrac12) = \frac{\Gamma(a+b-c)\Gamma(-c)}{\Gamma(b-c)\Gamma(a-c)} F(a,b;c+1;\tfrac12) + 2^c \frac{\Gamma(a+b-c)\Gamma(c)}{\Gamma(a)\Gamma(b)} F(a-c,b-c;1-c;\tfrac12), $$
  and the similar identity with~$c$ replaced by~$c-1$. We insert this in the left-hand side of~\eqref{eq:hypergeom-2}. Since
  $$ (c-a)(c-b)\frac{\Gamma(a+b+1-c)\Gamma(1-c)}{\Gamma(b+1-c)\Gamma(a+1-c)} + c(a+b-c) \frac{\Gamma(a+b-c)\Gamma(-c)}{\Gamma(b-c)\Gamma(a-c)} = 0, $$
  the terms involving~$F(\cdots; c; \tfrac12) F(\cdots; c+1; \tfrac12)$ cancel out, and so the left-hand side of~\eqref{eq:hypergeom-2} is equal to
  \begin{align*}
    & (c-a)(c-b) 2^{c-1} \frac{\Gamma(a+b+1-c)\Gamma(c-1)}{\Gamma(a)\Gamma(b)} F(a+1-c,b+1-c;2-c;\tfrac12) F(a,b;c+1;\tfrac12) \\
    &{} + c(a+b-c)2^c \frac{\Gamma(a+b-c)\Gamma(c)}{\Gamma(a)\Gamma(b)} F(a-c,b-c;1-c;\tfrac12) F(a,b;c;\tfrac12) \\
    ={}& 2^c \frac{\Gamma(a+b+1-c)\Gamma(c+1)}{\Gamma(a)\Gamma(b)} \Big(\frac{(c-a)(c-b)}{c(c-1)} \tfrac12 F(a+1-c,b+1-c;2-c;\tfrac12) F(a,b;c+1;\tfrac12)\\
    & \hspace{16em} - F(a-c,b-c;1-c;\tfrac12) F(a,b;c;\tfrac12)\Big).
  \end{align*}
  By Lemma~\ref{lem:hypergeom-1} at~$z=1/2$, the quantity inside the parentheses is equal to~$2^{a+b-c}$.
\end{proof}

\subsection{Computations relative to the functional equation}

In this section, we prove an identity for a quotient of hypergeometric functions, which was used in the proof of Lemma~\ref{lem:EF-Fk}.
First define
$$ G_\pm(a, b, c) := \frac{\Gamma(c)}{\Gamma(1-c+a+b)} F(a,b;1-c+a+b;\tfrac12) \pm \frac{\Gamma(c-a)}{\Gamma(1-c+b)} F(a,b;c;\tfrac12). $$
\begin{lemma}\label{lem:hypergeom-FE}
  Suppose~$a, b\not\in\Z_{\geq 1}$ and~$c\not\in\Z$. Then whenever~$a+b-2c \in \Z$, the quotient
  $$ Q_\pm(a,b,c) := \frac{G_\pm(a,b,c)}{G_\pm(1-a,1-b,1+c-a-b)} $$
  can be expressed in terms of elementary functions and~$\Gamma$ functions. More precisely, letting~$n = a+b-2c$, we have 
  $$ Q_\pm(a, b, c) = (-1)^n \frac{2^{a+b-1} \Gamma(1-a)\Gamma(1-b)}{\Gamma(\frac{2+n-a-b}2) \Gamma(\frac{2-n-a-b}2)} \Big(-1 \pm \frac{\sin(\frac\pi2(n+a-b))}{\sin(\frac\pi2(n+a+b))}\Big). $$
\end{lemma}
\begin{proof}
  We have
  $$ G_\pm(1-a, 1-b, 1+c-a-b) = \frac{\Gamma(1+c-a-b)}{\Gamma(2-c)} F(1-a, 1-b ; 2-c ; \tfrac12) \pm \frac{\Gamma(c-b)}{\Gamma(1-c+a)} F(1-a, 1-b ; 1+c-a-b; \tfrac12). $$
  Thus, letting
  \begin{align*}
    \omega_1 := {}& F(a,b;c;\tfrac12), &
    \omega_2 := {}& F(1-a,1-b;2-c;\tfrac12), \\
    \omega_3 := {}& F(a,b;1-c+a+b;\tfrac12), &
    \omega_4 := {}& F(1-a,1-b;1+c-a-b;\tfrac12),
  \end{align*}
  we find
  \begin{align*}
    G_\pm(a, b, c) ={}& \frac{\Gamma(c)}{\Gamma(1-c+a+b)} \omega_3 \pm \frac{\Gamma(c-a)}{\Gamma(1-c+b)} \omega_1, \\
    G_\pm(1-a, 1-b, 1+c-a-b) ={}& \frac{\Gamma(1+c-a-b)}{\Gamma(2-c)} \omega_2 \pm \frac{\Gamma(c-b)}{\Gamma(1-c+a)} \omega_4.
  \end{align*}
  By~\cite[(9.131.2) and (9.131.1)]{GradshteynRyzhik2007}, the functions~$\omega_j$ are related by
  \begin{align*}
    \omega_3 ={}& \frac{\Gamma(1-c+a+b)\Gamma(1-c)}{\Gamma(1-c+a)\Gamma(1-c+b)} \omega_1 + 2^{a+b-1} \frac{\Gamma(1-c+a+b)\Gamma(c-1)}{\Gamma(a)\Gamma(b)} \omega_2, \\    
    \omega_4 ={}& 2^{1-a-b}\frac{\Gamma(1+c-a-b)\Gamma(1-c)}{\Gamma(1-a)\Gamma(1-b)} \omega_1 + \frac{\Gamma(1+c-a-b)\Gamma(c-1)}{\Gamma(c-a)\Gamma(c-b)} \omega_2,
  \end{align*}
  and therefore
  \begin{align*}
    G_\pm(a, b, c) ={}& \lambda_1 \omega_1 + \lambda_2 \omega_2, &
    G_\pm(1-a, 1-b, 1+c-a-b) = {}& \mu_1\omega_1 + \mu_2 \omega_2,
  \end{align*}
  where
  \begin{align*}
    \lambda_1 ={}& \frac{\Gamma(c)\Gamma(1-c)}{\Gamma(1-c+a)\Gamma(1-c+b)} \pm \frac{\Gamma(c-a)}{\Gamma(1-c+b)}, &
    \lambda_2 ={}& 2^{a+b-1} \frac{\Gamma(c)\Gamma(c-1)}{\Gamma(a)\Gamma(b)}, \\
    \mu_1 ={}& \pm 2^{1-a-b} \frac{\Gamma(c-b)\Gamma(1+c-a-b)\Gamma(1-c)} {\Gamma(1-c+a)\Gamma(1-a)\Gamma(1-b)}, &
    \mu_2 ={}& \frac{\Gamma(1+c-a-b)}{\Gamma(2-c)} \pm \frac{\Gamma(c-1)\Gamma(1+c-a-b)}{\Gamma(1-c+a)\Gamma(c-a)}.
  \end{align*}
  By a straightforward computation using the complement formula, we obtain
  \begin{align*}
    \lambda_1 \mu_2 - \lambda_2 \mu_1 ={}& \pm \frac{\Gamma(1+c-a-b)\Gamma(c-a)}{\Gamma(2-c)\Gamma(1-c+b)} \Big(1 - \frac{\sin(\pi(c-a))^2}{\sin(\pi c)^2} + \frac{\sin(\pi a)\sin(\pi b)\sin(\pi(c-a))}{\sin(\pi c)^2 \sin(\pi(c-b))}\Big) \\
    ={}& \pm \frac{\Gamma(1+c-a-b)\Gamma(c-a)}{\Gamma(2-c)\Gamma(1-c+b)} \frac{\sin(\pi a)\sin(\pi(2c-a-b))}{\sin(\pi c)\sin(\pi(c-b))}.
  \end{align*}
  This is zero indeed since~$2c-a-b \in \Z$. We deduce that~$Q_\pm(a, b, c) = \lambda_1 / \mu_1$, which is equal, by the complement formula, to
  \begin{align*}
    \lambda_1 \mu_1^{-1} ={}& 2^{a+b-1} \frac{\Gamma(1-a)\Gamma(1-b)}{\Gamma(1-c)\Gamma(1+c-a-b)} \Big(\frac{\sin(\pi(c-b))}{\sin(\pi(c-a))} \pm \frac{\sin(\pi(c-b))}{\sin(\pi c)}\Big).
  \end{align*}
  Expressing this in terms of~$n=a+b-2c$ yields our formula as claimed.
\end{proof}

\providecommand{\bysame}{\leavevmode\hbox to3em{\hrulefill}\thinspace}
\providecommand{\MR}{\relax\ifhmode\unskip\space\fi MR }
\providecommand{\MRhref}[2]{%
  \href{http://www.ams.org/mathscinet-getitem?mr=#1}{#2}
}
\providecommand{\href}[2]{#2}

\end{document}